\documentclass[12pt,a4paper]{article}

\usepackage{amsmath}    
\usepackage{amssymb}
\usepackage{graphicx}   
\usepackage{verbatim}   
\usepackage{color}      

\usepackage{hyperref}   

\usepackage{enumerate}
\usepackage{mathrsfs}
\usepackage{empheq}
\usepackage{bbold}
\usepackage[british]{babel}
\usepackage{lscape}
\usepackage{cleveref}
\crefformat{footnote}{#2\footnotemark[#1]#3}
\usepackage[font={small,it}]{caption}

\usepackage[final]{showlabels}

\usepackage{amsthm}
\newtheorem{definition}{Definition}[section]
\newtheorem{lemma}{Lemma}[section]
\newtheorem{corollary}{Corollary}[section]
\newtheorem{theorem}{Theorem}[section]

\newtheorem{proposition}{Proposition}[section]

\newcommand{\del}{\partial}
\renewcommand{\theta}{\vartheta}
\renewcommand{\phi}{\varphi}
\newcommand{\vecc}[2]{\left ( \begin{array}{c}#1\\#2\\ \end{array}\right )}
\newcommand{\veccc}[3]{\left ( \begin{array}{c}#1\\#2\\#3\\ \end{array}\right )}

\newcommand{\dd}{\mathrm{d}}

\newcommand{\id}{\mathbb{1}}

\renewcommand{\div}{\mathrm{div\,}}
\renewcommand{\vec}{\mathbf}

\newcommand{\ii}{\mathbb{i}}

\renewcommand{\title}{Analysis of the multi-dimensional semi-discrete Active Flux method using the Fourier transform}

\newcommand{\authorOne}{Wasilij Barsukow\footnote{Bordeaux Institute of Mathematics, Bordeaux University and CNRS/UMR5251, Talence, 33405 France, wasilij.barsukow@math.u-bordeaux.fr}}
\newcommand{\authorTwo}{Janina~Kern\footnote{\label{fn:wue}Institute for Mathematics, University of Wurzburg, Emil-Fischer-Strasse 40, 97074 Wurzburg, Germany}} 
\newcommand{\authorThree}{Christian~Klingenberg\cref{fn:wue}}
\newcommand{\authorFour}{Lisa~Lechner\cref{fn:wue}}

\begin{document}

\begin{center} \Large
\title

\vspace{1cm}

\date{}
\normalsize

\authorOne, \authorTwo, \authorThree, \authorFour
\end{center}

\begin{abstract}

The degrees of freedom of Active Flux are cell averages and point values along the cell boundaries. These latter are shared between neighbouring cells, which gives rise to a globally continuous reconstruction. The semi-discrete Active Flux method uses its degrees of freedom to obtain Finite Difference approxi\-mations to the spatial derivatives which are used in the point value update. The averages are updated using a quadrature of the flux and making use of the point values as quadrature points. The integration in time employs standard Runge-Kutta methods. We show that this generalization of the Active Flux method in two and three spatial dimensions is stationarity preserving for linear acoustics on Cartesian grids, and present an analysis of numerical diffusion and stability.



\end{abstract}

\section{Introduction}

The classical Active Flux method has been introduced in \cite{eymann13}, based on a one-dimensional method from \cite{vanleer77}. From the beginning, it was conjectured that the continuous reconstruction, i.e. the absence of Riemann problems might help alleviate difficulties that traditional Finite Volume methods are facing in multiple dimensions. For instance, they are usually not preserving discrete involutions, they are not low Mach number compliant for the Euler equations, and not stationarity preserving. This situation does not improve even if the full multi-dimensional Riemann problem is solved and used in a Godunov method (\cite{barsukow17}). 

As Finite Volume methods, Active Flux evolves cell averages as discrete degrees of freedom, and additionally evolves point values. These latter are located at cell boundaries and are a way to ensure global continuity of the reconstruction, or -- in a Finite Element sense -- of the numerical solution. The first multi-dimensional system that Active Flux was applied to (in \cite{eymann13}) was linear acoustics. In multiple spatial dimensions, it is an interesting system of equations because it cannot be reduced to (multi-dimensional) advection. While the update of the average is immediately possible via quadrature along the cell boundary (using the point values as quadrature points), the update of the point values was initially achieved (e.g. in \cite{barsukow18activeflux}) using an exact evolution operator. This operator was derived for linear acoustics in \cite{barsukow17}. The reconstruction was serving as its initial datum. The structure preserving properties of this method have already been analyzed in \cite{barsukow18activeflux}. It has been found that indeed Active Flux is stationarity preserving. 

For nonlinear problems it is more difficult to obtain evolution operators of sufficient order of accuracy, let alone exact ones. This led \cite{abgrall20,abgrall22} to consider semi-discrete Active Flux methods, where the same degrees of freedom are used in order to discretize the spatial derivatives, while integration in time follows the method-of-lines strategy. As has been outlined in \cite{abgrall22proceeding}, this typically leads to reduced CFL conditions. The advantage of the semi-discrete approach is its immediate applicability to various kinds of systems of conservation laws. However, one would like not to lose structure preservation. 

This paper presents an analysis of the semi-discrete Active Flux method on linear problems and a comparison between the classical and the semi-discrete approaches, with an emphasis on structure preservation. As is shown below, the structure preserving properties of the semi-discrete approach are very similar to those of the classical one. In the context of the low Mach number limit for the Euler equations this has already been observed experimentally in \cite{barsukow24afeuler}.

The paper is organized as follows: Section \ref{sec:equations} introduces the equations and the analytical stationary states and Section \ref{sec:semidiscreteAF} presents the semi-discrete Active Flux method on two-dimensional and three-dimensional Cartesian grids. The discrete Fourier transform is introduced in Section \ref{sec:fourier} and is used to analyze stationarity preservation in Section \ref{sec:fourierAF}. Section \ref{sec:numdiff} presents an analysis of numerical diffusion and stability. Some numerical examples follow in Section \ref{sec:numerical}.

We denote by $P^k$ univariate polynomials of degree at most $k$, and by $P^{k,k}$ bivariate polynomials with degree at most $k$ in each variable. We also occasionally denote by $\mathbb R[x]$ or $\mathbb C[x, y]$ polynomials in $x$ with real coefficients, and polynomials in $x$, $y$ with complex coefficients, respectively. $d$ denotes the number of spatial dimensions, and objects with $d$ components are typeset in boldface. $\mathscr M^{a \times b}(\mathbb C)$ denotes matrices with complex entries of $a$ rows and $b$ columns. Indices never denote differentiation. $\id_m$ denotes the identity map on $\mathbb R^m$.

\section{Acoustic equations and their stationary states} \label{sec:equations}

\subsection{General linear systems}

Stationary states are solutions to evolutionary partial differential equations (PDEs) that remain constant in time. To achieve stability, numerical methods add numerical diffusion. A setup that is stationary according to the PDE might keep being diffused away and no longer be stationary in the discrete setting. Thus, the stationary states of the discretization are often a small, not representative subset of the stationary states of the PDE. It is preferable for numerical schemes to possess numerical stationary states that discretize \emph{all} the analytic stationary states of the underlying conservation law; a definition of what this really means is given in Section \ref{ssec:statpres}. In this section, the stationary states of the linear acoustic equations shall be analyzed. This is easiest done upon applying the Fourier transform to them.

Consider a $m \times m$ hyperbolic system of linear PDEs
\begin{align}
\partial_tq+\mathbf{J}\cdot\nabla q=0\qquad &q:\mathbb{R}_0^+\times\mathbb{R}^d\rightarrow\mathbb{R}^m \label{eq:conservationLaw}\\
&\mathbf{J}=(J_{1}, ..., J_{d}), J_{i}\in\mathscr M^{m\times m}(\mathbb{R}). \nonumber
\end{align}
Occasionally we will use notation $J_x \equiv J_1$, $J_y \equiv J_2$, $J_z \equiv J_3$ instead.
One Fourier mode of a function $q: \mathbb{R}_0^+\times\mathbb{R}^d\rightarrow \mathbb{R}^m$ is of the form
\begin{equation}
q(t, \mathbf{x})=\hat{q}(t, \vec k)\exp(\mathbb{i}\mathbf{k}\cdot\mathbf{x}).
\end{equation} 
Here, $\vec k \in \mathbb R^d$ is called the wave vector and determines the spatial frequency of the Fourier mode, while its amplitude $\hat{q}$ can be chosen differently for each $\vec k$. General solutions are obtained as linear combinations of such modes for different values of $\vec k$.
Inserting the mode into \eqref{eq:conservationLaw} yields 
\begin{align}
\frac{\dd}{\dd t} \hat q + \ii \mathbf{J}\cdot\vec k \hat q=0. \label{eq:Jk}
\end{align}
Hyperbolicity guarantees that the matrix $\mathbf{J}\cdot\vec k = J_{1} k_1 + \dots + J_{d} k_{d} $ is diagonalizable. The mode $\hat q \exp(\ii \vec k \cdot \vec x)$ is stationary if $\hat q$ is in the nullspace of $\vec J \cdot \vec k$. This can be achieved through particular choices of $\vec k$, trivially for $\vec k = 0$, which (since every Fourier mode depends on $\vec x$ as $\exp(\ii \vec k \cdot \vec x)$) corresponds to the data being uniformly constant. Non-trivial stationary states (\cite{barsukow17a}) are those for which no restriction on $\vec k$ is necessary, i.e. where for any $\vec k$ there exists a $\hat q_\text{stat}(\vec k) \in \mathbb C^m \backslash \{ 0 \}$ such that $(\vec J \cdot \vec k) \hat q_\text{stat}(\vec k) = 0$. A necessary condition for the existence of non-trivial stationary states is $\det (\vec J \cdot \vec k) = 0 \quad \forall \vec k \in \mathbb R^d$. Observe that the amplitude $\hat q_\text{stat}$ of a stationary mode generally depends on $\vec k$, see below for some examples.

\subsection{Acoustic equations}

The acoustic equations in $d$ spatial dimensions are given as
\begin{subequations}
\begin{align}
\partial_t\mathbf{v}+\nabla p=0\\
\partial_t p + \nabla\cdot\mathbf{v}=0
\end{align}\label{acousticEquationsOhneMatrizen}
\end{subequations}
with velocity $\mathbf{v} \colon \mathbb R^+_0 \times \mathbb R^d \to \mathbb{R}^d$ and pressure $p \colon \mathbb R^+_0 \times \mathbb R^d \to \mathbb{R}$.

In $d$ spatial dimensions, the matrix $\vec J \cdot \vec k$ reads
\begin{align}
 \left( \begin{array}{c|c} 0_{d \times d} & \vec k \\\hline \vec k^\text{T} & 0 \end{array} \right ) \in \mathscr M^{(d+1) \times (d+1)}
\end{align}
and any element $(\vec U, P)^\text{T} \in \mathbb C^{d+1}$ of its nullspace has to fulfill
\begin{align}
\vec k P &= 0 \\
\vec k^\text{T} \vec U &= 0
\end{align}
For any $\vec k$, the following is a nullspace of $\vec J \cdot \vec k$:
\begin{align}
 N_\text{non-trivial}^{(d)} := \{ (\vec U, 0)^\text{T} : \vec U \perp \vec k \}.
\end{align}
It is $d-1$ dimensional. In 3-d, one can choose
\begin{align}
 N_\text{non-trivial}^{(3)} = \mathrm{span}\Big\{(-k_y, k_x, 0, 0)^\text{T}, (-k_z, 0, k_x, 0)^\text{T}\Big\} \label{eq:3dker}
\end{align}
and in 2-d
\begin{align}
 N_\text{non-trivial}^{(2)} = \mathrm{span}\Big\{(-k_y, k_x, 0)^\text{T}\Big\}. \label{eq:2dker}
\end{align}

For special values of $\vec k$ there are additional elements in the kernel, which are referred to as \emph{trivial stationary states}. They are usually represented well by any kind of numerical method and of little interest in the following. On the contrary, the \emph{non-trivial stationary states} in $N_\text{non-trivial}^{(d)}$ are usually poorly represented by numerical methods (see \cite{barsukow17a} for more details). The following will show that they are well represented by the semi-discrete Active Flux method.

\section{Semi-discrete Active Flux on 2-d and 3-d Cartesian Grids} \label{sec:semidiscreteAF}

\subsection{Degrees of Freedom}\label{sec:degreesOfFreedom}
In contrast to classical Finite Volume methods that only involve the cell average as degree of freedom, the Active Flux method additionally uses point values distributed on the cell boundaries as degrees of freedom. 

\subsubsection{General remarks}

In two spatial dimensions the computational grid consists of cells 
\begin{equation}
\mathcal{C}_{ij}=\left[\left(i-\frac{1}{2}\right)\Delta x, \left(i+\frac{1}{2}\right)\Delta x\right]\times \left[\left(j-\frac{1}{2}\right)\Delta y, \left(j+\frac{1}{2}\right)\Delta y\right] \subset \mathbb{R}^2
\end{equation}
and in 3-d
\begin{align}
\mathcal{C}_{ijk}=\left[\left(i-\frac{1}{2}\right)\Delta x, \left(i+\frac{1}{2}\right)\Delta x\right]&\times \left[\left(j-\frac{1}{2}\right)\Delta y, \left(j+\frac{1}{2}\right)\Delta y\right]\\\nonumber&\times \left[\left(k-\frac{1}{2}\right)\Delta z, \left(k+\frac{1}{2}\right)\Delta z\right] \subset \mathbb{R}^3.
\end{align}
We denote by $\vec x_{ij}$ and $\vec x_{ijk}$ the cell centroid.

There is one average in every cell and a certain number of point values located at the boundary of the cell. These latter are shared. We make a distinction between degrees of freedom that \emph{belong} to a cell, whose number we will denote by $N^\text{dof}$, and those that are \emph{accessible} to a cell:

\begin{definition}\label{def:belong}
 The degrees of freedom \emph{accessible} to a cell $\mathcal C_{ij}$ (or $\mathcal C_{ijk}$) are those located in it or along its boundary, or that are averages over it. Their number per cell is denoted by $N^\text{dof}_\text{acc}$ and they are denoted by 
 \begin{align}
  &q_{0,ij}, q_{1,ij}, \ldots, q_{N^\text{dof}_\text{acc}-1, ij} \text{ and}\\
  &q_{0,ijk}, q_{1,ijk}, \ldots, q_{N^\text{dof}_\text{acc}-1, ijk}
 \end{align}
 in 2-d and 3-d, respectively, with a numbering that is arbitrary but fixed once for all cells, and we reserve the index 0 for the cell average. Locations of the point values are denoted by $\vec x_{r,ij}$ and $x_{r,ijk}$, $r \in 1, \dots, N^\text{dof}_\text{acc}-1$ and we define the relative coordinate
 \begin{align}
  \vec x_r := \vec x_{r,ij} - \vec x_{ij} \in \left[-\frac{\Delta x}{2}, \frac{\Delta x}{2} \right] \times \left[-\frac{\Delta y}{2}, \frac{\Delta y}{2} \right]
 \end{align}
 or
 \begin{align}
  \vec x_r := \vec x_{r,ijk} - \vec x_{ijk} \in \left[-\frac{\Delta x}{2}, \frac{\Delta x}{2} \right] \times \left[-\frac{\Delta y}{2}, \frac{\Delta y}{2} \right] \times \left[-\frac{\Delta z}{2}, \frac{\Delta z}{2} \right]
 \end{align}
 in two and three spatial dimensions, respectively.
 \end{definition}
 
 Examples are given below.
 
 \begin{definition}
 The degrees of freedom \emph{belonging} to a cell are a minimal set such that all the accessible degrees of freedom can be obtained from them through shifts by $\Delta x$, $\Delta y$ (and $\Delta z$) in the two/three directions.  We agree in the following on choosing the same ones in each cell. Their number per cell is denoted by $N^\text{dof}$ and we denote them by $q_{ij}^X$ where $X$ takes values in a set of identifiers defined below. Locations of degrees of freedom that are point values are denoted by $\vec x_{ij}^X$ and we define the relative coordinate
 \begin{align}
  \vec x^X := \vec x_{ij}^X - \vec x_{ij}.
 \end{align}
 \end{definition}
 
In other words, consider an equivalence relation $\sim$ between any two degrees of freedom on an infinite grid, with two degrees of freedom $p_1, p_2$ being equivalent if there exist $A_x,A_y,A_z \in \mathbb Z$ such that the shift 
 \begin{align}
 (x, y, z) \mapsto (x + A_x \Delta x, y + A_y \Delta y, z + A_z \Delta z) \label{eq:intshift} 
 \end{align}
 maps $p_1$ onto $p_2$. For example, all point values at nodes and all the cell averages are equivalent. The degrees of freedom belonging to one cell is the quotient $\mathbb{DOF}/\sim$ of all the degrees of freedom $\mathbb{DOF}$ on the grid by $\sim$. Finally, instead of the equivalence classes we speak of representative elements chosen according to some conventions for definiteness (explained above and shown on Figure \ref{fig:distributionpointvalues}). Degrees of freedom that belong to a cell become important in the context of the Fourier analysis of Section \ref{sec:fourier}, since discrete Fourier modes are the eigenbasis of the operator defined in \eqref{eq:intshift}.

Consider for every accessible degree of freedom $r = 0, \dots, q_{N^\text{dof}_\text{acc}-1}$ the number of cells that share it and denote this number by $\alpha_{r} \in \mathbb N$. Then
\begin{align}
 N^\text{dof} = \sum_{r = 0}^{N^\text{dof}_\text{acc}-1} \frac{1}{\alpha_r}.
\end{align}

The distinction between degrees of freedom belonging to a cell and those accessible to it is due to global continuity and the fact that degrees of freedom are shared. In DG methods, for example, the accessible degrees of freedom are just the ones that belong to the cell.

Definition \ref{def:belong} states that each of the degrees of freedom that belong to a cell forms a lattice with spacing $\Delta x, \Delta y$ or $\Delta x, \Delta y$ and $\Delta z$, respectively. This will become important below in the context of the Fourier transform. 

\newcommand{\nodeonly}{\mathrm{N}}
\newcommand{\evonly}{\mathrm{E}_\mathrm{V}}
\newcommand{\ehonly}{\mathrm{E}_\mathrm{H}}
\newcommand{\avgonly}{\mathrm{A}}

\newcommand{\exonly}{\mathrm{E}_x}
\newcommand{\eyonly}{\mathrm{E}_y}
\newcommand{\ezonly}{\mathrm{E}_z}

\newcommand{\fxonly}{\mathrm{F}_x}
\newcommand{\fyonly}{\mathrm{F}_y}
\newcommand{\fzonly}{\mathrm{F}_z}

\newcommand{\node}{^{\nodeonly}}
\newcommand{\ev}{^{\evonly}}
\newcommand{\eh}{^{\ehonly}}
\newcommand{\avg}{^{\avgonly}}

\newcommand{\ex}{^{\exonly}}
\newcommand{\ey}{^{\eyonly}}
\newcommand{\ez}{^{\ezonly}}

\newcommand{\fx}{^{\fxonly}}
\newcommand{\fy}{^{\fyonly}}
\newcommand{\fz}{^{\fzonly}}

\subsubsection{Degrees of freedom in two spatial dimensions}

\begin{figure}
 \centering
 \includegraphics[width=0.4\textwidth]{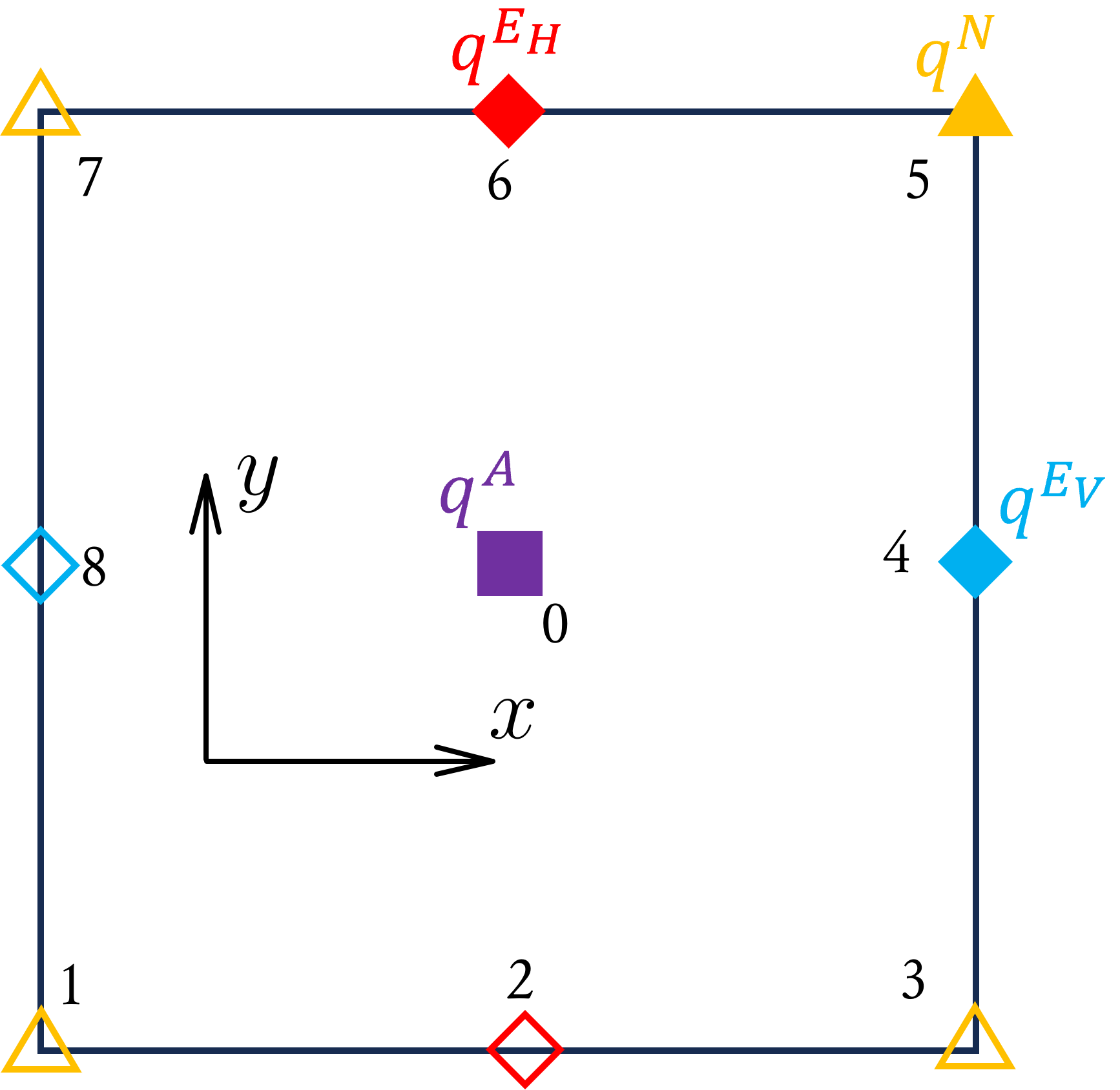} \qquad\qquad \includegraphics[width=0.4\textwidth]{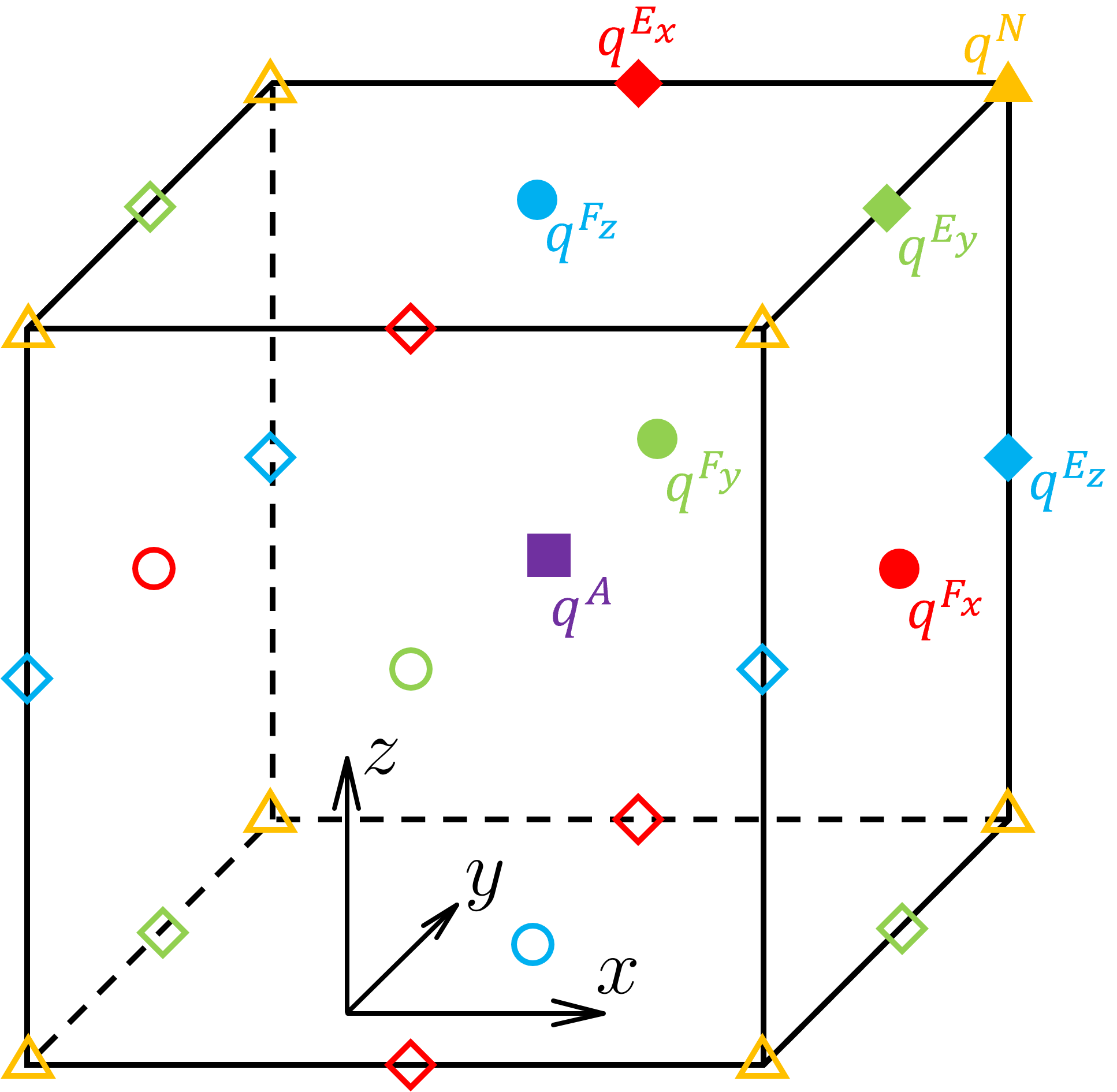}
 \caption{Distribution of the point values along the cell boundary for a cell of a two-dimensional (left) and three-dimensional (right) Cartesian grid. The cell average is depicted with a square located at the center of the cell.}
 \label{fig:distributionpointvalues}
\end{figure}

In two spatial dimensions (see Figure \ref{fig:distributionpointvalues}) we consider a classical distribution of in total 8 point values located at the corners and the midpoints of the cell edges. The $N^\text{dof}_\text{acc} = 9$ degrees of freedom accessible to cell $\mathcal C_{ij}$, denoted by $q_{r,ij}$, $r = 0, \dots, 8$, are
\begin{align}
 &q_{7,ij} := q\node_{i-1,j}     &&q_{6,ij} := q\eh_{ij}     &&q_{5,ij} := q\node_{ij} \label{eq:2ddofs1}\\
 &q_{8,ij} := q\ev_{i-1,j} &&q_{0,ij} := q\avg_{ij}         && q_{4,ij} := q\ev_{ij}\\
 &q_{1,ij} := q\node_{i-1,j-1}   &&q_{2,ij} := q\eh_{i,j-1}  &&q_{3,ij} := q\node_{i,j-1} \label{eq:2ddofs3}
\end{align}

Recall that the numbering is arbitrary but needs to be fixed once.

Each corner of a cell is likewise a corner to three other adjacent cells and each midpoint of a cell edge is shared with one other cell, i.e. $N^\text{dof} = 4 \cdot \frac14 + 2 \cdot 2 \cdot \frac12 + 1 = 4$. Therefore, 4 degrees of freedom belong to a cell $\mathcal C_{ij}$ and we have chosen them in the top right corner as
\begin{enumerate}[--]
\item cell average $q\avg_{ij} = q_{0,ij}$,
\item node value $q\node_{ij} = q_{5,ij}$,
\item horizontal edge value $q\eh_{ij} = q_{6,ij}$,
\item vertical edge value $q\ev_{ij} = q_{4,ij}$.
\end{enumerate}
as shown in Figure \ref{fig:distributionpointvalues}.

\subsubsection{Degrees of freedom in three spatial dimensions}

The natural extension of the previous two-dimensional situation to three spatial dimensions results in 26 point values accessible to each cell, i.e. $N^\text{dof}_\text{acc} = 27$. Eight of them are located at the corners of a cell, 12 at the midpoints of the edges and six at the midpoints of the faces of a cell, see Figure \ref{fig:distributionpointvalues}. Again, as in the 2-d case, several cells are sharing these point values, so that one is left with $N^\text{dof} = 8$ degrees of freedom that belong to cell $\mathcal{C}_{ijk}$:
\begin{enumerate}[--]
\item one cell average $q\avg_{ijk} = q_{0,ijk}$,
\item one node value $q\node_{ijk}$,
\item three edge values $q\ex_{ijk}$, $q\ey_{ijk}$, $q\ez_{ijk}$, one per edge parallel to the three axes,  
\item three face values $q\fx_{ijk}$, $q\fy_{ijk}$, $q\fz_{ijk}$, one per face orthogonal to the three axes.
\end{enumerate}

Observe that the values on edges and faces are point values at the respective centroids, and not averages over these entities.

\subsection{Update of the Average}

The update of the average will first be described for a generic nonlinear conservation law
\begin{align}
 \del_t q + \nabla \cdot \vec f(q) = 0, \label{eq:conslaw}
\end{align}
since generalization improves clarity in this case. At the end of the Section the special case of linear systems is explained.

For updating the cell average, the conservation law \eqref{eq:conslaw}
is integrated over one cell ($\mathcal C_{ij}$ or $\mathcal C_{ijk}$), and Gauss' law is applied. This gives
\begin{align}
\frac{\dd}{\dd t} q\avg_{ij} + \frac{1}{\Delta x \Delta y} \sum_{e \subset \del \mathcal{C}_{ij}} \int_{e} \vec n_e \cdot \vec f(q) \,\dd s &= 0 &&\text{(2-d)} \label{eq:updateaverage2d}\\
\frac{\dd}{\dd t} q\avg_{ijk} + \frac{1}{\Delta x \Delta y\Delta z} \sum_{f \subset \del \mathcal{C}_{ijk}} \int_{f} \vec n_f \cdot \vec f(q)\, \dd s &= 0 &&\text{(3-d)}
\end{align}
with 
\begin{align}
q\avg_{ij}(t) = \frac{1}{\Delta x \Delta y} \int_{\mathcal{C}_{ij}} q(t, x, y)\, \dd x\,\dd y    &&\text{(2-d)}\\
q\avg_{ijk}(t) = \frac{1}{\Delta x \Delta y \Delta z} \int_{\mathcal{C}_{ijk}} q(t, x, y, z) \,\dd x\,\dd y\, \dd z &&\text{(3-d)}
\end{align}
In Active Flux, the locations of the point values are chosen such that they can be used as quadrature points for the flux averages through the edges $e$/faces $f$. The distribution of point values as described in the previous Section allows to use Simpsons' rule, such that
\begin{subequations}
\begin{align}
\frac{\dd}{\dd t} q\avg_{ij} &+ \frac{\hat f^x_{i+\frac12,j} - \hat f^x_{i-\frac12,j}}{\Delta x} + \frac{\hat f^y_{i,j+\frac12} - \hat f^y_{i,j-\frac12}}{\Delta y} = 0 &&\text{(2-d)} \\
\frac{\dd}{\dd t} q\avg_{ijk} &+ \frac{\hat f^x_{i+\frac12,j,k} - \hat f^x_{i-\frac12,j,k}}{\Delta x} + \frac{\hat f^y_{i,j+\frac12,k} - \hat f^y_{i,j-\frac12,k}}{\Delta y} \\
\nonumber &\phantom{mmmmmmmmmmmmm}+ \frac{\hat f^z_{i,j,k+\frac12} - \hat f^z_{i,j,k-\frac12}}{\Delta z} = 0 &&\text{(3-d)}
\end{align}\label{eq:updateAverages}
\end{subequations}
with the numerical fluxes 
\begin{align}
 \hat f^x_{i+\frac12,j} &= \frac{1}{6} \left( f^x(q\node_{ij})+ 4f^x(q\ev_{ij} )  + f^x(q\node_{i,j-1})  \right )\\
 \hat f^y_{i,j+\frac12} &= \frac{1}{6} \left( f^y(q\node_{i-1,j})+ 4 f^y(q\eh_{ij})   + f^y(q\node_{ij} ) \right )\\
 \hat f^x_{i+\frac12,j,k} &= \frac1{36} \left(  f^x(q\node_{ijk})+ 4 f^x(q\ez_{ijk}) + f^x(q\node_{i,j,k-1}) \right . \\ &\phantom{mm}
  \nonumber +  \left. 4f^x(q\ey_{ijk})+ 16 f^x(q\fx_{ijk}) +4f^x(q\ey_{i,j,k-1}) \right . \\ &\phantom{mm} 
  \nonumber +  \left.f^x(q\node_{i,j-1,k})+ 4 f^x(q\ez_{i,j-1,k}) + f^x(q\node_{i,j-1,k-1} )\right ) \qquad \text{etc.}
\end{align}

Observe that the numerical fluxes do not come from a Riemann solver, but are simply quadratures of the physical flux. Observe also that the only approximation made here is the replacement of integrals by quadratures. While these formulas are valid for nonlinear conservation laws, below they are used to discretize \eqref{eq:conservationLaw}, i.e. with $f^x(q) = J_1 q$, $f^y(q) = J_2 q$, $f^z(q) = J_3 q$.

\subsection{Update of the Point Values}
The following description of the update of the point values holds on 2-d and 3-d Cartesian grids. For simplicity, the formulas are mostly given for the 2-d case only.

For updating the point values, a biparabolic reconstruction $$q_{\mathrm{recon}, ij}:\left[-\frac{\Delta x}{2}, \frac{\Delta x}{2}\right]\times\left[-\frac{\Delta y}{2}, \frac{\Delta y}{2}\right]\rightarrow\mathbb{R}^m, \qquad q_{\mathrm{recon}, ij} \in (P^{2,2})^m$$ is built for each cell. It has to satisfy
\begin{equation} 
\begin{aligned}
\frac{1}{\Delta x \Delta y}\int_{-\frac{\Delta x}{2}}^{\frac{\Delta x}{2}}\int_{-\frac{\Delta y}{2}}^{\frac{\Delta y}{2}}q_{\mathrm{recon}, ij}(\mathbf{x})\,\mathrm{d}\mathbf{x}&=q_{0, ij}\\
q_{\mathrm{recon}, ij}(\mathbf{x}_{r})&=q_{r, ij} \qquad r=1, ..., 8
\end{aligned}
\end{equation}
where $q_{0, ij}$ is the average over cell $\mathcal{C}_{ij}$ and $r$ is indexing all eight point values on the boundary of this cell, starting in the lower left corner with $q_{1, ij}$ and running through all point values counterclockwise (see Figure \ref{fig:distributionpointvalues}). As is customary for Finite Elements, one can compute shape functions and write the reconstruction as a linear combination of these:
\begin{equation}
q_{\mathrm{recon}, ij}(x, y)\coloneqq\sum_{r\in \{0, 1, ..., 8\}}q_{r, ij}B_r(x, y) \in (P^{2,2})^m. \label{eq:reconstruction}
\end{equation}
The basis functions involved in equation \eqref{eq:reconstruction} are given in Appendix \ref{app:basisfcts}; see Figure \ref{fig:basisfunctions2d} for an illustration of $B_7$, $B_8$ and $B_0$ as well as for an example of the reconstruction \eqref{eq:reconstruction}. 

\begin{figure}
 \centering
 \includegraphics[width=0.8\textwidth]{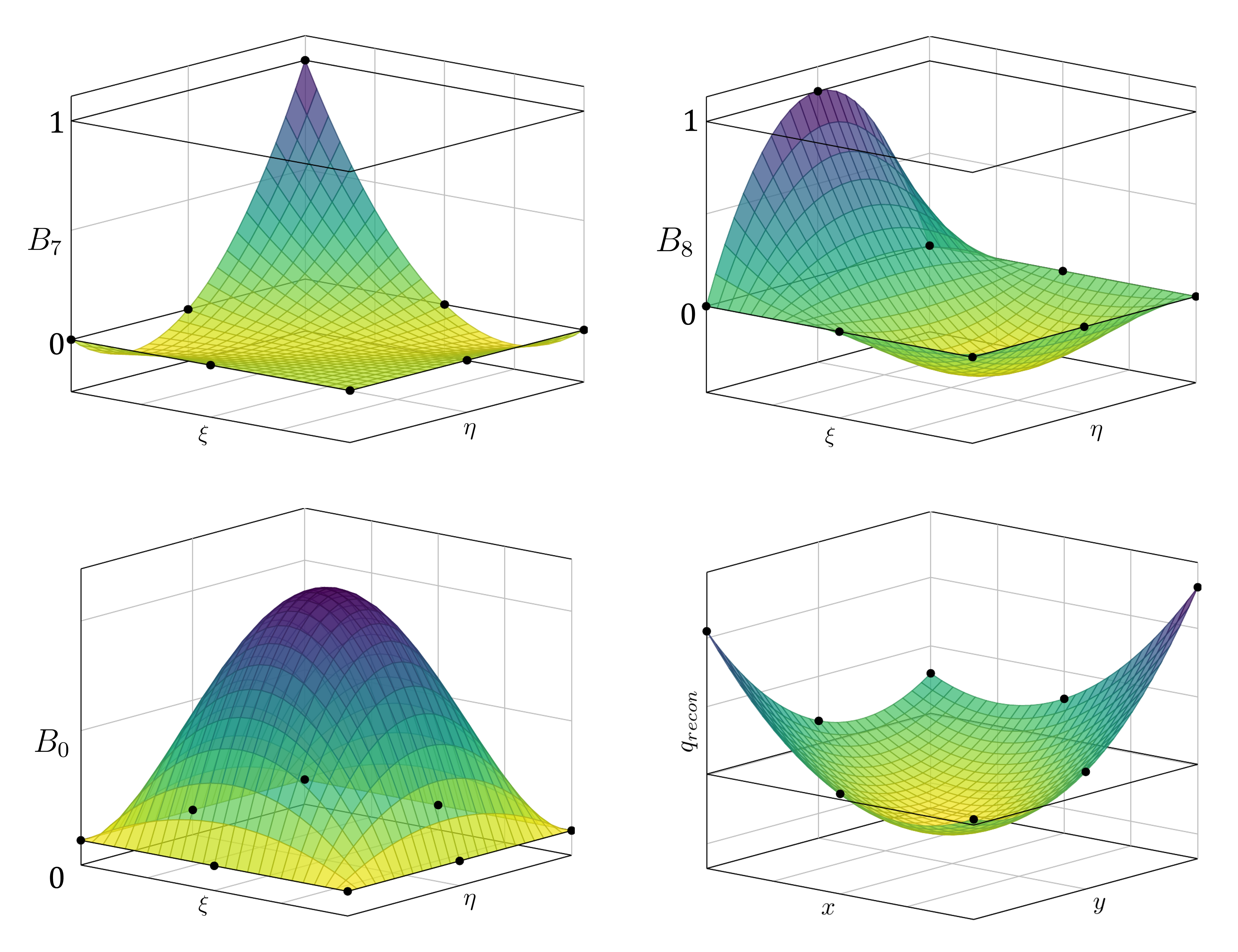} 
 \caption{Shape functions $B_7$ (\emph{top left}), $B_8$ (\emph{top right}), $B_0$ (\emph{bottom left}) in two spatial dimensions and a reconstruction devised by using these shape functions (\emph{bottom right}).}
 \label{fig:basisfunctions2d}
\end{figure}

Observe that a biparabolic reconstruction is parabolic along any of the edges of the (Cartesian) cell and that the parabola is uniquely defined by the three point values located at any edge. As the point values are shared with adjacent cells, one concludes that the above biparabolic reconstructions form a globally continuous reconstruction. 

For the point value update, we are aiming at a semi-discretization of \eqref{eq:conservationLaw}, i.e. a discretization in space is needed.
A Finite Difference approximation to the derivative that uses a compact stencil and is of maximal order of accuracy can be obtained by differentiating the reconstruction at the location of the respective point value. The reconstruction is not continuously differentiable in the direction perpendicular to the edges/faces, which is beneficial, as it allows to include upwinding: the derivative is taken from the cell in upwind direction.

Observe that the choice of degrees of freedom that belong to the cell $\mathcal C_{ij}$ privileges the right upper corner. This means that for positive $x$, $y$ velocity components, $\mathcal C_{ij}$ is the upwind cell for all degrees of freedom $q^P_{ij}$, $P \in \{ \nodeonly, \ehonly, \evonly   \} $. We therefore denote by $D_x \big \rvert^{P}_{ij} q$ the $x$-derivative of the reconstruction \eqref{eq:reconstruction} in cell $\mathcal{C}_{ij}$ evaluated at the location of the point value $P$:

\begin{align}
 D_x \big \rvert^{P}_{ij}q := \frac{\del}{\del x} q_{\text{recon},ij}\big \rvert_{\vec x = \vec x^P} \in \mathbb R^m. \label{eq:findiffderivplus}
\end{align}

For negative velocity components, in some cases one differentiates the reconstruction from cell $\mathcal{C}_{i+1, j}$ at the same location, e.g. for the $x$-derivative
\begin{align}
D_x^* \big \rvert^{P}_{ij}q := \begin{cases} \displaystyle \frac{\del}{\del x} q_{\text{recon},i+1,j} \big \rvert_{\vec x = \vec x_{ij}^P - x_{i+1,j}} & P \in \{ \nodeonly, \evonly \} \\
                                \displaystyle D_x \big \rvert^{P}_{ij}q & P \in \{ \ehonly \}
                               \end{cases} \label{eq:findiffderivminus}
\end{align}

\begin{figure}
 \centering
 \includegraphics[width=0.7\textwidth]{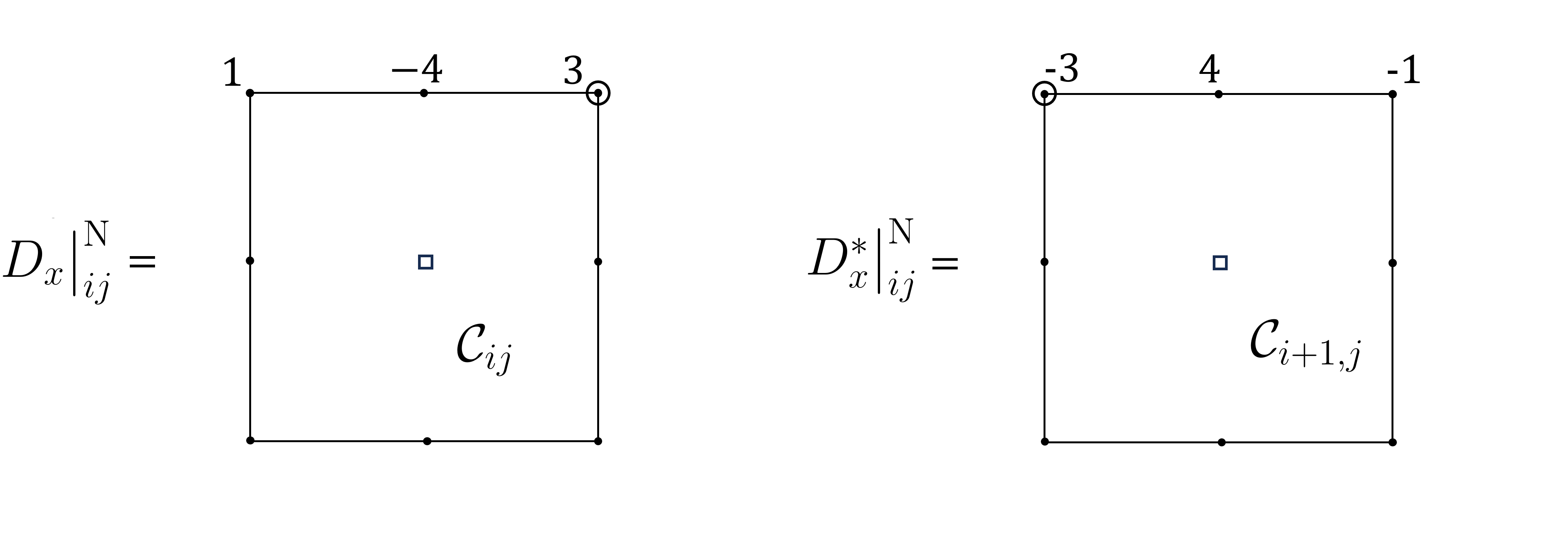}\\
 \includegraphics[width=0.7\textwidth]{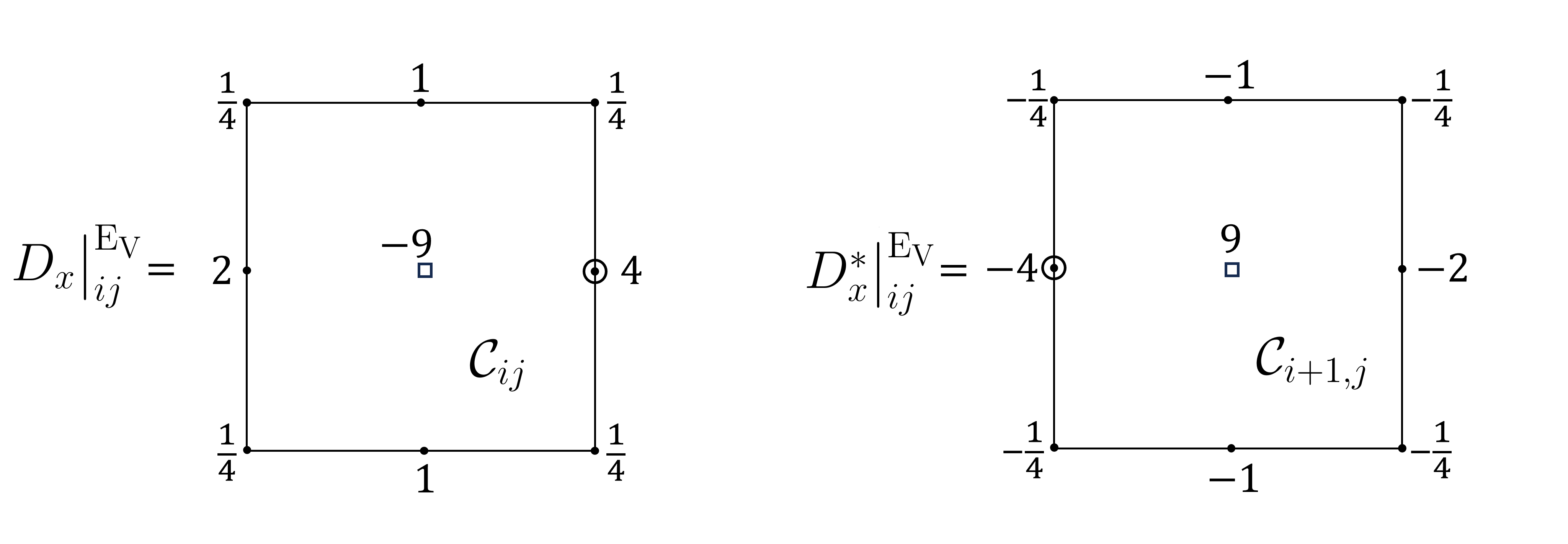}\\
 \includegraphics[width=0.7\textwidth]{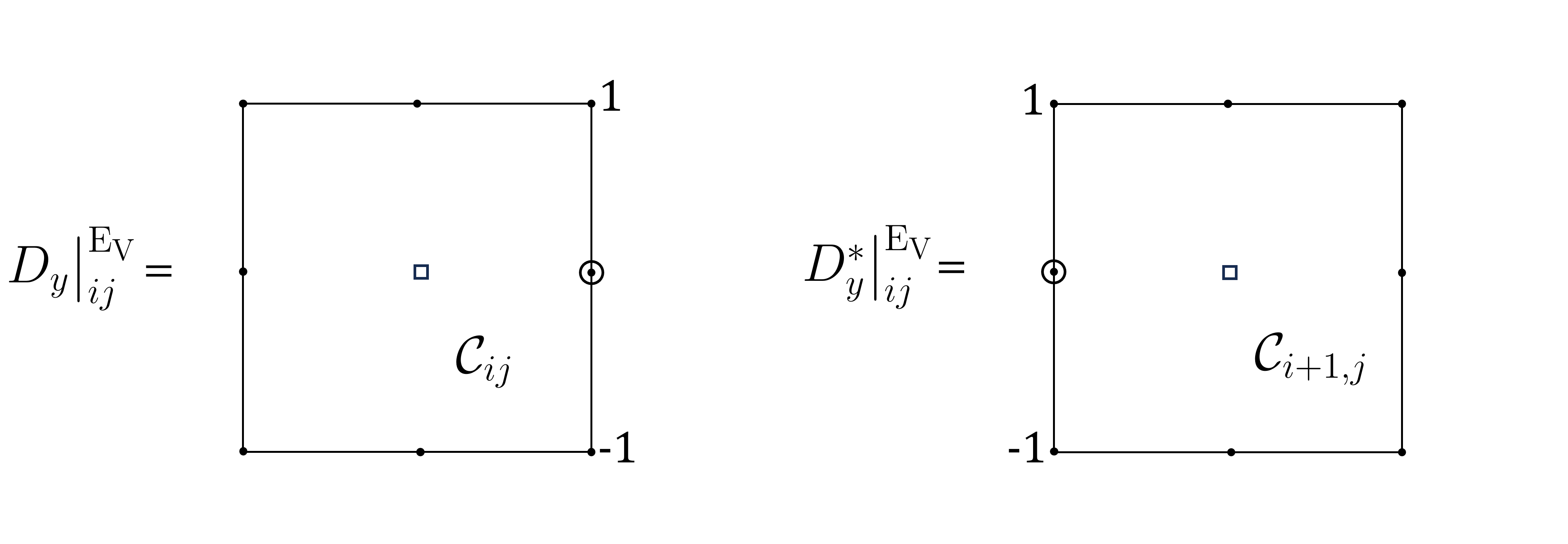} \\
 \includegraphics[width=0.7\textwidth]{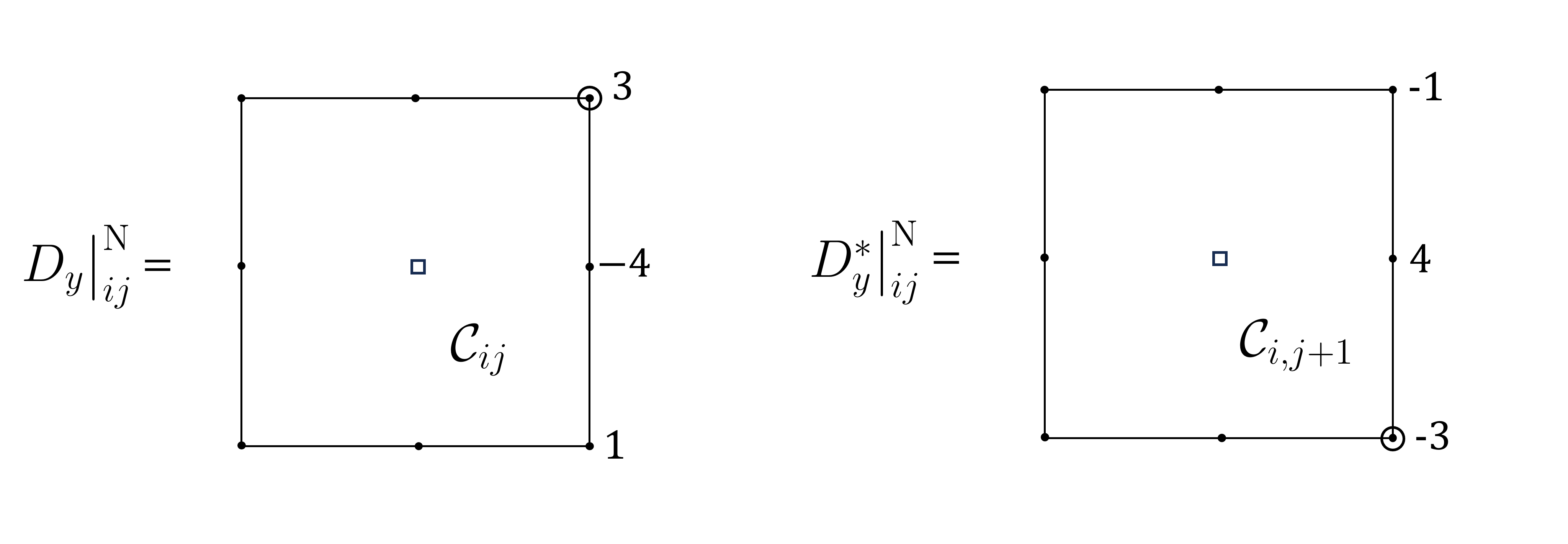}
 \caption{Illustration of the coefficients occurring in the difference formulas for node values and vertical edge values. They are placed at the location of the corresponding point value or at the cell center for the average. The circled dot marks the location of the point value intended to be updated. Cells $\mathcal{C}_{ij}$ (\emph{left}) and $\mathcal{C}_{i+1, j}$ (rows 1--3) / $\mathcal{C}_{i, j+1}$ (row 4) (\emph{right}) are shown. \emph{Upper row:} Difference formulas $D_x \big \rvert\node_{ij}$ and $D^*_x \big \rvert\node_{ij}$. \emph{Second row:} Difference formulas $D_x\big \rvert\ev_{ij}$ and $D^*_x\big \rvert\ev_{ij}$. \emph{Third row}: Difference formulas $D_y\big \rvert\ev_{ij}$ and $D^*_y\big \rvert\ev_{ij}$, which are equivalent. \emph{Fourth row:} Difference formulas $D_y\big \rvert\node_{ij}$ and $D^*_y\big \rvert\node_{ij}$.   }
 \label{fig:difference_formulas_2d}
\end{figure}

The same applies to the formulas $D_y\big \rvert^{P}_{ij}$ and $D_y^*\big \rvert^{P}_{ij}$, and analogously in 3-d. Figure \ref{fig:difference_formulas_2d} shows the coefficients occurring in the difference formulas assigned to the corresponding point value or the cell average.

The update of a point value is in 2-d
\begin{subequations}
\begin{align}
\frac{\mathrm{d}}{\mathrm{d}t}q^{P}_{ij}&=-\left(J_x^{+}  D_x \big \rvert^{P}_{ij}q+J_x^{-}D^*_x \big \rvert^{P}_{ij}q\right)-\left(J_y^{+}D_y \big \rvert^{P}_{ij}q+J_y^{-}D^{*}_{y}\big \rvert^{P}_{ij}q\right),\\
&\phantom{mm} \nonumber P \in \{ \nodeonly, \ehonly, \evonly   \} 
\end{align}
and in 3-d 
\begin{align}
\frac{\mathrm{d}}{\mathrm{d}t}q^{P}_{ijk}&=-\left(J_x^{+}  D_x \big \rvert^{P}_{ijk}q+J_x^{-}D^*_x \big \rvert^{P}_{ijk}q\right)-\left(J_y^{+}D_y \big \rvert^{P}_{ijk}q+J_y^{-}D^{*}_{y}\big \rvert^{P}_{ijk}q\right) \\
&\nonumber \phantom{mmmmmmmmmmmm}-\left(J_z^{+}D_z \big \rvert^{P}_{ijk}q+J_z^{-}D^{*}_{z}\big \rvert^{P}_{ijk}q\right), \\
&\phantom{mm} \nonumber P \in \{ \nodeonly, \exonly, \eyonly, \ezonly, \fxonly, \fyonly, \fzonly  \} 
\end{align}   \label{eq:updatePointValues}
\end{subequations}
where $J_x^\pm$, $J_y^\pm$ are the positive/negative parts of $J_x$, $J_y$. Given the diagonali\-zation
\begin{align}
J_x&=R \mathrm{diag}(\lambda_1, ..., \lambda_m)R^{-1} 
\end{align}
one defines
\begin{align}
J_x^{\pm}&:=R\mathrm{diag}(\lambda_1^{\pm}, ..., \lambda_m^{\pm})R^{-1} 
\label{eq:jacobiansplit}
\end{align}
and analogously for $J_y$.

The above point value update is a Jacobian splitting inspired by the 1-d upwind method. One might consider an alternative point value update
\begin{align}\label{eq:alternativePointValueUpdate}
\tilde J_x^\pm & := J_x\pm a_x \id_m & \tilde J_y^\pm &:= J_y\pm a_y\id_m
\end{align} 
where $a_x=\max(\rvert\lambda(J_x)\rvert)$, $a_y=\max(\rvert\lambda(J_y)\rvert)$ with $\lambda$ being an eigenvalue of $J_{x}$ or $J_y$. This Jacobian splitting is inspired by the Rusanov method. As will be seen later, when it comes to structure preservation it is significantly inferior to the upwind splitting.

\subsection{Integration in Time}

The update equations \eqref{eq:updateAverages} and \eqref{eq:updatePointValues} for the average and the point values are evolved in time using Runge-Kutta time integration. One can expect a maximal CFL number of about 0.2, half of what is known for the 1-d situation (\cite{abgrall22proceeding}). The analysis of Section \ref{sec:numdiff} indicates that for linear acoustics it actually is $\sim$0.28. The stationary states of the numerical method are by defnition those for which the space derivative vanishes, i.e. one is left in this case with the ordinary differential equation $q'(t) = 0$. Any Runge-Kutta method is able to integrate it exactly, and the choice of the time integrator is thus not relevant for stationarity preservation.

\section{The discrete Fourier transform} \label{sec:fourier}

\subsection{Fourier modes}

Stationarity preservation of linear numerical schemes for linear hyperbolic systems \eqref{eq:conservationLaw} can be investigated using the discrete Fourier transform, i.e. by expressing the spatial dependence of any grid function $q_{ij}$ (or $q_{ijk}$), $q \colon \mathbb Z^d \to \mathbb R^m$ as the linear combination 
\begin{align}
q_{ij} &= \sum_{\vec k} \hat q(\vec k) \exp(\ii k_x i \Delta x + \ii k_y j \Delta y) \qquad \text{or}\\
q_{ijk} &= \sum_{\vec k} \hat q(\vec k) \exp(\ii k_x i \Delta x + \ii k_y j \Delta y + \ii k_z k \Delta z).
\end{align}
The exponentials are discrete versions of $\exp(\ii \vec k \cdot \vec x)$, having written $\vec k = (k_x, k_y)$, or $\vec k = (k_x, k_y, k_z)$. This is possible because the method under consideration is linear and because we assume an equidistant Cartesian grid. The coefficients $\hat q(\vec k) \in \mathbb C^m$ of the linear combination are called the \emph{discrete Fourier transform of $q$}.

According to Definition \ref{def:belong}, each of the degrees of freedom that belong to a cell forms a lattice with spacings $\Delta x$, $\Delta y$ (and $\Delta z$, in 3-d), i.e. can be seen as a grid function. When investigating the semi-discrete Active Flux method using the discrete Fourier transform, each lattice is associated with its own Fourier mode. 
We express every degree of freedom belonging to a cell $\mathcal C_{ij}$ in 2-d as 
\begin{align}
 q^X_{ij}(t) &= \sum_{\vec k} \hat q^X(t, \vec k) \exp(\ii i k_x \Delta x + \ii j k_y \Delta y) && X \in \{ \avgonly, \nodeonly, \ehonly, \evonly \} 
\end{align}
and to a cell $\mathcal C_{ijk}$ in 3-d as
\begin{align}
 q^X_{ijk}(t) &= \sum_{\vec k}  \hat q^X(t, \vec k) \exp(\ii i k_x \Delta x + \ii j k_y \Delta y + \ii k k_z \Delta z) \\\nonumber &X \in \{ \avgonly, \nodeonly, \exonly, \eyonly, \ezonly, \fxonly, \fyonly, \fzonly \}
\end{align}
Here $\hat q^X(t, \vec k) \in \mathbb C^m$. 
It is useful to introduce the following
\begin{definition}[Translation factor]
The \emph{translation factors} are given as $t_x(k_x) = \exp(\ii k_x \Delta x)$, $t_y(k_y) = \exp(\ii k_y \Delta y)$,$t_z(k_z) = \exp(\ii k_z \Delta z)$. We drop explicit mention of the parameters in the following.
\end{definition}
Index shifts lead to multiplications with shift factors:
\begin{align}
 q^X_{i+I,j+J}(t) &= \sum_{\vec k}t_x^I t_y^J \cdot \hat q^X(t, \vec k) \exp(\ii i k_x \Delta x + \ii j k_y \Delta y) \\
 &\phantom{mmmmmmmmmmmmm}\nonumber X \in \{ \avgonly, \nodeonly, \ehonly, \evonly \} \\
 q^X_{i+I,j+J,k+K}(t) &= \sum_{\vec k} t_x^I t_y^J t_z^K \cdot \hat q^X(t, \vec k) \exp(\ii i k_x \Delta x + \ii j k_y \Delta y + \ii k k_z \Delta z) \\
 &\phantom{mmmmmmmmmmmmm}\nonumber X \in \{ \avgonly, \nodeonly, \exonly, \eyonly, \ezonly, \fxonly, \fyonly, \fzonly \},
\end{align}
where $t_x^I$ is the $I$-th power of the unit complex number $t_x = \exp(\ii k_x \Delta x)$, etc.

Both the summation and the global exponential factor $\exp(\ii i k_x \Delta x + \ii j k_y \Delta y)$ will eventually drop out of all expressions:

\begin{proposition} \label{prop:cancleoutfactor}
 Consider grid functions $q, Q \colon \mathbb Z^2 \to \mathbb R^m$ and assume that $Q$ is linear in $q$, i.e. an expression of the form
 \begin{align}
  Q_{ij} = \sum_{(I,J) \in \mathbb Z^2} \alpha_{IJ} q_{i+I, j+J} \qquad \alpha_{IJ} \in \mathscr M^{m\times m}(\mathbb R) \quad \forall I,J.
 \end{align}
 Then by writing 
 \begin{align}
 q_{ij} &:= \sum_{t_x,t_y} \hat q(t_x,t_y) t_x^i t_y^j, \quad\hat q(t_x,t_y) \in \mathbb C^m & Q_{ij} &:=\sum_{t_x,t_y}  \hat Q(t_x,t_y) t_x^i t_y^j \label{eq:commenfactor}
 \end{align}
 one obtains
 \begin{align}
 \hat Q = \sum_{(I,J) \in \mathbb Z^2} t_x^I t_y^J \, \alpha_{IJ} \hat q  \in \mathbb C^m.
 \end{align}
 The 3-d case is analogous.
\end{proposition}
\begin{proof}
 By direct computation and using linearity.
\end{proof}

Note that $\hat Q$ depends on $\vec k$, or equivalently on $t_x,t_y$, but that we will neglect explicit mention of these parameters.

\begin{definition}[Discrete Fourier transform] \label{def:fourier}
 Given $q_{ij}$ as in Proposition \ref{prop:cancleoutfactor}, we call $\hat q$ the \emph{discrete Fourier transform of $q$}. We will neglect the common factor $t_x^i t_y^j$ in \eqref{eq:commenfactor} as it is bound to cancel in the end and simply write $q_{ij} \equiv \hat q$, $q_{i+1,j} \equiv \hat q t_x$ etc. as a formal replacement rule.
 \end{definition}

For the degrees of freedom accessible to one cell for 2-d Active Flux one finds the Fourier transforms by applying the replacement rule of Definition \ref{def:fourier} to Equations \eqref{eq:2ddofs1}--\eqref{eq:2ddofs3}:
\begin{align}
 &q_{7,ij} \equiv \hat q\node/t_x &&q_{6,ij} \equiv \hat q\eh &&q_{5,ij} \equiv \hat q\node\\
 &q_{8,ij} \equiv \hat q\ev/t_x &&q_{0,ij} \equiv \hat q\avg && q_{4,ij} \equiv \hat q\ev \\
 &q_{1,ij} \equiv \hat q\node/(t_xt_y) &&q_{2,ij} \equiv \hat q\eh/t_y &&q_{3,ij} \equiv \hat q\node / t_y \label{eq:replacementrulefourier}
\end{align}

The Fourier transforms for Active Flux are collected in a block-vector $\hat{q}$, where
\begin{equation}
\hat{q}=(\hat{q}\avg, \hat{q}\eh, \hat{q}\ev, \hat{q}\node) \in \mathbb C^{m N^\text{dof}} =\mathbb C^{m \cdot 4}
\end{equation}
 for the two-dimensional case and 
 \begin{equation}
\hat{q}=(\hat{q}\avg, \hat{q}\ex, \hat{q}\ey, \hat{q}\ez, \hat{q}\fx, \hat{q}\fy, \hat{q}\fz, \hat{q}\node)  \in \mathbb C^{m N^\text{dof}} = \mathbb C^{m \cdot 8}
\end{equation}
for the three-dimensional case. Here $\hat{q}\avg \in \mathbb C^m$ is the Fourier transform of the average, $\hat{q}^{\text{E}_i} \in \mathbb C^m$ is the Fourier transform of the point value on an edge parallel to the $i$-axis and $\hat{q}^{\text{F}_i} \in \mathbb C^m$ is the Fourier transform of the point value on the face orthogonal to the $i$-axis.

\subsection{The Fourier transform of the reconstruction} \label{sec:fourierrecon}

\begin{definition}
 Consider the reconstruction $q_{\text{recon}, ij} \in (P^{k,k}(\mathbb R^2))^m$ (i.e. the one in \eqref{eq:reconstruction}) as a function of the degrees of freedom accessible to cell $\mathcal C_{ij}$, i.e. write 
 \begin{align}
  q_{\text{recon}, ij}(x, y; q_{0, ij}, q_{1,ij}, \ldots, q_{N^\text{dof}_{\text{acc}}-1, ij}).
 \end{align}
 We denote by $\widehat{q_{\text{recon}}} \in (P^{k,k}(\mathbb C^2))^m$ the same polynomial obtained by replacing the degrees of freedom by their Fourier transforms according to the replacement rule of Definition \ref{def:fourier}. For example, in 2-d we have, with the replacement rule \eqref{eq:replacementrulefourier}
 \begin{align}
  \widehat{q_{\text{recon}}}&(x, y; \hat q\avg, \hat q\eh, \hat q\ev,\hat q\node) \\&:= q_{\text{recon}, ij}\left(x, y;\hat q\avg , \frac{\hat q\node}{t_xt_y},   \frac{\hat q\eh}{t_y} ,  \frac{\hat q\node }{t_y}, \hat q\ev,  \hat q\node, \hat q\eh,  \frac{\hat q\node}{t_x} ,\frac{\hat q\ev}{t_x} \right).
 \end{align}
 We will generally omit explicit mention of the parameters in the following. 
The 3-d case is analogous.
 
\end{definition}

The reconstruction polynomial $q_{\text{recon},ij}(x, y)$ depends on all parameters \emph{accessible} to the cell (9 in case of \eqref{eq:reconstruction}), but its discrete Fourier transform $\widehat{q_\text{recon}}(x, y)$ depends only on those \emph{belonging} to one cell (4 in case of \eqref{eq:reconstruction}), i.e. $N^\text{dof}$. $\widehat{q_\text{recon}}$ can also simply be seen as the discrete Fourier transform pointwise at each $(x, y)$.

\subsection{The evolution matrix}

Equations \eqref{eq:updateAverages} and \eqref{eq:updatePointValues} for the average and the point value updates can be combined into
\begin{align}
 \frac{\dd}{\dd t} q_{I}^X + \sum_{Y \in \Sigma} \sum_{S\in[-N, N]^d\subset\mathbb{Z}^d} \alpha^X_{Y,S} q^Y_{I+S} &= 0\qquad \forall X \in \Sigma \label{eq:finiteDifference}
\end{align}
where 
\begin{align}
 \Sigma &= \{ \avgonly, \nodeonly, \ehonly, \evonly \} & I &= (i,j) \qquad (\text{2-d})\\
 \Sigma &= \{ \avgonly, \nodeonly, \exonly, \eyonly, \ezonly, \fxonly, \fyonly, \fzonly \} & I &= (i,j,k) \qquad (\text{3-d})
\end{align}
and $\alpha^X_{Y,S} \in \mathscr M^{m \times m}(\mathbb R)$.
Upon the Fourier transform one obtains
\begin{align}
 \frac{\dd}{\dd t} \hat q^X + \sum_{Y \in \Sigma} \sum_{S\in[-N, N]^d\subset\mathbb{Z}^d} \alpha^X_{Y,S} \hat q^Y \prod\limits_{m=1}^d t_m^{S_m} &= 0. \label{eq:generalevolutionmatrix}
\end{align}

For the Active Flux method applied to linear acoustics one has, for example for the node value:
\begin{align}
0=\frac{\mathrm{d}}{\mathrm{d}t}\hat{q}\node&+J_x^+\frac{1}{\Delta x}\left(3\hat{q}\node-4\hat{q}\eh+\hat{q}\node\frac{1}{t_x}\right)
+J_x^-\frac{1}{\Delta x}\left(-3\hat{q}\node+4\hat{q}\eh t_x-\hat{q}\node t_x\right)\\
&+J_y^+\frac{1}{\Delta y}\left(\hat{q}\node\frac{1}{t_y}-4\hat{q}\ev+3\hat{q}\node\right)
+J_y^-\frac{1}{\Delta y}\left(-\hat{q}\node t_y+4\hat{q}\ev t_y-3\hat{q}\node\right).
\end{align}
Further such equations are given in Appendix \ref{app:fourierequations}.

\begin{definition}[Evolution matrix]
The evolution matrix associated to the Finite Difference scheme \eqref{eq:finiteDifference} is the block matrix $\mathcal{E}(\mathbf{k}) \in \mathscr M^{mN^\text{dof} \times mN^\text{dof}}(\mathbb C)$ with its block-entry $(X,Y)$ being
\begin{equation}
\Big ( \mathcal{E}(\mathbf{k}) \Big)_{X,Y}= \sum_{S\in[-N, N]^d\subset\mathbb{Z}^d} \alpha^X_{Y,S}   \prod\limits_{m=1}^d t_m^{S_m}
\end{equation}
such that \eqref{eq:generalevolutionmatrix} can be written as
\begin{equation}
\frac{\dd}{\dd t}\hat{q}+\mathcal{E}(\mathbf{k})\hat{q}=0.
\label{eq:semidiscreteFourier}
\end{equation}\end{definition}
We will frequently omit the argument and simply write $\mathcal E$.

For 2-dimensional linear acoustics we choose to write the evolution block matrix as follows
\begin{equation}
\mathcal{E}=\left(
\begin{array}{cccc}
\mathcal{E}_{\avgonly\avgonly} \hspace{2mm} &\mathcal{E}_{\avgonly\ehonly} \hspace{2mm} & \mathcal{E}_{\avgonly\evonly} \hspace{2mm}& \mathcal{E}_{\avgonly\nodeonly} \vspace{3mm} \\  
\mathcal{E}_{\ehonly\avgonly} \hspace{2mm} & \mathcal{E}_{\ehonly\ehonly} \hspace{2mm}& \mathcal{E}_{\ehonly\evonly} \hspace{2mm}& \mathcal{E}_{\ehonly\nodeonly} \vspace{3mm}\\
\mathcal{E}_{\evonly\avgonly}  \hspace{2mm} & \mathcal{E}_{\evonly\ehonly} \hspace{2mm}& \mathcal{E}_{\evonly\evonly} \hspace{2mm}& \mathcal{E}_{\evonly\nodeonly} \vspace{3mm}\\
\mathcal{E}_{\nodeonly\avgonly} \hspace{2mm} & \mathcal{E}_{\nodeonly\ehonly} \hspace{2mm}& \mathcal{E}_{\nodeonly\evonly} \hspace{2mm}& \mathcal{E}_{\nodeonly\nodeonly} 
\end{array}\right)
\label{eq:evolutionMatrixBlock}
\end{equation}
where e.g. $\mathcal{E}_{\ehonly\evonly}$ states the influence of the point value at the $\evonly$ edge on the point value at a $\ehonly$ edge. The explicit form of the blocks can be found in Appendix \ref{sec:2devomatricesacoustics}. 

By comparison with \eqref{eq:Jk} one observes that the evolution matrix $\mathcal{E}(\mathbf{k})$ plays the role of $\ii \vec J \cdot \vec k$. Formulation \eqref{eq:semidiscreteFourier} of the method therefore is ideal for studying the numerical stationary states, as they are given by the kernel of $\mathcal E(\vec k)$. 

While the right kernel allows to identify stationary states, the left kernel of $\mathcal E$ is related to involutions: If there exists an $\omega(\vec k) \in \mathbb C^{m N_\text{dof}}$ such that
\begin{align}
 \omega(\vec k)^\text{T} \mathcal E &= 0 \qquad \forall \vec k,
\end{align}
then
\begin{align}
 \frac{\dd}{\dd t} \omega(\vec k)^\text{T} \hat q &= 0.
\end{align}
In the case of linear acoustics $\omega(\vec k)^\text{T} \hat q$ would be the Fourier transform of a discretization of the vorticity $\nabla \times \vec v$. However, we find that in practice it is very difficult to explicitly compute the left kernel of $\mathcal E$. The existence of a right kernel guarantees the existence of a left kernel of the same dimension, but we are unable to say anything about the nature of the discrete vorticity beyond its existence.

\section{Using the discrete Fourier transform for the analysis of stationarity preservation} \label{sec:fourierAF}

\subsection{Stationarity preservation} \label{ssec:statpres}

For both stationarity preserving and stationarity non-preserving methods, on e.g. periodic domains the solution becomes stationary after long times, since von Neumann stability prohibits Fourier modes that grow in time: they can only decay or remain stationary. The difference between stationarity preserving and stationarity non-preserving methods is about how this final stationary state looks like.
If the discrete stationary states of the method discretize all the analytic stationary states of the PDE, then the method is called \emph{stationarity preserving}.
For linear acoustics, stationarity preserving methods have stationary states that are characterized by a consistent discretization of a divergence-free velocity without any further constraints and a vanishing pressure. Classical schemes mostly are not stationarity preserving: their discrete stationary state is instead a discretization of $\del_x u = 0$, $\del_yv= 0$ (see \cite{barsukow17a},\cite{barsukow18thesis}). Such stationary states are called trivial and are a small subset of divergencefree velocities, and e.g. vortices are not contained in it. Thus, stationarity non-preserving methods are only able to preserve shear flows, but cannot consistently discretize vortices, or generally speaking they are not able to offer discretizations of \emph{all} the stationary states of the PDE. Moreover, for stationarity non-preserving methods grid refinement only slows down the (exponentially quick) transition to the ``bad'' stationary state, but does not improve the stationary state itself, i.e. one can say that they lose consistency at stationary state.

A discrete Fourier mode with spatial frequency $\vec k$ is stationary if it is in the kernel of $\mathcal E(\vec k)$ (see \eqref{eq:semidiscreteFourier}).
The focus on non-trivial stationary states means that one is keen on identifying stationary modes for any $\vec k$ (i.e. for general $t_x, t_y, t_z$). The kernel thus obtained shall then be compared to the kernel of $\vec J \cdot \vec k$. 

For Finite Difference methods, the condition for a method to be stationarity preserving is (\cite{barsukow17a})
\begin{align}
 \min_{\vec k} \dim \ker \mathcal E(\vec k) = \min_{\vec k} \dim \ker \vec J \cdot \vec k
\end{align}
For numerical methods with multiple degrees of freedom per cell, such as Active Flux, the dimension of the space on which $\mathcal E$ operates is $N^\text{dof}$ times larger than that of $\vec J \cdot \vec k$. One might thus impose a correspondingly higher $\dim \ker \mathcal E$ as a condition for stationarity preservation. We shall, however, be modest here:

\begin{definition}
 A linear numerical method with $N^\text{dof}$ degrees of freedom per cell per variable is called stationarity preserving if 
 \begin{align}
 S \leq \min_{\vec k} \dim \ker \mathcal E(\vec k) \leq  S N^\text{dof}
\end{align}
where $S := \min_{\vec k} \dim \ker \vec J \cdot \vec k$.
\end{definition}

\subsection{A review of stationarity preservation of the classical Active Flux method} \label{sec:traditionalfourier}

For 2-d linear acoustics, the matrix $\mathcal E$ of Active Flux is a (dense) $12 \times 12$ complex-valued matrix. Although only its kernel is required (and not a full diagonalization) it remains a formidable task to compute it. In \cite{barsukow18activeflux} therefore a slightly easier approach has been used, and the techniques used therein allow some deeper understanding of the kernel. They are useful here as well. 

The object under consideration in \cite{barsukow18activeflux} is classical Active Flux, where point values are updated by means of an evolution operator. The latter maps the reconstruction (used as initial data) directly to the value at the new time, and the \emph{exact} evolution operator from \cite{barsukow17} has been used in \cite{barsukow18activeflux}. Clearly then, a numerical stationary state arises if the reconstruction ($=$ initial data for the evolution operator) is of constant pressure and divergence-free velocity. One thus needs to analyze the subset of divergence-free reconstructions.

The reconstruction is in $P^{2,2}$ in every variable and thus the divergence is in
\begin{align}
 D_\text{br}^{2,2} := \{ v \in L^\infty : v \vert_{C_{ij}} \in \mathrm{span}\{1,x,x^2,y,xy, x^2y, y^2, xy^2\} \}, \label{eq:Dbroken}
\end{align}
i.e. in an 8-dimensional space (observe the absence of $x^2y^2$). Due to continuity, the 9 degrees of freedom accessible to a cell, which are used in writing the reconstruction polynomial, cannot all be chosen independently in every cell. It is better at this stage to switch to the discrete Fourier transform $\widehat{q_\text{recon}}$ of the reconstruction, introduced in Section \ref{sec:fourierrecon}. This allows to factor out continuity, and let appear the dependence of the reconstruction on only 4 degrees of freedom per variable that belong to a cell. These can indeed now be chosen freely in each cell.

Seeking a divergence-free reconstruction in 2-d, one thus ends up with 8 equations for 8 free parameters. Without rank defect the only divergence-free reconstruction would be trivial, and not representative of the richness of divergence-free vector fields. The linear system turns out to be \emph{not} of full rank, however. This remained without a clear explanation in \cite{barsukow18activeflux}, and is elucidated below in Theorem \ref{thm:kernel2dacoustics}. It was found (see equations (6.25) and (6.26) in \cite{barsukow18activeflux}) that as long as the 12 Fourier modes are parallel to 

{\footnotesize
\begin{align}
  \hat Q = \bigg( &-\frac23 \frac{1 + 4 t_x + t_x^2}{t_x} \cdot \frac{(t_y-1)(t_y+1)}{\Delta y t_y} 
   , &&\frac23 \frac{1 + 4 t_y + t_y^2}{t_y}  \cdot \frac{(t_x-1)(t_x+1)}{\Delta x t_x}, &&0, \nonumber  \\ \nonumber
   &- \frac{1 + 6 t_x + t_x^2}{t_x} \cdot  \frac{t_y-1}{\Delta y},
   &&2 \frac{(t_x-1) (t_x+1)}{\Delta x t_x} (t_y+1) , &&0, \\\nonumber
   &-2 (t_x+1) \frac{(t_y-1)(t_y+1)}{\Delta y t_y},
   &&\frac{t_x-1}{\Delta x}  \cdot \frac{1 + 6 t_y + t_y^2}{t_y}, &&0, \\
   &- 4 (t_x+1)\frac{t_y-1}{\Delta y}
   , &&4 \frac{t_x-1}{\Delta x} (t_y+1) ,&&0 \bigg ) \label{eq:stationaryeigenvector2dac}
 \end{align}}
 the divergence in each cell vanishes. Here, the variables are sorted as $(\hat q\avg,  \hat q\eh    , \hat q\ev,\hat q\node)$. One can easily see that this condition indeed can be rewritten as vanishing Finite Difference discretizations of the divergence. Writing $\hat q^X = (\hat u^X, \hat v^X, \hat p^X)$ with $X \in \{ \avgonly, \nodeonly, \ehonly, \evonly \}$ one finds e.g.

\begin{align}
\frac{1+4t_y+t_y^2}{t_y}\frac{(t_x-1)(1+t_x)}{t_x\Delta x}\hat{u}\avg+\frac{1+4t_x+t_x^2}{t_x}\frac{(t_y-1)(1+t_y)}{t_y\Delta y}\hat{v}\avg=0  \label{eq:fourierVelAverage}\\
\frac{(t_x-1)(1+t_y)}{\Delta x}\hat{u}\node+\frac{(1+t_x)(-1+ty)}{\Delta y}\hat{v}\node=0 \label{eq:fourierVelNode}\\
\frac{(1+6t_y+t_y^2)(t_x-1)}{\Delta x t_y}\hat{u}\eh+\frac{(1+6t_x+t_x^2)(t_y-1)}{\Delta y t_y}\hat{v}\ev=0 \label{eq:fourierVelHorizontal}\\
\frac{t_x-1}{\Delta x t_x}\hat{u}\ev+\frac{t_y-1}{\Delta y t_y}\hat{v}\eh=0 \label{eq:fourierVelVertical}\\
\frac{t_x-1}{\Delta x t_x} \hat{u}\ev + \frac{(t_y-1)(t_x+1)}{2\Delta y t_x t_y} \hat{v}\node = 0 \label{eq:fourierdiv5}\\
\frac{(1 + 6 t_y + t_y^2)(t_x-1)}{8 \Delta x t_y} \hat{u}\ev + \frac{(t_y^2 -1)(t_x+1)}{4 \Delta y t_y} \hat{v}\ev = 0 \label{eq:fourierdiv6}\\
\frac{(1 + 4 t_y + t_y^2)(t_x-1)}{6 \Delta x t_x t_y} \hat{u}\ev + \frac{t_y^2 -1}{2 \Delta y t_y} \hat{v}\avg = 0 \label{eq:fourierdiv7}
\end{align}
These relations are Fourier transforms of 
\begin{align}
\frac{\langle[u\avg]_{i\pm1}\rangle_{j}^{(4)}}{\Delta x}+\frac{[\langle v\avg\rangle_{i}^{(4)}]_{j\pm1}}{\Delta y}&=0 & \frac{\{[u\node]_{i+\frac{1}{2}}\}_{j+\frac{1}{2}}}{\Delta x}+\frac{[\{v\node\}_{i+\frac{1}{2}}]_{j+\frac{1}{2}}}{\Delta y}&=0 \label{eq:discreteStatStatesVelocity1}\\
\frac{\langle[u\eh]_{i+\frac12}\rangle_{j}^{(6)}}{\Delta x}+\frac{[\langle v\ev\rangle_{i}^{(6)}]_{j+\frac12}}{\Delta y}&=0 & \frac{[u\ev]_{i-\frac12,j}}{\Delta x}+\frac{[ v\eh_{i}]_{j-\frac12}}{\Delta y}&=0\label{eq:discreteStatStatesVelocity2}\\
\frac{ [u\ev ]_{i+\frac12,j+1} }{\Delta x} + \frac{ [\{ v\node \}_{i+\frac12}]_{j+\frac12}}{2 \Delta y} &= 0 & \frac{ \langle [u\ev]_{i+\frac12} \rangle_j^{(6)}   }{8 \Delta x} + \frac{[ \{ v\ev \}_{i+\frac12}]_{j\pm1} }{4 \Delta y} &= 0 \label{eq:discreteStatStatesVelocity3}\\
\frac{ \langle [u\ev]_{i+\frac12} \rangle_j^{(4)}   }{6 \Delta x} + \frac{[ v\avg_{i+1}]_{j\pm1} }{2 \Delta y} &= 0 \label{eq:discreteStatStatesVelocity4}
\end{align}
having introduced the notation
\begin{align}
[q]_{i+\frac12}&=q_{i+1}-q_{i} & \{q\}_{i+\frac12}&=q_{i+1}-q_i \\
[q]_{i\pm1}&=q_{i+1}-q_{i-1} & \langle q\rangle_i^{(\alpha)}&=q_{i-1}+\alpha q_i+q_{i+1}.
\end{align}

\begin{theorem} \label{thm:biparabolicdivfree}
 Consider the set $V$ of globally continuous vector fields $\vec v$, biparabolic in each cell: 
 \begin{align}
 V = \left\{ \vec v \colon \mathbb R^2 \to \mathbb R^2 : \vec v\Big \rvert_{\mathcal C_{ij}} \in (P^{2,2})^2\right \}
 \end{align}
 with the degrees of freedom of Active Flux and periodic boundaries. Assume the grid to contain $N_\text{cells}$ cells. Then,
 \begin{enumerate}[(i)]
 \item $V \simeq \mathbb R^{8N_\text{cells}}$ \label{it:allvectors}
 \item $V_{\div} := \{ \vec v \in V : \nabla \cdot \vec v = 0 \} \simeq \mathbb R^{N_\text{cells} + \mathcal O(\sqrt{N_\text{cells}})}$ \label{it:vectorsdiv}
 \item Any element in $V_{\div}$ has continuous normal derivatives, i.e. at every cell interface with normal $\vec n$, $\vec v \cdot \vec n$ is continuously differentiable. \label{it:vectorsdivnormal}
 \end{enumerate}
\end{theorem}
\begin{proof}
 (\ref{it:allvectors}) and (\ref{it:vectorsdiv}) are clear from what has been said before; $\mathcal O(\sqrt{N_\text{cells}})$ takes into account that the number of divergence-free vector fields is slightly larger than one per cell, as there are some trivial ones missing in the previous analysis (indeed, no assumptions were made on $\vec k$ or $t_x,t_y$). A particular value of, say, $t_x$ would mean that at most, any grid function $q_{ij}$ is a function of $j$ only, i.e. that one can specify $\mathcal O(\sqrt{N_\text{cells}})$ values independently.
 
 The proof of (\ref{it:vectorsdivnormal}) is obtained by explicit computation. Denote first the $u$ and $v$ components of $\widehat{q_\text{recon}}$ by $\widehat{u_\text{recon}}$ and $\widehat{v_\text{recon}}$, respectively. If $\hat q$ is parallel to $\hat Q$ from \eqref{eq:stationaryeigenvector2dac}, then the $x$-derivative of $\widehat{u_\text{recon}}$ is proportional to
 \begin{align}
   \del_x \widehat{u_\text{recon}}(x,y) &= - 2\frac{(t_x-1)(t_y-1)}{\Delta x \Delta y t_x t_y} \left( (t_x+1) + 2 (t_x-1) \frac{x}{\Delta x} \right) \left( (t_y+1) + 2 (t_y-1) \frac{y}{\Delta y} \right )
 \end{align}
 One then immediately confirms that
 \begin{align}
   \del_x \widehat{u_\text{recon}} \left(\frac{\Delta x}{2}, y\right) &= t_x \del_x \widehat{u_\text{recon}} \left(-\frac{\Delta x}{2}, y\right),
 \end{align}
 i.e.
 \begin{align}
   \del_x u_{\text{recon},ij} \left(\frac{\Delta x}{2}, y\right) &= \del_x u_{\text{recon},i+1,j} \left(-\frac{\Delta x}{2}, y\right).
 \end{align}
 
 An analogous statement holds for $\del_y \widehat{v_\text{recon}}$.

\end{proof}

\subsection{Stationarity preservation of the semi-discrete method}

The following sections prove that the semi-discrete (generalized) Active Flux is stationarity preserving for the acoustic equations in two and three spatial dimensions. Additionally, the elements of the nullspace of the evolution matrix and the corresponding numerical stationary states are given and analyzed. Table \ref{tab:overview} sums up and compares key data of semi-discrete Active Flux for above-mentioned equations.\newline

\begin{table}[h]
 \centering
 \begin{tabular}[h]{l|c|c|c}
                                                                      & 2-d Acoustics & 3-d Acoustics  \\ \hline\hline 
 \# variables $ = m$ &                                                    3 &              4             \\\hline
 \# $\min_{\vec k} \dim \ker \vec J \cdot \vec k=S$ &                                           1 &              2            \\ \hline
 \# dof used in reconstruction $= N^\text{dof}_\text{acc}$ &                                        9    &27      \\\hline
 \# Fourier modes per variable $=N^\text{dof}$ &                                 4 &     8      \\\hline
 \# rows/columns in $\mathcal E = m N^\text{dof}$                       &        12 &               32         \\\hline\hline
 $\min_{\vec k} \dim \ker \mathcal E(\vec k)$ &                           1 &              5           
\end{tabular}
\caption{Overview of dimensions relevant for stationarity preservation of semi-discrete Active Flux for linear acoustics in 2 and 3 spatial dimensions.}
\label{tab:overview}
\end{table}

\subsubsection{Acoustic Equations in Two Spatial Dimensions}

Recall that the finite difference operators used in the update of the point values have been obtained by differentiating the reconstruction. 

Consider instead of \eqref{eq:updatePointValues} first the \textbf{central Active Flux} method, which reads
\begin{align*}
\frac{\mathrm{d}}{\mathrm{d}t}q^{P}_{ij}&=- J_x \frac12 \left( D_x \big \rvert^{P}_{ij}+D^*_x \big \rvert^{P}_{ij}\right)q-J_y\frac12\left(D_y \big \rvert^{P}_{ij}+D^{*}_{y}\big \rvert^{P}_{ij}\right)q\\
&= - \left ( \begin{array}{ccc} 0 & 0 & \frac12 \left( D_x +D^*_x \right)\big \rvert^{P}_{ij} \\
            0 & 0 & \frac12\left(D_y +D^{*}_{y}\right)\big \rvert^{P}_{ij} \\
            \frac12 \left( D_x +D^*_x \right)\big \rvert^{P}_{ij} & \frac12\left(D_y +D^{*}_{y}\right)\big \rvert^{P}_{ij} & 0 \end{array} \right ) \veccc{\phantom{\big \rvert^{P}_{ij}}\!\!\!\!\!\!\!u}{\phantom{\big \rvert^{P}_{ij}}\!\!\!\!\!\!\!v}{\phantom{\big \rvert^{P}_{ij}}\!\!\!\!\!\!\!p}
           \\
&\phantom{mm} \nonumber P \in \{ \nodeonly, \ehonly, \evonly   \} 
\end{align*}

The divergence relevant for the update of $q_{ij}^P$, $P \in \{ \nodeonly, \ehonly,\evonly \}$ is obtained by considering the velocity reconstructions in the two/four adjacent cells, taking their divergences at the location of $P$ (which gives two/four different values) and finally taking the mean of these values.
Recall that the image of $(P^{2,2})^2$ under the (weak) divergence is the broken space $D^{2,2}_\text{br}$ defined in \eqref{eq:Dbroken}. Obviously computing the mean of the two/four values restores continuity, i.e. afterwards there is only one value of the divergence associated with each point. This motivates the following

\begin{definition} \label{def:p22prjecteddiv}
 The \emph{$P^{2,2}$-projected divergence} $w \in P^{2,2}$ is defined by \eqref{eq:reconstruction} (with $m=1$) with the degrees of freedom given by
 \begin{subequations}
 \begin{align}
  w_{ij}\node &:= \frac14 \Big( W_{ij}(\vec x_{ij}\node - \vec x_{ij})+ W_{i+1,j}(\vec x_{ij}\node - \vec x_{i+1,j}) \\ \nonumber& \phantom{mmm}+ W_{i,j+1}(\mathbf x_{ij}\node - \vec x_{i,j+1}) + W_{i+1,j+1}(\mathbf x_{ij}\node - \vec x_{i+1,j+1})  \Big)\\
  w_{ij}\ev &:= \frac12 \left( W_{ij}(\mathbf x_{ij}\ev - \vec x_{ij})+ W_{i+1,j}(\mathbf x_{ij}\ev - \vec x_{i+1,j}) \right) \\
  w_{ij}\eh &:= \frac12 \left( W_{ij}(\mathbf x_{ij}\eh - \vec x_{ij})+ W_{i,j+1}(\mathbf x_{ij}\eh - \vec x_{i,j+1}) \right) \\
  w_{ij}\avg &:= \frac{1}{\Delta x \Delta y} \int_{-\frac{\Delta x}{2}}^{\frac{\Delta x}{2}} \int_{-\frac{\Delta y}{2}}^{\frac{\Delta y}{2}} W_{ij}(\vec x) \,\dd \vec x
 \end{align} \label{eq:averageingdoffordiv}
 \end{subequations}
 where $W_{ij}(\vec x) := (\del_x u_{\text{recon},ij} + \del_y v_{\text{recon},ij})(\vec x)$ is the divergence of the velocity reconstruction $(u_{\text{recon},ij}, v_{\text{recon},ij}) \in (P^{2,2})^2$ in $\vec x \in \left[-\frac{\Delta x}{2},\frac{\Delta x}{2}\right]\times \left[-\frac{\Delta y}{2},\frac{\Delta y}{2}\right]$.
\end{definition}

\begin{theorem}
 If the discrete data is such that the $P^{2,2}$-projected divergence of $\vec v$ vanishes (as a polynomial), and if $p$ is uniformly constant, then the central Active Flux method keeps this data stationary.
\end{theorem}
\begin{proof}
 By unisolvence of the $P^{2,2}$ space with the degrees of freedom specified in Definition \ref{def:p22prjecteddiv}, the vanishing of the divergence as a polynomial is equivalent to it being zero at the four degrees of freedom. By virtue of the definitions of the finite difference formulas (Equations \eqref{eq:findiffderivplus}--\eqref{eq:findiffderivminus}), the updates of $p_{ij}^P$, $P \in \{ \nodeonly, \ehonly,\evonly \}$ are the degrees of freedom of the $P^{2,2}$-projected divergence in $P$. The update equation \eqref{eq:updateaverage2d} of the pressure average is, by Gauss' law and the fact that the quadrature along the edges is exact for parabolas, equivalent to evaluating the average of the divergence $W_{ij}$ over the cell $(i,j)$. But this average is assumed to vanish as well, which concludes the proof.
\end{proof}

Central Active Flux can thus be associated with projecting the reconstruction of the divergence back into the space $P^{2,2}$. The map \eqref{eq:averageingdoffordiv} from the velocity variables to values of the $P^{2,2}$-projected divergence in its degrees of freedom is surjective, therefore the 4 equations for the 8 velocity variables leave a four-dimensional kernel of the evolution matrix. This sufficient condition of stationarity can be verified (using \textsc{mathematica}) to also be necessary, which proves
\begin{corollary}
 Central Active Flux for linear acoustics in 2 dimensions is stationarity preserving and the kernel of its evolution matrix is 4-dimensional.
\end{corollary}

Recall \eqref{eq:jacobiansplit}
with
\begin{align}
 \vert J_x\vert &= \left( \begin{array}{ccc} 1 & 0 &0 \\ 0&0&0 \\ 0&0&1 \end{array} \right) &
 \vert J_y\vert &= \left( \begin{array}{ccc} 0 & 0 &0 \\ 0&1&0 \\ 0&0&1 \end{array} \right)
\end{align}

The \textbf{upwind Active Flux} method \eqref{eq:updatePointValues} then writes
\begin{align}
\frac{\mathrm{d}}{\mathrm{d}t}q^{P}_{ij}
&= - \left ( \begin{array}{ccc} \frac12 \left( D_x -D^*_x \right)\big \rvert^{P}_{ij} & 0 & \frac12 \left( D_x +D^*_x \right)\big \rvert^{P}_{ij} \\
            0 & \frac12 \left( D_y-D^*_y \right)\big \rvert^{P}_{ij} & \frac12\left(D_y +D^{*}_{y}\right)\big \rvert^{P}_{ij} \\
            \frac12 \left( D_x +D^*_x \right)\big \rvert^{P}_{ij} & \frac12\left(D_y +D^{*}_{y}\right)\big \rvert^{P}_{ij} & \bowtie \end{array} \right ) \veccc{\phantom{\big \rvert^{P}_{ij}}\!\!\!\!\!\!\!u}{\phantom{\big \rvert^{P}_{ij}}\!\!\!\!\!\!\!v}{\phantom{\big \rvert^{P}_{ij}}\!\!\!\!\!\!\!p} \label{eq:updatePointValuesexplicit}
           \\
&\phantom{mm} \nonumber P \in \{ \nodeonly, \ehonly, \evonly   \} 
\end{align}
where the bow tie denotes further terms of no importance for the present discussion.
To gain insight into its stationarity preservation properties, it is useful to start with the
\begin{lemma}\label{lem:D22unsiolvence}
 $D^{2,2}_\text{br}$ is unisolvent with respect to the 8 pointwise degrees of freedom located at
 \begin{align}
  \left(x_{i\pm\frac12},y_{j\pm\frac12}\right), \hfill \left(x_{i\pm\frac12},y_j\right), \hfill \left(x_i,y_{j\pm\frac12}\right)
 \end{align}
\end{lemma}
\begin{proof}
 By explicit computation.
\end{proof}
Observe that these degrees of freedom are just the point values of Active Flux accessible to each cell, but without continuity.

\begin{theorem} \label{thm:kernel2dacoustics}
The semi-discrete upwind Active Flux method \eqref{eq:updatePointValues} for the linear acoustic equations in two spatial dimensions is stationarity preserving. A basis of the 1-dimensional nullspace of its evolution matrix is given by 
\begin{align}
\hat Q=\bigg(
&-\frac{\Delta x (1 + t_x (4 + t_x)) (-1 + t_y)}{
 6 \Delta y (-1 + t_x) t_x t_y},
 &&\frac{(1 + t_x) (1 + t_y (4 + t_y))}{6 t_x t_y (1 + 
   t_y)},&&0, \\ \nonumber
   &-\frac{\Delta x (1 + t_x (6 + t_x)) (-1 + t_y)}{
 4 \Delta y (-1 + t_x) t_x (1 + t_y)}, && \frac{1 + t_x}{2 t_x}, &&0, \\ \nonumber
 &-\frac{
 \Delta x (1 + t_x) (-1 + t_y)}{2 \Delta y (-1 + t_x) t_y},&& -\frac{-1-6t_y-t_y^2}{4t_y(1+t_y)}, &&0, \\ \nonumber
 &-\frac{\Delta x (1 + t_x) (-1 + t_y)}{\Delta y (-1 + t_x) (1 + t_y)},&& 1,&&0\bigg)^T\in\mathbb{C}^{12}. \label{eq:nullspaceAcoustics2d}
 \end{align} 
 The variables are ordered as $(\hat{q}\avg, \hat{q}\eh, \hat{q}\ev, \hat{q}\node)$. The divergence of the reconstruction vanishes iff $\hat q$ is parallel to $\hat Q$.
\end{theorem}
\begin{proof}
A sufficient condition for stationarity is obtained by assuming that $p$ is constant at the stationary state.
Then, additionally to the vanishing (central) divergence 
\begin{align}
 \frac12 \left( D_x -D^*_x \right)\big \rvert^{P}_{ij} u + \frac12 \left( D_y-D^*_y \right)\big \rvert^{P}_{ij} v
\end{align}
the normal derivatives (i.e. the $x$-derivative of $u$ and the $y$-derivative of $v$) need to be continuous. This brings four additional constraints: one jump across, respectively, a horizontal and a vertical edge
\begin{align}
 \left( D_x -D^*_x \right)\big \rvert^{P}_{ij}u &= 0 & \left( D_y-D^*_y \right)\big \rvert^{P}_{ij}v &= 0
\end{align}
and two jumps at the node. Together with the three pointwise conditions on the central divergence these are 7 equations for 8 variables.

The stationarity of the cell average of $p$ does not contribute an independent equation for the following reason. Assume the jumps of the normal derivatives of the velocity to vanish, i.e. the divergence of the velocity reconstruction to be continuous across the edges. Then it does not matter whether $D_x$ or $D^*_x$ is used to evaluate a derivative, as they give the same result. One can thus simply modify \eqref{eq:updatePointValuesexplicit} such that all derivatives are evaluated in cell $(i,j)$. Stationarity of $p_{ij}^P$ for all point values $P$ implies that the divergence of the velocity reconstruction vanishes at all 8 point values along the boundary of cell $(i,j)$. By Lemma \ref{lem:D22unsiolvence}, these 8 conditions for a polynomial in the 8-dimensional space $D_\text{br}^{2,2}$ imply that the polynomial itself vanishes. (This is the problem considered and solved in \cite{barsukow18activeflux} for the classical Active Flux, as has been outlined in Section \ref{sec:traditionalfourier}.) As the divergence vanishes, its average over the cell vanishes, and by Gauss' law the cell average $p_{ij}\avg$ is automatically stationary.
The kernel of the evolution matrix for upwind Active Flux is thus at least 1-dimensional.

The $12\times 12$ matrix $\mathcal{E}(\mathbf{k})$ is given in Appendix \ref{sec:2devomatricesacoustics}. That there are no other linearly independent elements in the kernel has been verified using \textsc{mathematica} due to the excessive length of computations.
\end{proof}

One observes that upwind Active Flux operates on the divergence of the velocity reconstruction in each cell, without projection (as was the case for central Active Flux earlier).

\begin{corollary}
The numerical stationary states \eqref{eq:nullspaceAcoustics2d} of semi-discrete Active Flux on 2-d Cartesian grids are the same as the stationary states \eqref{eq:stationaryeigenvector2dac} of classical Active Flux on 2-d Cartesian grids, see Section \ref{sec:traditionalfourier}. \label{remark:statStatesFullyVsSemiDiscrete}
\end{corollary}

Thus, in terms of stationarity preservation there is no difference between the behaviour of third-order semi-discrete (generalized) or fully discrete (classical) Active Flux. With the above reasoning in mind, this is not surprising at all, and is mostly a consequence of the choice of the reconstruction space.

A final comment is due concerning the rank defect mentioned in Section \ref{sec:traditionalfourier}. As is clear from the above discussion, a vanishing element in $D_\text{br}^{2,2}$ amounts to 8 equations. However, here we do not deal with \emph{any} kind of element, but with an element obtained as the divergence of an element in $(P^{2,2})^2$. The map from the velocity variables to an element of $D_\text{br}^{2,2}$, that is induced by taking the divergence of the reconstruction, is not surjective. In the proof of Theorem \ref{thm:kernel2dacoustics} it is shown that stationarity of all point values amounts to only 7 equations, and that the stationarity of the average turns out to be redundant.

The very same argument does not hold true for the central Active Flux method: the divergence can be understood as being projected onto $P^{2,2}$, which vanishes upon stationarity, but the average update continues to involve the cell average of the ``true'' divergence in $D_\text{br}^{2,2}$. This update equation for the cell average is not redundant then.

The choice of upwinding matters for stationarity preservation as well. The semi-discrete Active Flux method with the upwind splitting only involves normal derivatives such as $\del_x^2 u$ and $\del_y^2 v$. The alternative point value update \eqref{eq:alternativePointValueUpdate} also was investigated for the acoustic equations in two spatial dimensions. Clearly, this upwinding generates terms such as $\del_x^2 v$ and $\del_y^2 u$ in the numerical diffusion. In this case stationarity implies a larger number of supplementary conditions, further restricting the set of stationary states. Indeed, non-trivial stationary states no longer exist as, in general, $\det\mathcal{E}(\mathbf{k})\neq0$. For $t_x=t_y=-1$ ($k_x \Delta x = k_y \Delta y = \pi$), for example, it is
\begin{equation}
\det\mathcal{E}(\mathbf{k})=\frac{110592 c^6 (\Delta x + \Delta y)^4 (\Delta x \Delta y + 
   2 c^2 (\Delta x^2 - \Delta x \Delta y + \Delta y^2)))}{\Delta x^9 \Delta y^9} \neq 0.
\end{equation}

\subsubsection{Acoustic equations in three spatial dimensions}

\begin{theorem}
 \begin{enumerate}[(i)]
 \item The semi-discrete Active Flux method for the linear acoustic equations \eqref{acousticEquationsOhneMatrizen} in three spatial dimensions is stationarity preserving. The kernel of its evolution matrix, given in Appendix \ref{sec:3dkernelacoustics}, is 5-dimensional for general $\vec k$.
 \item The divergence of the reconstruction vanishes iff $\hat q$ is in the above-mentioned kernel.
 \end{enumerate}
\end{theorem}
\begin{proof}
 The argumentation is similar to the 2-dimensional case. There are $3 \cdot 8 = 24$ velocity variables
 \begin{align}
   \hat u\avg, \hat u\node, \hat u\ex, \dots, \hat v\avg, \hat v\node, \dots, \hat w\fz  .
 \end{align}
While central Active Flux amounts to 7 equations for the point values (plus 1 for the average, i.e. a kernel of dimension 16), the upwind method additionally enforces continuity of normal derivatives. These are $3  \cdot 1$ conditions on the faces, $3 \cdot 2$ conditions on the edges and $1 \cdot 3$ conditions on the node, which makes a kernel of dimension at least $24 - 7 - 3 - 6 - 3 = 5$. Despite the excessive length of the expressions under consideration ($\mathcal E$ is a $32 \times 32$ matrix, given in Appendix \ref{sec:3dkernelacoustics}) we confirm using \textsc{mathematica} that the kernel is indeed 5-dimensional. 
 
 One defines $D^{2,2,2}_\text{br}$ as the image of $(P^{2,2,2})^2$ under the (weak) divergence. As the jumps vanish, stationarity again implies that the divergence vanishes at the 26 locations of point values in every cell. A unisolvence Lemma similar to Lemma \ref{lem:D22unsiolvence} holds again, that we omit. This implies the second statement.
\end{proof}

Thus, the stationary states of this Active Flux method are those that give rise to a divergence-free reconstruction of the velocity and a uniformly constant pressure. The 5 elements $Q_1, \ldots, Q_5$ of $\ker \mathcal E$ are given in Appendix \ref{sec:3dkernelacoustics}.
We find that the projection onto just $Q_1$ is consistent with $f$ that is chosen from the subspace spanned by $(-k_y, k_x, 0, 0)$ alone and $Q_2$ is consistent with $f$ from the subspace parallel to $(-k_z, 0, k_x, 0)$, i.e. it so happens that our choice of basis of $\ker \mathcal E$ even is consistent with the choice of basis for $\vec J \cdot \vec k$ from \eqref{eq:3dker}. The remaining 3 elements of $\ker \mathcal E$ correspond to higher order terms.


\section{Analysis of Numerical Diffusion}   \label{sec:numdiff}

\subsection{Linear advection in one spatial dimension}

\begin{figure}
\centering
\includegraphics[width=\textwidth]{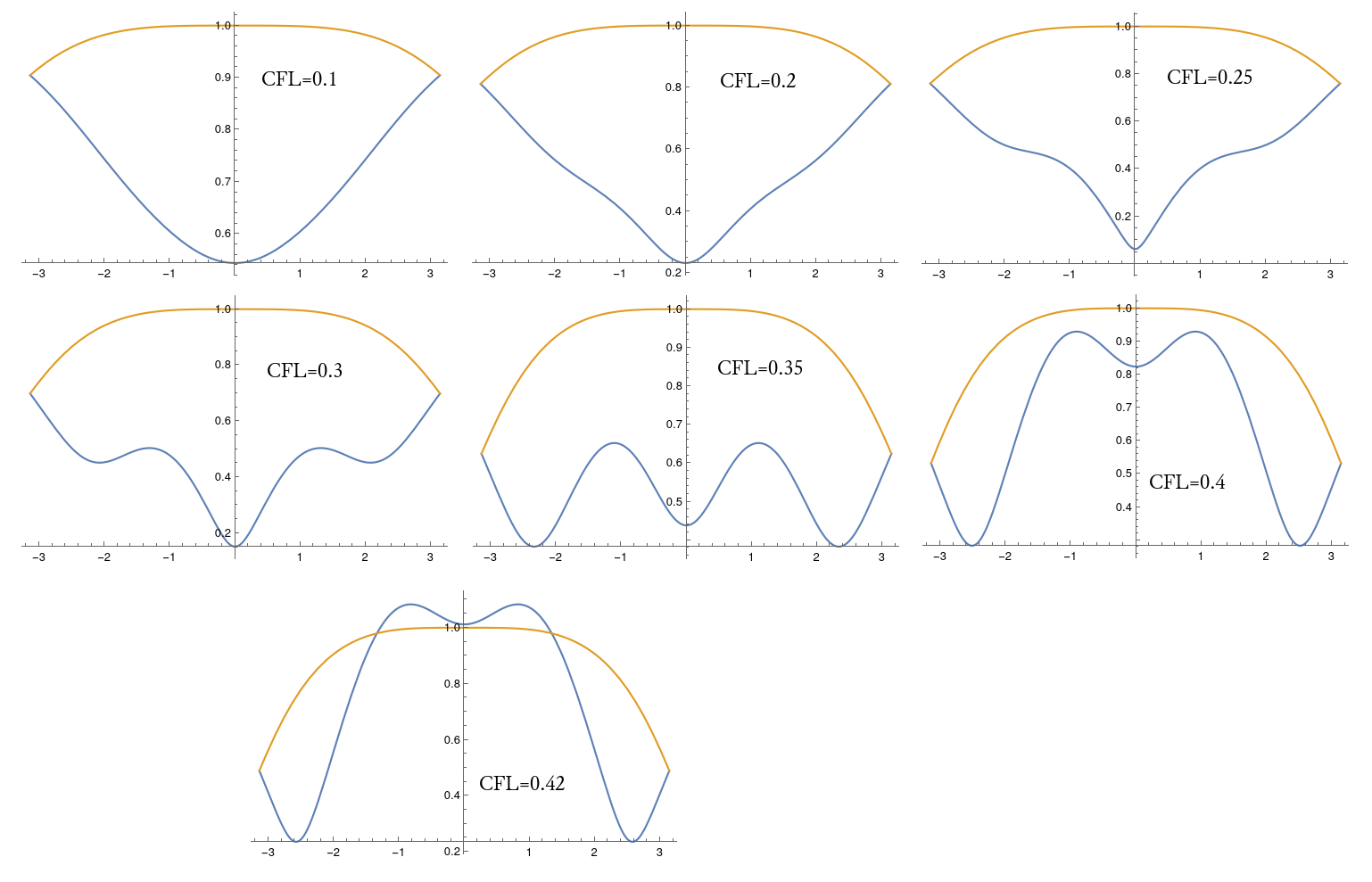}
\caption{Numerical diffusion analysis of Active Flux for 1-d linear advection. The absolute value of the two eigenvalues of the $2 \times 2$ amplification matrix is shown as a function of $\beta := k \Delta x \in [-\pi,\pi]$. Observe the different scales on the vertical axes.}
\label{fig:diffusiondispersion}
\end{figure}

An analysis of numerical diffusion for multi-dimensional acoustics is rather difficult. Therefore, we first consider the one-dimensional case of Active Flux applied to the linear advection equation $\del_t q + c \del_x q = 0$ ($c>0$). Active Flux has two degrees of freedom per cell (one average and one point value) in this case (see e.g. \cite{eymann13,kerkmann18,barsukow19activeflux,abgrall22}), and upon the Fourier transform in space it can be rewritten as
\begin{align}
 \del_t \hat q + \mathcal E(k) \hat q &= 0
\end{align}
with a $2\times 2$ evolution matrix $\mathcal E$. For a reasonable analysis of numerical diffusion, a time discretization needs to be specified, and a Runge-Kutta method of order 3 seems well-suited. Thus, the Fourier transform of the fully discrete method reads
\begin{align}
 \hat q^{n+1} &= \mathcal A(k) \hat q^n
\end{align}
with the amplification matrix
\begin{align}
 \mathcal A &= \id - \Delta t \mathcal E + \frac12 \Delta t^2 \mathcal E^2 - \frac{1}{6} \Delta t^3 \mathcal E^3.
\end{align}
Figure \ref{fig:diffusiondispersion} shows the modulus of the two eigenvalues of $\mathcal A$ as a function of $\beta := \Delta x k \in [-\pi, \pi]$ for different values of $\text{CFL} := \frac{c\Delta t}{\Delta x}$. 

The reason for having two eigenvalues is the following: The discrete Fourier transform distinguishes between the modes associated to the different types of degrees of freedom, i.e. there is a mode for the averages and an independent mode for the point values in the one-dimensional case. For example $k=0$ corresponds to the point values having all the same value $Q_\text{p} \in \mathbb R$ and the averages having all the same value $Q_\text{a} \in \mathbb R$, but the two values can be different in general. The point update will see this and generally not remain stationary: while the ``physical'' eigenvalue corresponds to the eigenvector $(1,1)^\text{T}$, which means $Q_\text{a} = Q_\text{p}$, the other eigenvalue corresponds to the eigenvector $(0,1)^\text{T}$, i.e. $Q_\text{a} = 0$, $Q_\text{p} = 1$. This highly oscillatory function evolves in time, decaying for small CFL numbers and exploding for those high enough, as is discussed next.

As linear advection amounts merely to a translation, the exact amplification factor is a unit complex number. Von Neumann stability therefore requires all the eigenvalues of the amplification matrix to not be greater than 1 in absolute value. Of the two eigenvalues shown in Figure \ref{fig:diffusiondispersion}, one observes that the ``physical'' eigenvalue remains, in absolute value, below and rather close to 1. 

It has been elucidated in \cite{roe21}, that for semi-discrete methods stability quite generally tends to be governed by the ``non-physical'' eigenvalue. Also in the case of semi-discrete one-dimensional Active Flux, one can observe it grow as the CFL number is chosen larger and larger, until for $\text{CFL} = 0.42$ its absolute value surpasses 1. This is consistent with previous findings (\cite{zeng19,abgrall22proceeding}) of 0.41 as the stability limit for semi-discrete Active Flux in one spatial dimension.

\subsection{Acoustic equations in two spatial dimensions}

\begin{figure}
\centering
\includegraphics[width=\textwidth]{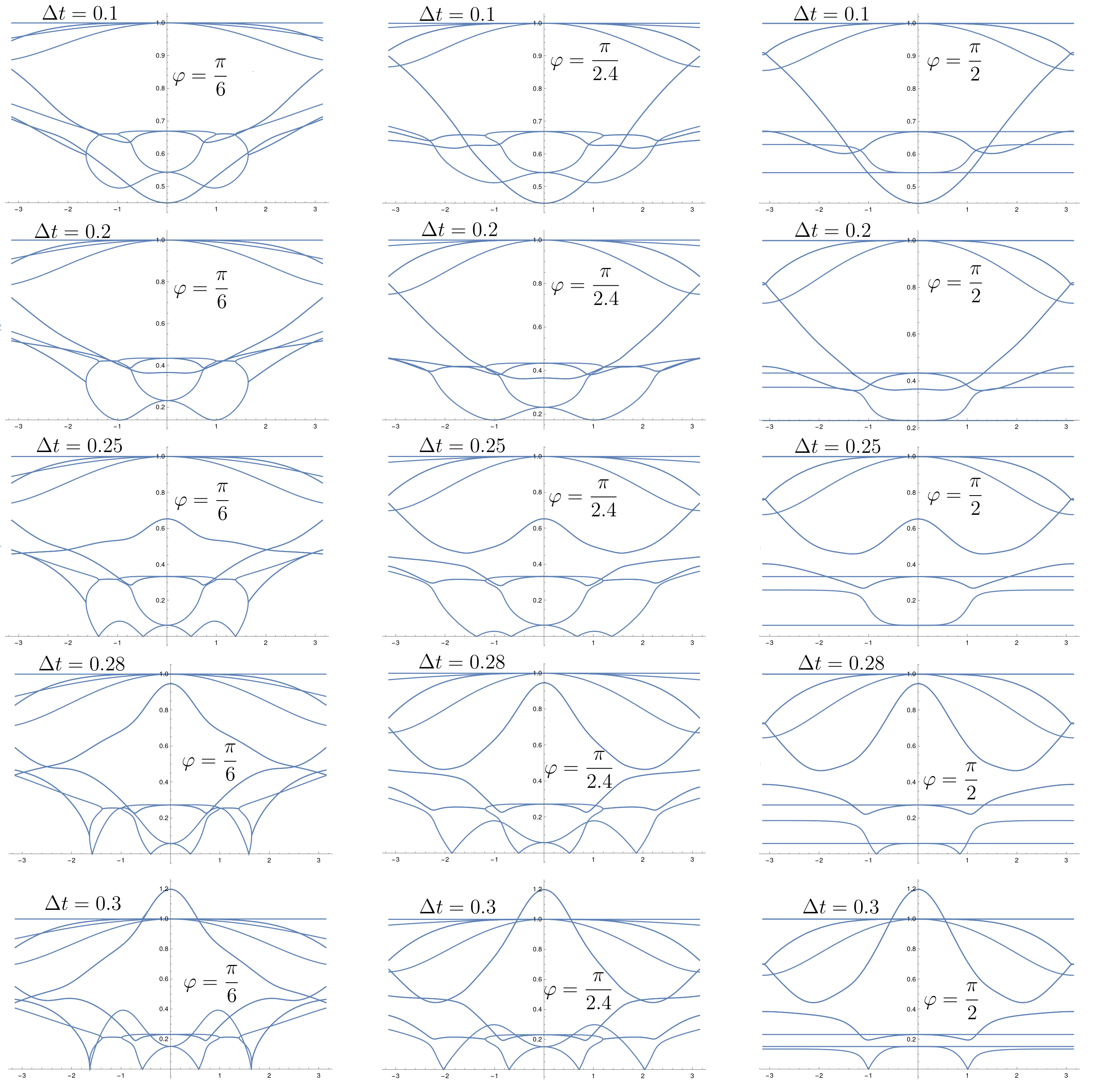}
\caption{Numerical diffusion analysis of Active Flux for 2-d linear acoustics. The absolute value of the 12 eigenvalues of the amplification matrix is shown as a function of $s\in [-\pi,\pi]$ (see \eqref{eq:kparametrization}) for different values of $\Delta t$ (\emph{rows}) and $\phi$ (\emph{columns}). Here, $\Delta x = \Delta y = 1$. }
\label{fig:diffusionacoustic}
\end{figure}

Finally, to complete the picture, the same analysis is performed for the two-dimensional Active Flux method for linear acoustics. The absolute values of the 12 eigenvalues are shown in Figure \ref{fig:diffusionacoustic} (some of them lie on top of each other). For simpler presentation, here $\Delta x = \Delta y = 1$ and the wave vector $\vec k$ is parametrized as 
\begin{align}
 \vec k = s \vecc{\cos \phi}{\sin \phi} \label{eq:kparametrization}
\end{align}
To perform a comparison between the behaviour of the numerics and that of the solutions to the PDE, one now needs to know the analytical value of the amplification matrix. From \eqref{eq:Jk}, writing $\Lambda := R^{-1} (\vec J \cdot \vec k) R$, one obtains
\begin{align}
\hat q(t + \Delta t) &= R \exp( - \ii \Lambda t) R^{-1}\hat q(t)
\end{align}
The eigenvalues $\Lambda$ of $(\vec J \cdot \vec k)$ are real by hyperbolicity of \eqref{eq:conservationLaw} (they are $0$ and $\pm c\rvert\vec k\rvert$). Thus, the eigenvalues of the analytical amplification matrix $R e^{- \ii \Lambda t} R^{-1}$ are all of modulus one.

Despite the increased complexity, the behaviour is similar to that of linear advection. The instability is governed by the non-physical eigenvalue(s), whose norm becomes larger than 1 between $\Delta t = 0.28$ and $\Delta t = 0.3$. We have not observed a strong $\phi$-dependence of this value.

\section{Numerical examples} \label{sec:numerical}
\subsection{Acoustic Equations in 2-d}

\subsubsection{Traveling waves} \label{sec:convergence2d}

\begin{figure}
 \centering
 \includegraphics[width=0.32\textwidth]{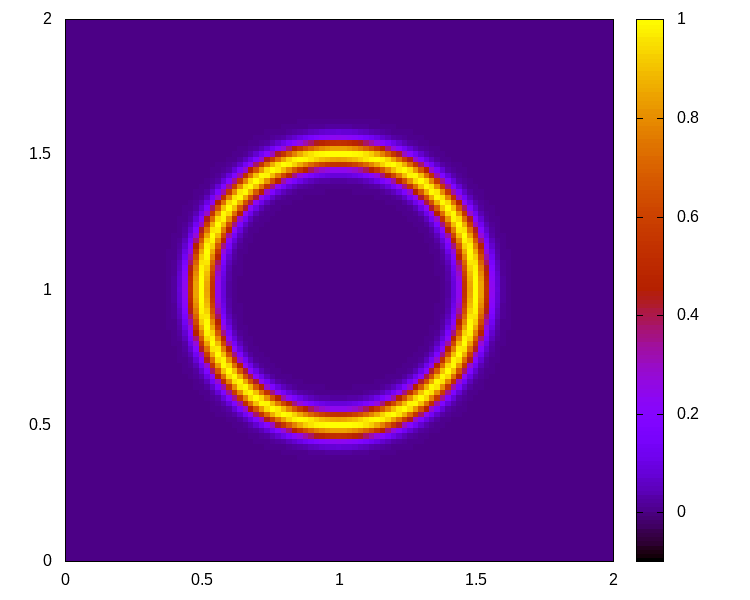} \includegraphics[width=0.32\textwidth]{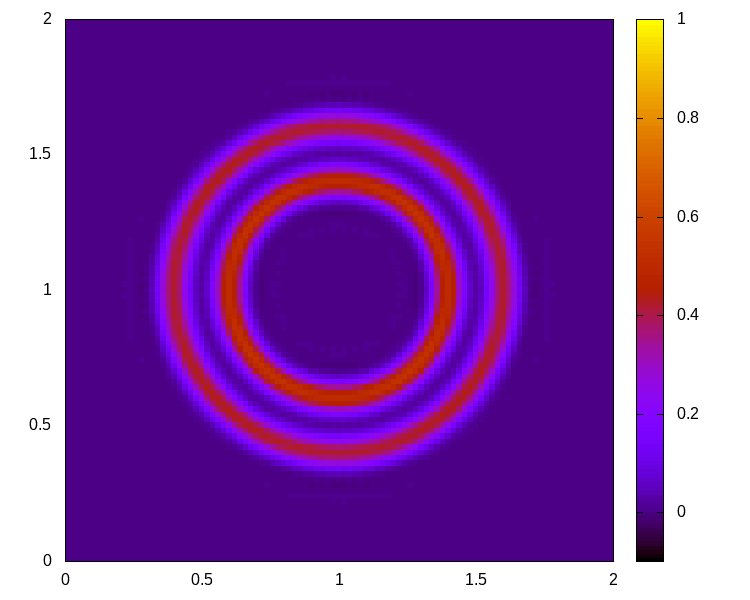} \includegraphics[width=0.32\textwidth]{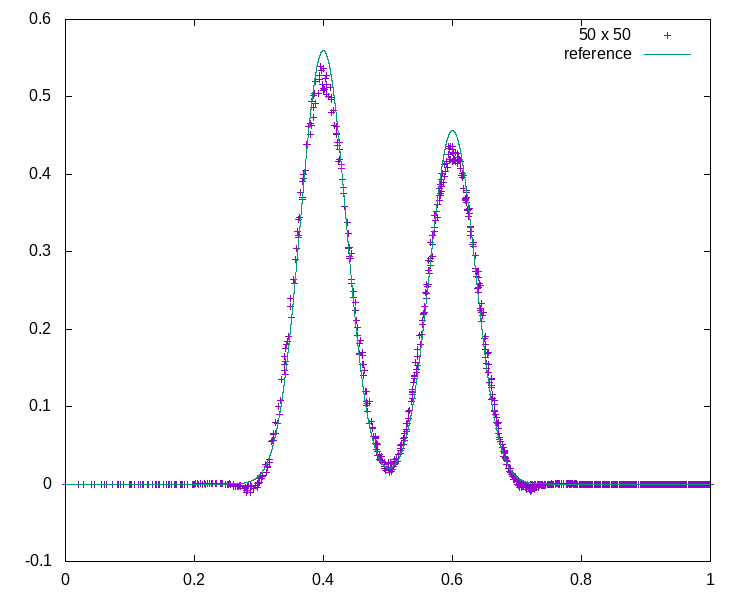}
 \caption{The setup of a spherical Gaussian in the pressure. \emph{Left}: Initial setup. \emph{Center}: Pressure at time $t=0.1$ solved on a $50 \times 50$ grid. \emph{Right}: Scatter plot of the pressure at time $t=0.1$ together with a reference solution (solid line).}
 \label{fig:gaussianspherical}
\end{figure}

\begin{figure}
 \centering
 \includegraphics[width=0.32\textwidth]{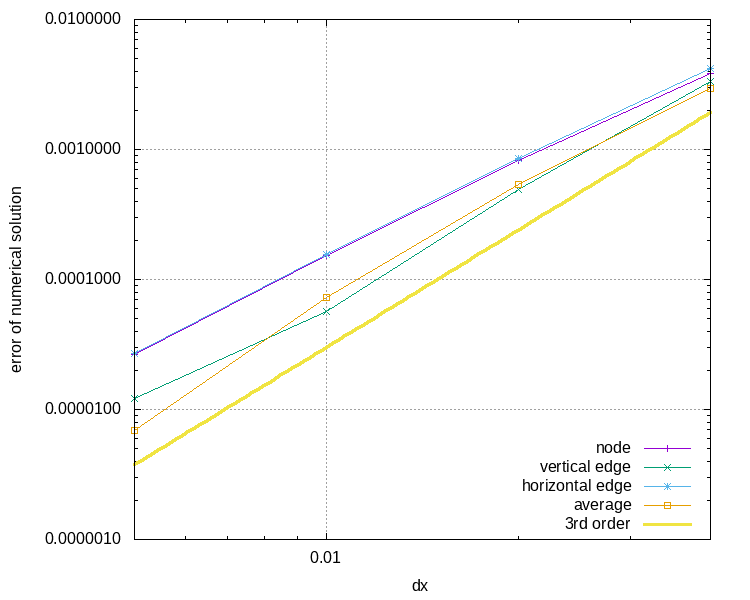} \includegraphics[width=0.32\textwidth]{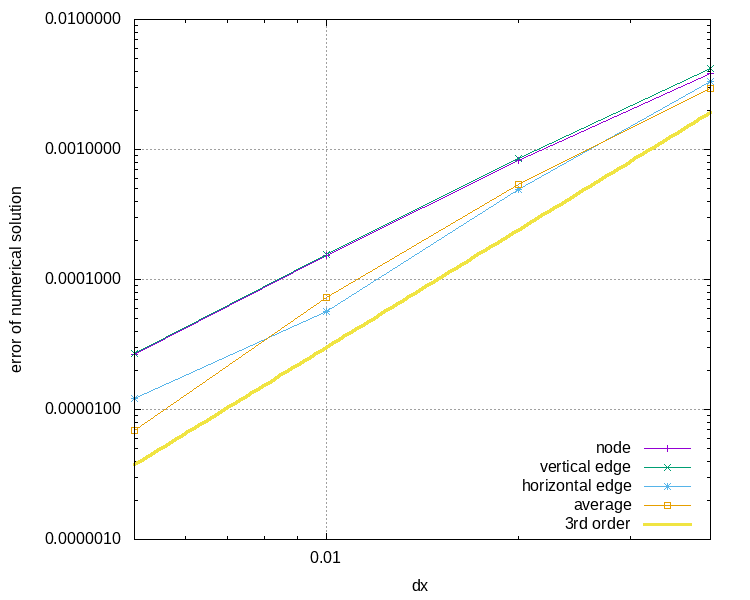} \includegraphics[width=0.32\textwidth]{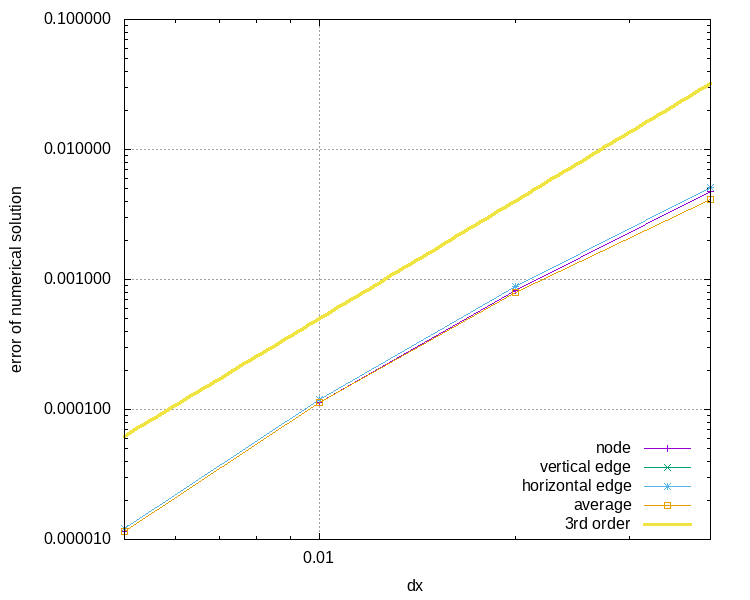}
 \caption{The error of the numerical solution for the spherical Gaussian wave at time $t = 0.1$ upon grid refinement. \emph{From left to right}: $u$, $v$, $p$.}
 \label{fig:gaussiansphericalerror}
\end{figure}

We consider a truly multi-dimensional setup of a spherical Gaussian in the pressure:
\begin{align}
 p_0(\vec x) &= \exp\left( \frac{ \sqrt{(x-1)^2 + (y-1)^2} - r_0)^2}{w^2} \right) & \vec v_0(\vec x) &= 0
\end{align}
with $w = 0.05$, $r_0 = \frac12$. The setup is shown in Figure \ref{fig:gaussianspherical}. It is solved on a domain of $[0,2]^2$ using RK3 with a CFL number of 0.2 with periodic boundaries. The $L^1$ error of the numerical solution is shown in Figure \ref{fig:gaussiansphericalerror} as a function of the spatial discretization length. The reference solution has been obtained with a first-order method solving the radial equations on a grid of $4 \cdot 10^5$ points. We observe third order accuracy.


\subsubsection{Stationary vortex}

We consider a stationary solution of the acoustic equations:
\begin{align}
 p_0(\vec x ) &= 0 & \vec v_0(\vec x) &= \vecc{-(y-\frac12)/r}{(x-\frac12)/r} \cdot \begin{cases}   
                                          5 r & r < 0.2 \\
                                          \max(0, 2 - 5r) & \text{else}
                                         \end{cases}
\end{align}
with $r = \sqrt{\left( x-\frac12\right )^2  + \left(y-\frac12 \right )^2}$.
This setup is solved using the semi-discrete Active Flux method (3rd order) on a grid of $50 \times 50$, using RK3 and a CFL = 0.2, with zero-gradient boundary conditions. 

\begin{figure}
 \centering
 \includegraphics[width=\textwidth]{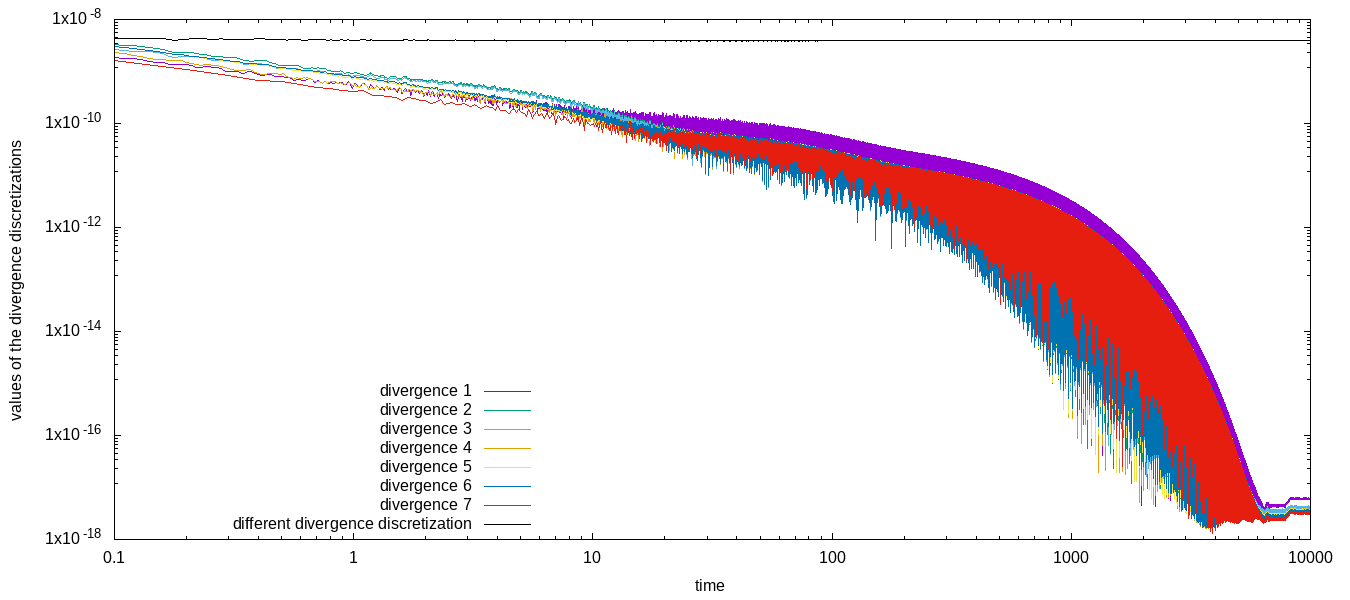}
 \caption{Stationary vortex test case. Decay of the 7 discrete divergences \eqref{eq:discreteStatStatesVelocity1}--\eqref{eq:discreteStatStatesVelocity4} that characterize the stationary state. Their values reach machine precision after long times. For comparison, the behaviour of some other discrete divergence is shown; it does not decay to machine zero. The curves lie partly on top of each other, but we refrain from showing details as the Figure is merely intended to show that all the divergences reach machine zero, and not how.}
 \label{fig:vortexstatdiv}
\end{figure}

The aim of this test is to show experimentally that a discrete steady state is indeed parallel to \eqref{fig:vortexstatdiv}, and we deliberately choose a setup involving many different Fourier modes. First, one needs to wait until the setup becomes stationary. Then, instead of computing the Fourier transform, we compute the values of the discrete divergences \eqref{eq:fourierVelAverage}--\eqref{eq:fourierdiv7}.

\begin{figure}
 \centering
 \includegraphics[width=0.45\textwidth]{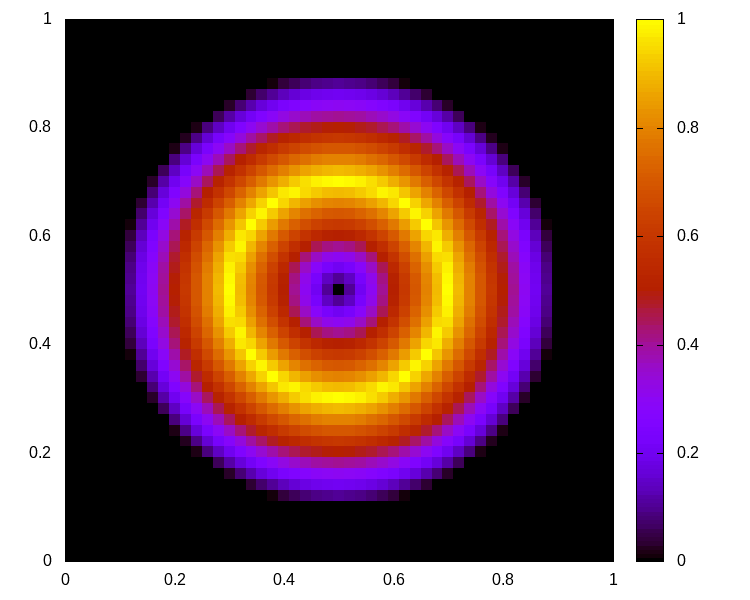} \hfill \includegraphics[width=0.45\textwidth]{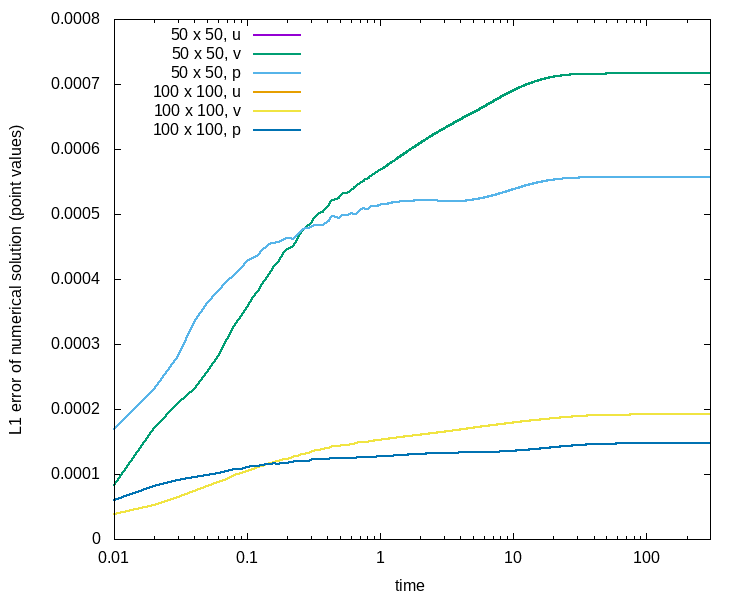}
 \caption{\emph{Left}: Initial data of the stationary vortex test. The magnitude of the velocity is color-coded. \emph{Right}: Error of the numerical solution as a function of time. One observes the stationarization of the setup (the velocity components are on top of each other).}
 \label{fig:vortexstat}
\end{figure}

Figure \ref{fig:vortexstat} demonstrates that the setup becomes stationary in all variables at about $t = 100$. Unsurprisingly this happens later for finer grids (less numerical diffusion) and the stationary state is also closer to the initial one. 
Figure \ref{fig:vortexstatdiv} shows the time evolution of the discrete divergences \eqref{eq:fourierVelAverage}--\eqref{eq:fourierdiv7} (\eqref{eq:discreteStatStatesVelocity1}--\eqref{eq:discreteStatStatesVelocity4}) on a grid of $50 \times 50$. As our datum is not well-prepared, they are not zero initially, but as the setup becomes stationary, they attain values at the level of machine precision. This indicates that $\ker \mathcal E$ is indeed one-dimensional and parallel to \eqref{fig:vortexstatdiv}. For comparison, the Figure shows the time evolution of a different divergence discretization
\begin{align}
 \frac{1+2t_y+t_y^2}{t_y}\frac{(-1+t_x)(1+t_x)}{t_x\Delta x}\hat{u}\avg+\frac{1+2t_x+t_x^2}{t_x}\frac{(-1+t_y)(1+t_y)}{t_y\Delta y}\hat{v}\avg
\end{align}
that corresponds to
\begin{align}
\frac{\langle[u\avg]_{i\pm1}\rangle_{j}^{(2)}}{\Delta x}+\frac{[\langle v\avg\rangle_{i}^{(2)}]_{j\pm1}}{\Delta y}
\end{align}
which can be observed to become stationary (as the entire setup) but not to decay to machine zero.


\begin{figure}
 \centering
 \includegraphics[width=0.45\textwidth]{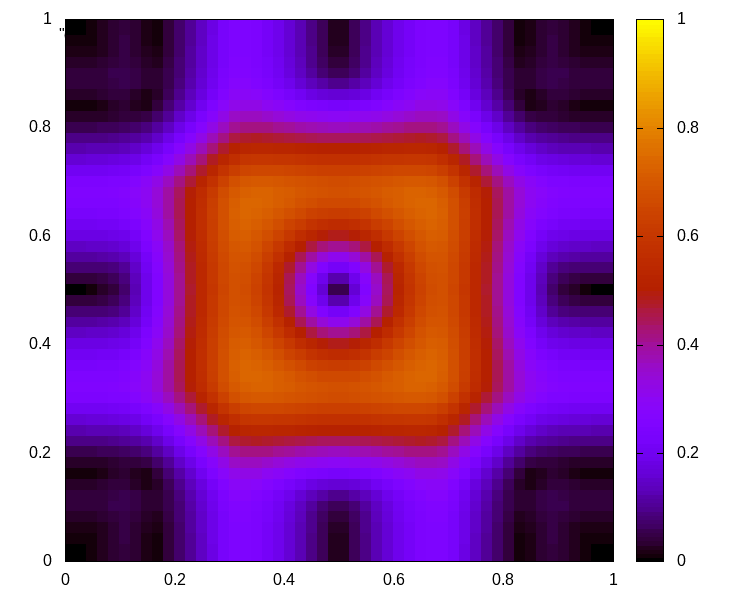} \hfill \includegraphics[width=0.45\textwidth]{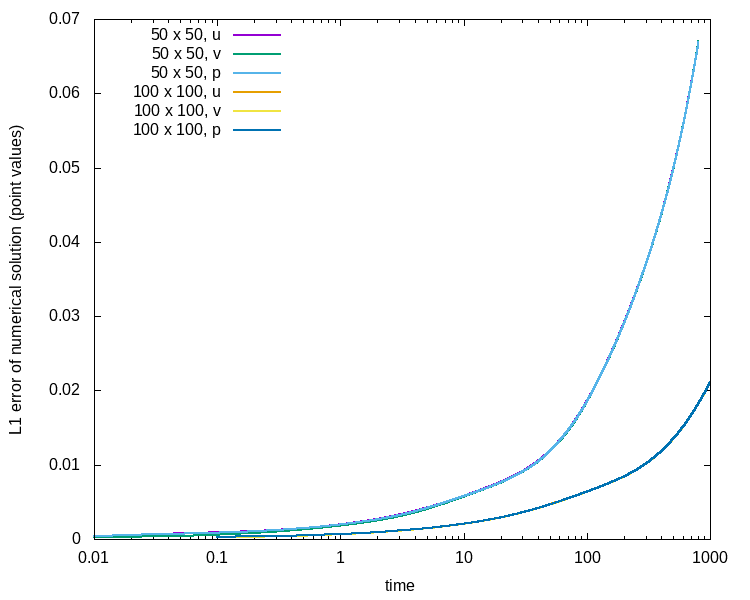}
 \caption{Rusanov-type Jacobian split \eqref{eq:alternativePointValueUpdate} does not lead to a stationarity preserving method. \emph{Left}: Solution of the stationary vortex test at $t=1000$, the magnitude of the velocity is color-coded. The vortex has become two overlapping shear flows. \emph{Right}: Error of the numerical solution as a function of time. It takes much longer than shown for the setup to stationarize, and the final state is very far from the initial condition (compare the errors to those in Fig. \ref{fig:vortexstat}).}
 \label{fig:vortexstatrusanov}
\end{figure}

\subsubsection{Well-prepared stationary mode}

\begin{figure}
 \centering
 \includegraphics[width=0.49\textwidth]{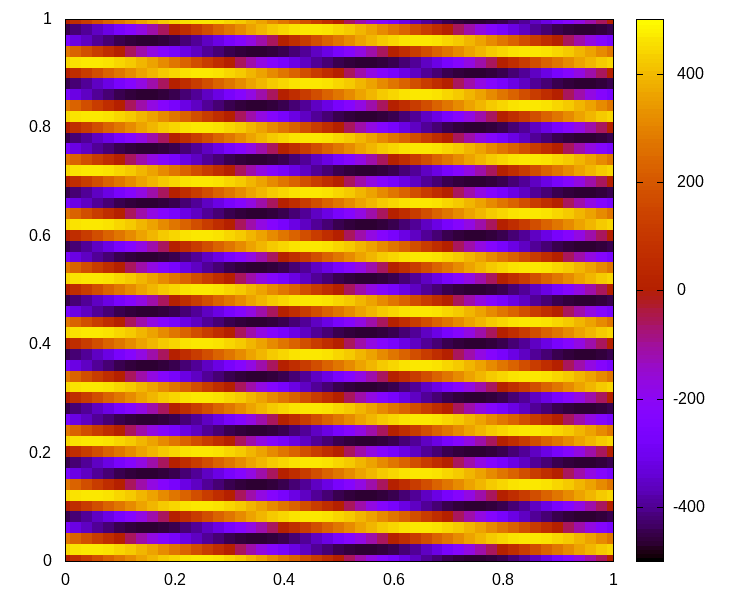} \hfill \includegraphics[width=0.49\textwidth]{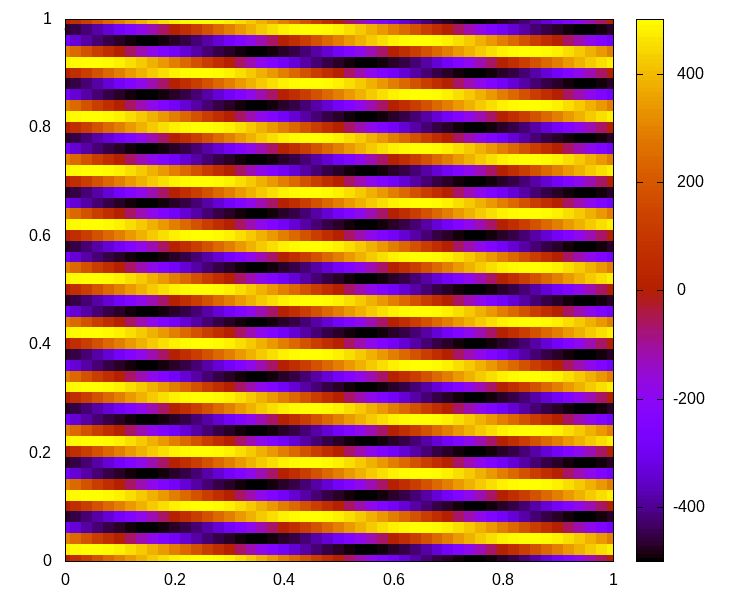}\\
 \includegraphics[width=0.49\textwidth]{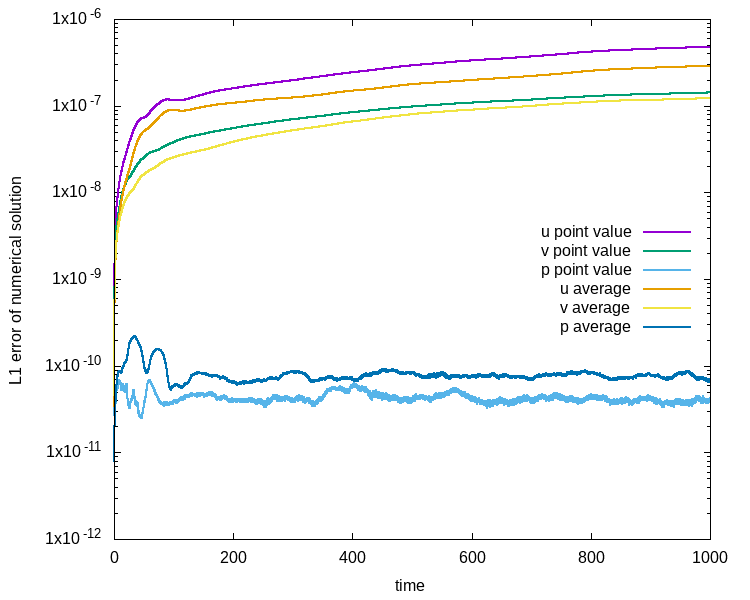} \hfill  \includegraphics[width=0.49\textwidth]{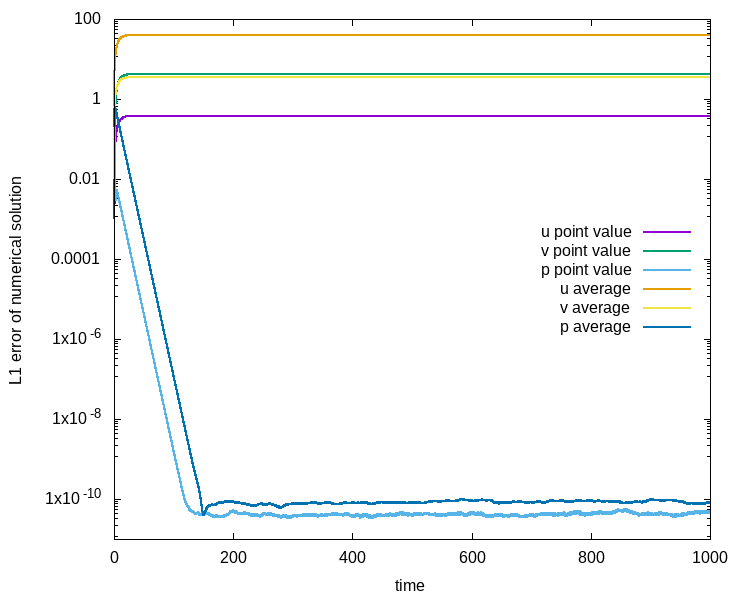}
 \caption{Stationary mode for 2-d acoustics. \emph{Top}. Initial data in $u$, nodal values are shown. \emph{Left}: Well-prepared data \eqref{eq:discretestationaryrealmode1}. \emph{Right}: Not well-prepared data directly from \eqref{eq:mode5}. There is no visible difference. \emph{Bottom}: Error of the numerical solution is shown as a function of time. Note the different scales of the vertical axis. \emph{Left}: Well-prepared data show only linear growth due to round-off errors; the pressure errors are very low. \emph{Right}: Not discretely well-prepared data. The solution rapidly deviates from the initial state. It stationarizes after an initial layer at around $t\simeq 150$. The pressure, initially zero, attains finite values in the first step of the calculation and decays exponentially quickly towards zero again.}
 \label{fig:statmodeacoustic}
\end{figure}

The next aim is to verify the preservation of the discrete stationary state given by \eqref{eq:nullspaceAcoustics2d}. To this end, the discrete data are well-prepared in a way similar to \cite{barsukow20hypproceeding}, starting from the Fourier mode $\hat Q \exp(\ii k_x i \Delta x + \ii k_y j \Delta y)$. For it to be stationary, $\hat Q$ has to be parallel to \eqref{eq:nullspaceAcoustics2d}. This would result in complex-valued data. Therefore, we take the average of two such modes with opposite signs of $\vec k$, which results in a real-valued grid function:
\begin{subequations}
\begin{align}
 q_{ij}\avg &= \left( \begin{array}{c} 
\displaystyle \frac{8 (2+\cos (2 \Delta x  \pi )) \sin (20 \Delta y  \pi ) \sin (2 \pi  (i \Delta x+10  j \Delta y))}{3 \Delta y }\\
\displaystyle -\frac{8 (2+\cos (20 \Delta y  \pi )) \sin (2 \Delta x  \pi ) \sin (2 \pi  (i \Delta x+10  j \Delta y))}{3 \Delta x }\\
\displaystyle 0 
\end{array} \right ) \\
q_{ij}\eh &=  \left( \begin{array}{c} 
\displaystyle \frac{4\sin (10 \Delta y  \pi ) \sin (\pi (2 i \Delta x + 20  j \Delta y + 10 \Delta y )) (3+\cos (2 \Delta x  \pi ))}{\Delta y }\\
\displaystyle -\frac{8 \cos (10 \Delta y  \pi ) \sin (2 \Delta x  \pi ) \sin (2 \pi  (5 \Delta y +i \Delta x+10  j \Delta y))}{\Delta x }\\
\displaystyle 0
\end{array} \right )\\
q_{ij}\ev &=  \left( \begin{array}{c} 
\displaystyle \frac{8 \cos (\Delta x  \pi ) \sin (20 \Delta y  \pi ) \sin (\pi  (\Delta x +2 i \Delta x+20  j \Delta y))}{\Delta y }\\
\displaystyle -\frac{4 (3+\cos (20 \Delta y  \pi )) \sin (\Delta x  \pi ) \sin (\pi  (\Delta x +2 i \Delta x+20  j \Delta y))}{\Delta x }\\
\displaystyle 0
\end{array} \right )\\
q_{ij}\node &=  \left( \begin{array}{c} 
\displaystyle \frac{16 \cos (\Delta x  \pi ) \sin (10 \Delta y  \pi ) \sin (2 \pi (\Delta x /2 + 10  j \Delta y + i \Delta x + 5 \Delta y ))}{\Delta y }\\
\displaystyle -\frac{16 \sin (\Delta x  \pi ) \cos (10 \Delta y  \pi ) \sin (2 \pi (\Delta x /2 + 10 j \Delta y + i \Delta x + 5 \Delta y ))}{\Delta x }\\
\displaystyle 0
\end{array} \right ) 
\end{align} \label{eq:discretestationaryrealmode1}
\end{subequations}
The wave numbers are chosen as in \cite{barsukow20hypproceeding}: $k_x = 2\pi$ and $k_y = 10 \cdot 2 \pi$. The setup is solved 
using CFL = 0.2 on a $50 \times 50$ grid using periodic boundary conditions. The results are shown in Figure \ref{fig:statmodeacoustic}, and one observes that the well-prepared setup remains stationary, up to round-off errors.

Setting $i \Delta x = x$ and $j \Delta y = y$, as $\Delta x, \Delta y \to 0$, this mode can be seen to converge to
\begin{align}
 q_{ij}\avg &\to \tilde q(x, y ) \label{eq:mode1} & 
 q_{ij}\eh &\to \tilde q\left(x, y + \frac{\Delta y}{2}\right) \\
 q_{ij}\eh &\to \tilde q\left(x, y + \frac{\Delta y}{2}\right)  &
 q_{ij}\node &\to \tilde q\left(x + \frac{\Delta x}{2}, y + \frac{\Delta y}{2}\right)
\end{align}
with 
\begin{align}
 \tilde q(x, y) &= \veccc{10}{-1}{0} 16 \pi \sin(2 \pi (x + 10 y)). \label{eq:mode5}
\end{align}

For comparison, one might therefore consider discrete initial data that have been obtained by the direct evaluation of the degrees of freedom from \eqref{eq:mode5}, i.e. by taking point values and cell averages of \eqref{eq:mode5}. In this case one observes an initial layer ($t \leq 150$) during which the pressure decays after having attained finite-size values in the first step of the calculation. The seemingly high value $\sim 40$ of the error of $u$ is due to generally high values of this component: the maximum exact value of $u$ is $160 \pi \simeq 503$, i.e. the error actually is only $8\%$.

\subsection{Acoustic Equations in 3-d}

\subsubsection{Traveling waves}

A test analogous to that of Section \ref{sec:convergence2d} is performed to assess the order of accuracy of the three-dimensional method. The initial data vanish for all velocity components, and the initial pressure is given by
\begin{align}
 p_0(\vec x) &= \exp\left( \frac{(\lvert\vec x\rvert - r_0)^2}{w^2}  \right).
\end{align}

\begin{figure}
 \centering
 \includegraphics[width=0.7\textwidth]{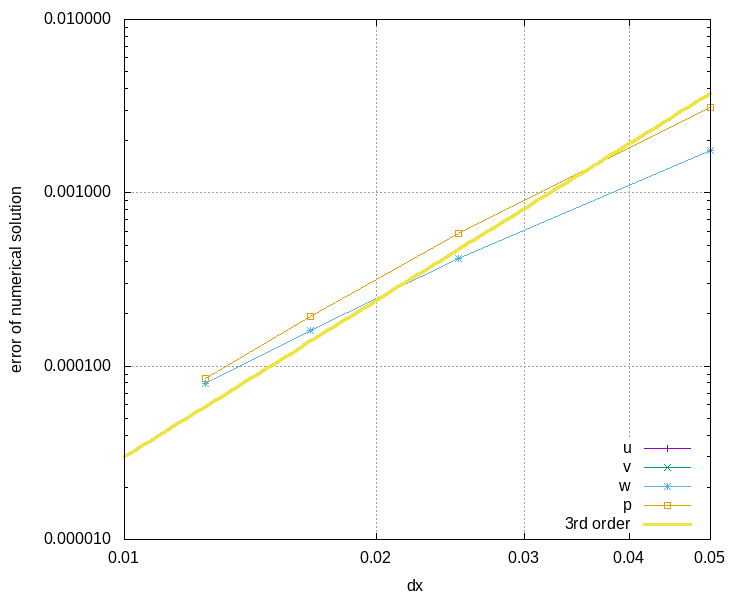}
 \caption{Convergence study for the spherical Gaussian in the pressure. Computations on finer meshes are too expensive to run, but the data shown on the coarse meshes indicate a convergence order of about 2.5. The error curves of the three velocities are on top of each other.}
 \label{fig:convergence3d}
\end{figure}

Again, $r_0 = \frac12$ and $w = 0.05$, and the domain is $[-1,1]^3$. We use RK3 with a CFL number of 0.1. The $L^1$ numerical error as a function of the spatial discretization length is shown in Figure \ref{fig:convergence3d}, which was computed using a reference solution obtained from a 1-d spherically symmetric code on $4 \cdot 10^6$ points. The three-dimensional computations are rather expensive and do not allow to increase the number of cells as much as one requires to see third-order accuracy clearly. The data shown indicate a convergence rate of about 2.5.

\subsubsection{Stationary mode} \label{sec:stationarymode3d}

\begin{figure}
 \centering
 \includegraphics[width=0.49\textwidth]{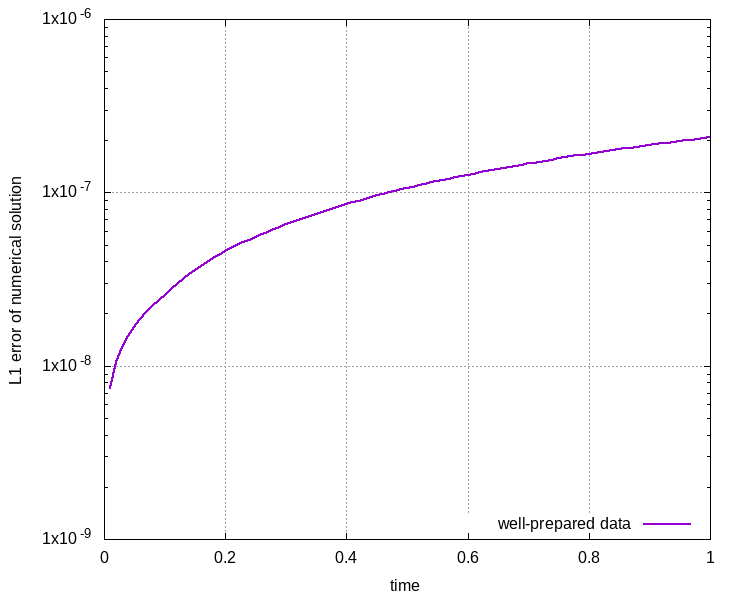} \includegraphics[width=0.49\textwidth]{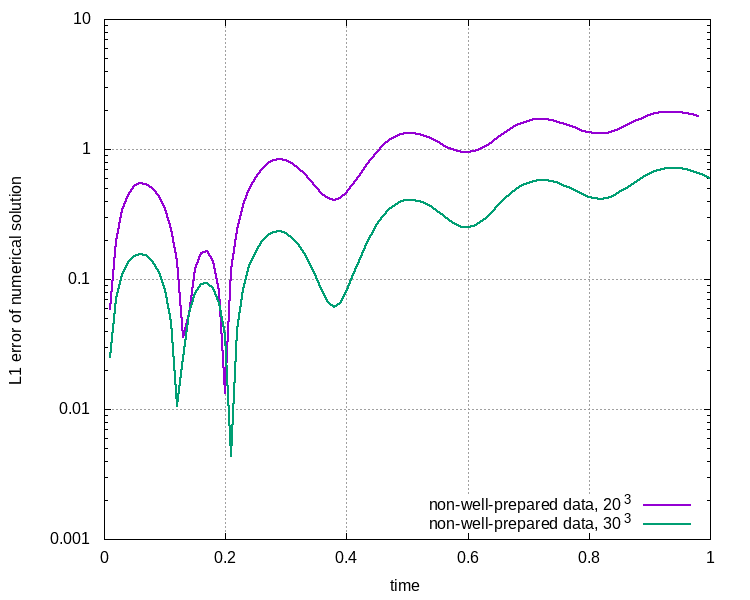}
 \caption{The $L^1$ distance between the initial data of a Fourier mode in the kernel of the evolution matrix $\mathcal E$ for Active Flux in three spatial dimensions shown as function of time. \emph{Left}: Non-well-prepared initialization according to \eqref{eq:convergeed3dmode}. The mode is not discretely stationary, but upon grid refinement the error decreases at approximately third order. \emph{Right}: Well-prepared data according to \eqref{eq:statmode3d}. One observes only linear growth due to round-off errors; note also the different value of the error.}
 \label{fig:fouriermodeerror}
\end{figure}

We consider a linear combination of the elements of the kernel, given in Appendix \ref{sec:3dkernelacoustics}
\begin{align}
 Q_{ijk} = \sum_{r = 1}^5 a_r \hat Q_r \exp(\ii k_x i \Delta x + \ii k_y j \Delta y + \ii k_z k \Delta z) \label{eq:statmode3d}
\end{align}
with $a_i$ also given in Appendix \ref{sec:3dkernelacoustics}. The reason for the coefficients is that a Finite Difference formula corresponds to a Laurent polynomial in $t_x,t_y,t_z$, and all the denominators appearing in the basis elements $\hat Q_i$ need to be removed first. The domain is $[0,1]^3$ with periodic boundaries, discretized by $20^3$ cells and 
\begin{align}
 k_x &= 2 \pi & k_y &= 8 \pi & k_z &= -4 \pi
\end{align}
We evolve only its real part, which is possible due to the linearity of the entire problem. Figure \ref{fig:fouriermodeerror} shows that indeed, the error is initially at machine precision and grows only linearly due to round-off errors.

To highest order in $\Delta x, \Delta y, \Delta z$, the values at the nodes are
\begin{align}
 (2 k_z (k_y - \ii), -2 k_x k_z, 2 \ii k_x, 0)^\text{T} \exp(\ii k_x x + \ii k_y y + \ii k_z z) \label{eq:convergeed3dmode}
\end{align}
One easily verifies that this is a divergence-free Fourier mode. However, if this mode is used as initial data simply by pointwise evaluation, it is not discretely stationary.

\subsubsection{Vortex ring}

\begin{figure}
 \centering
 \includegraphics[width=0.7\textwidth]{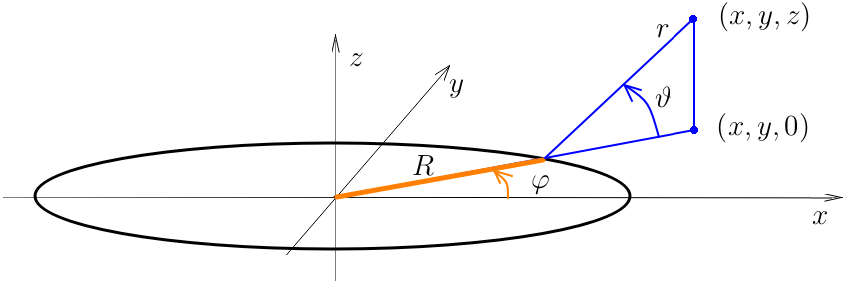}
 \caption{Coordinate setup for the vortex ring test case.}
 \label{fig:vortex-coordinates}
\end{figure}

We consider a vortex ring centred at the origin, with its centerline a circle of radius $R$ in the $x$-$y$-plane (see Figure \ref{fig:vortex-coordinates}).

Define first 
\begin{align}
 r := \sqrt{\left(\sqrt{x^2 + y^2} - R\right )^2 + z^2}
\end{align}
and angles in the plane and perpendicular to it:
\begin{align}
 \tan \phi &= \frac{y}{x}  &
 \sin\theta &= \frac{z}{r}
\end{align}

The flow happens only in the plane spanned by 
\begin{align}
 \veccc{0}{0}{1} \quad \text{and} \quad \veccc{\cos \phi}{\sin \phi}{0}
\end{align}

\begin{figure}
 \centering
 \includegraphics[width=0.49\textwidth]{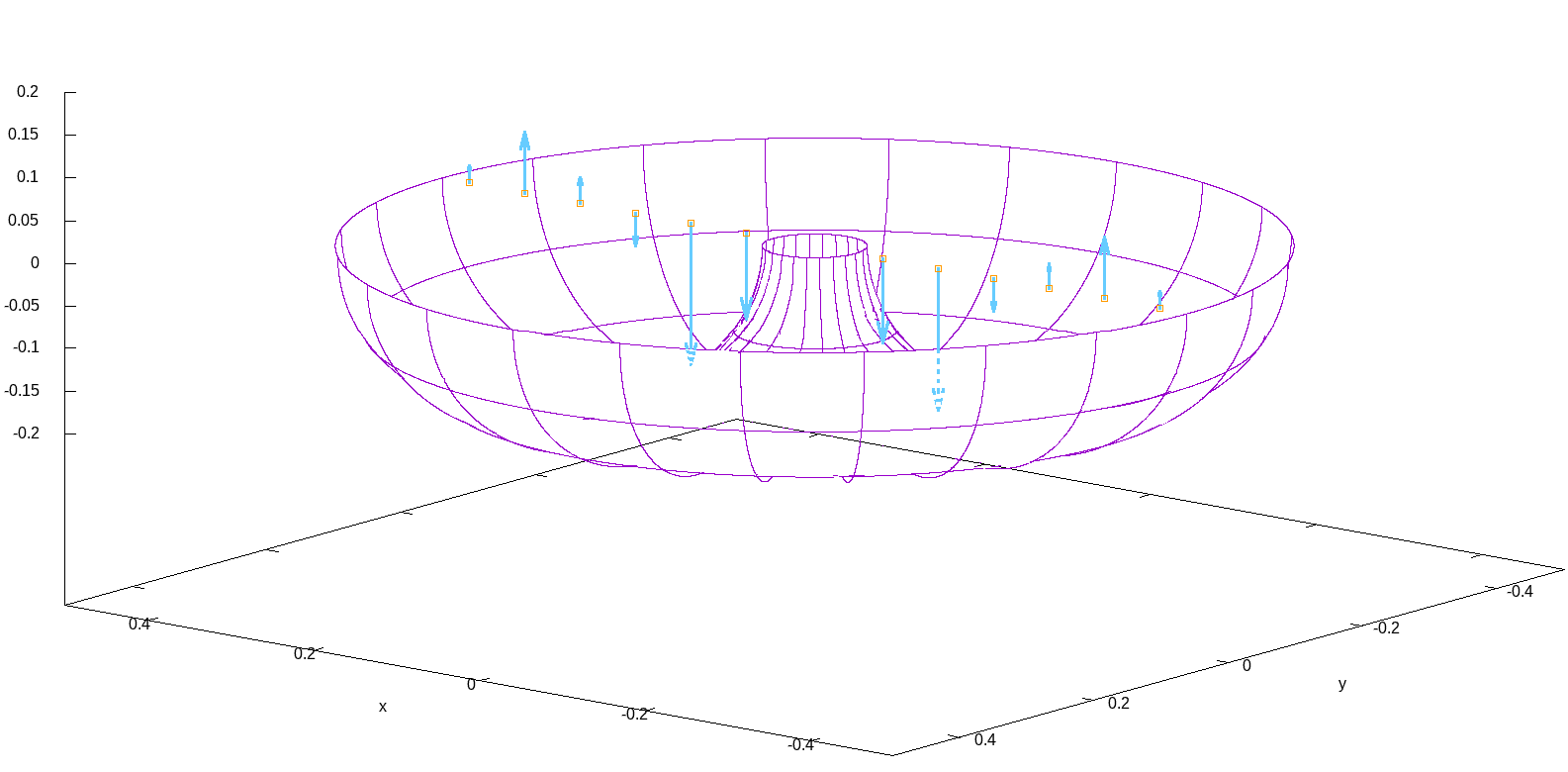} \includegraphics[width=0.49\textwidth]{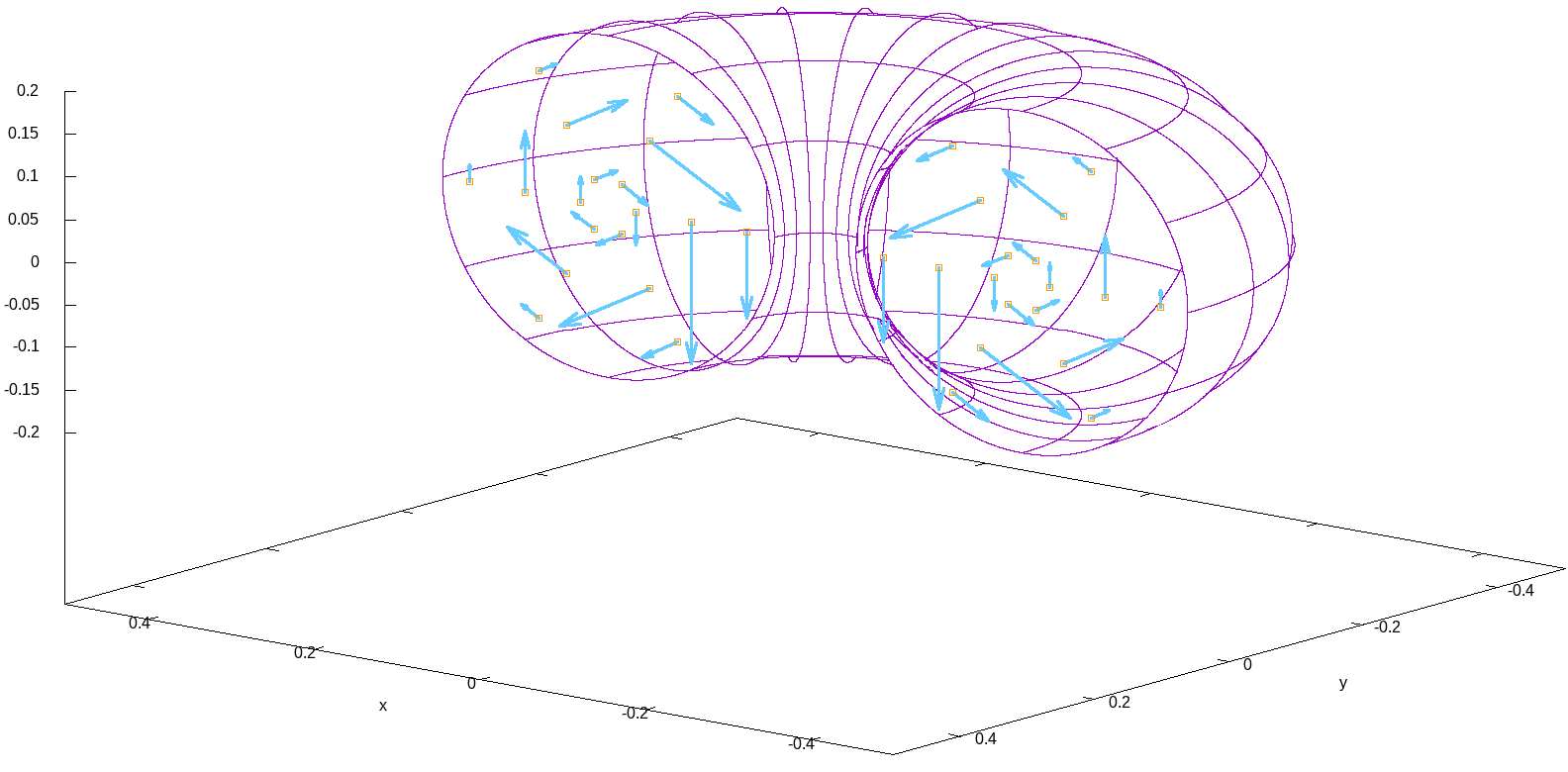}
 \caption{View of the vortex ring setup, showing the toroidal boundary of the support of the velocity, cut open along $z = 0$ (\emph{left}) and $y=0$ (\emph{right}). Arrows show the velocity in these planes.}
 \label{fig:torus-views}
\end{figure}

\begin{figure}
 \centering
 \includegraphics[width=0.49\textwidth]{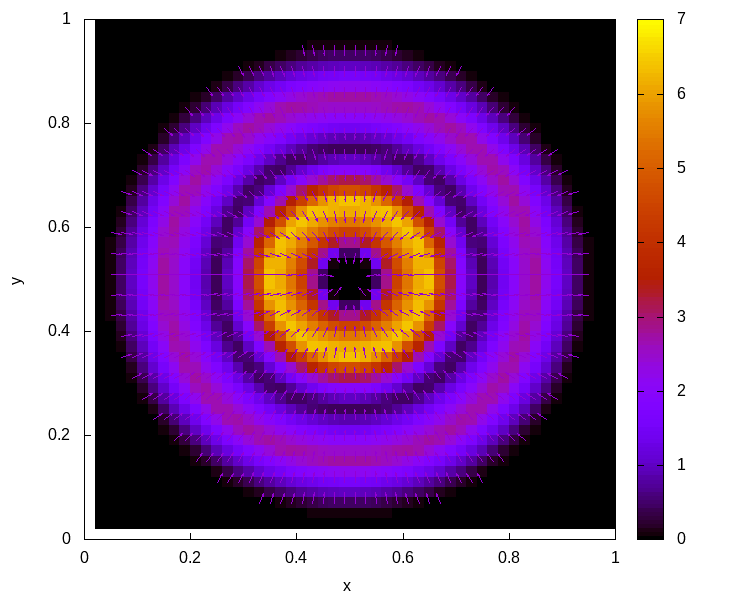} \hfill \includegraphics[width=0.49\textwidth]{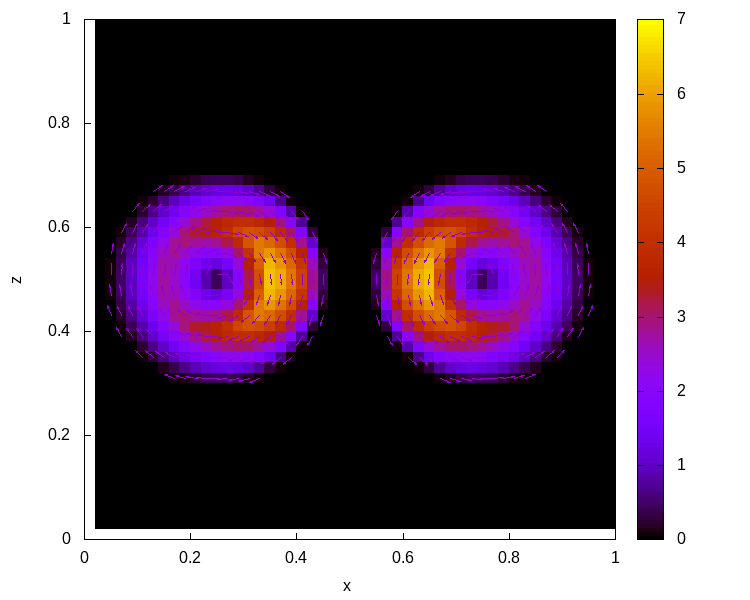} \\
 \includegraphics[width=0.49\textwidth]{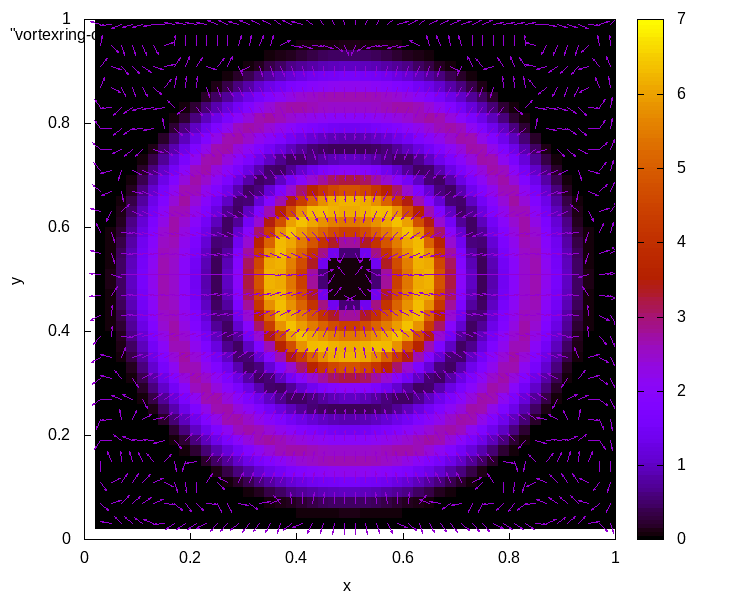} \hfill \includegraphics[width=0.49\textwidth]{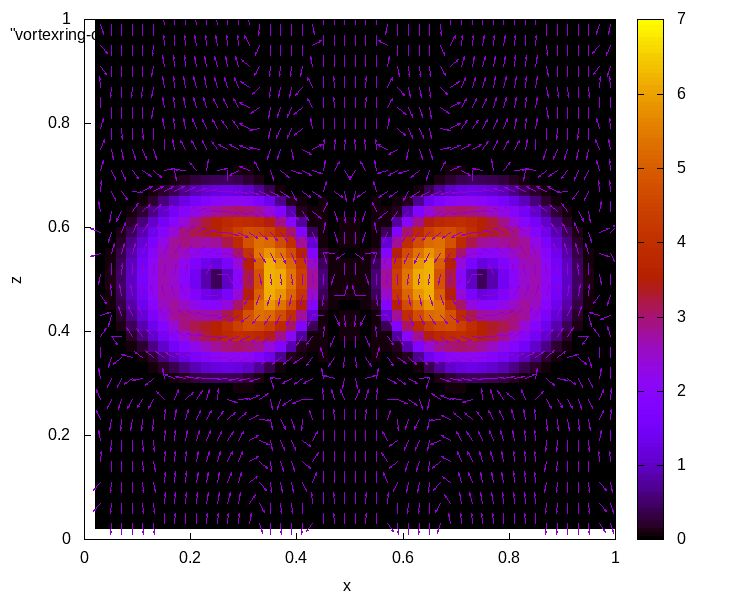}
 \caption{Vortex ring test case. \emph{Top}: Cut of the initial data through the plane $z = 0$ (\emph{left}) and $y = 0$ (\emph{right}). Color coded is the magnitude of the velocity, and the arrows (normalized in length) indicate the direction. \emph{Bottom}: The same for the numerical solution at time $t = 25$. The presence of non-zero velocity in regions that initially were at rest shows that the numerical stationary solution is slightly different from the initial data (which is natural as they were not well-prepared). The arrows indicate only the direction of the velocity, its magnitude is visually indistinguishable from zero in regions initially at rest.}
 \label{fig:vortexring}
\end{figure}

\begin{figure}
 \centering
 \includegraphics[width=0.49\textwidth]{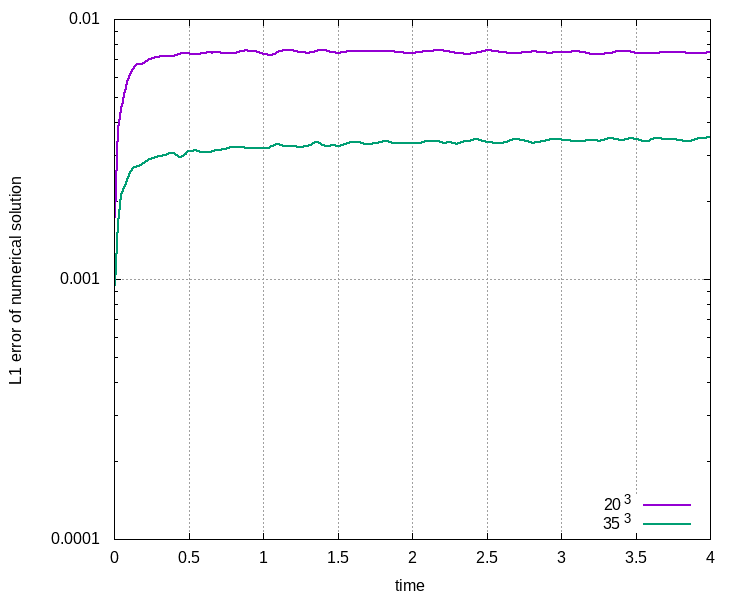} 
 \caption{Vortex ring test case. $L^1$ error of the numerical solution on grids of $20^3$ and $35^3$ cells, computed on the nodal point values. One observes the stationarization of the setup, as it converges towards the numerical stationary state.}
 \label{fig:vortexringerror}
\end{figure}

The vector field, divergence-free by construction, is then given by
\begin{align}
    u &= - \sin\theta\cos\phi \frac{V(r)}{\sqrt{x^2 + y^2}} \\
    v &= - \sin\theta\sin\phi \frac{V(r)}{\sqrt{x^2 + y^2}} \\
    w &= \cos\theta \frac{V(r)}{\sqrt{x^2 + y^2}}
\end{align}
Here $V$ is the radial velocity profile of the vortex ring and is an arbitrary differentiable, for simplicity compactly supported, function of $r$ only. We choose
\begin{align}
 V(r) = \begin{cases} 10r & r < 0.1 \\ \max(0, 2 - 10r) & \text{else} \end{cases}
\end{align}
and $R = \frac14$. The pressure $p$ vanishes identically. The setup is shown in Figure \ref{fig:vortexring} (top), in slices indicated in Figure \ref{fig:torus-views}. We solve it on a $50^3$ grid covering $[0,1]^3$ using RK3 and a CFL of 0.1. Results at a later time ($t=25$) are shown in Figure \ref{fig:vortexring} (bottom), they are virtually indistinguishable. Figure \ref{fig:vortexringerror} shows the numerical error as a function of time. One observes an initial layer (as the data are not discretely well-prepared) and the subsequent stationarization of the setup. A stationarity \emph{non}-preserving method will diffuse the vortex ring very quickly, converging for $t\to\infty$ to a discrete stationary state that is not a consistent discretization of a vortex in any way.

\subsection{Comparison to the fully discrete Active Flux method}
As stated in Corollary \ref{remark:statStatesFullyVsSemiDiscrete}, the numerical stationary states of the semi-discrete (generalized) and classical Active Flux are the same. The following numerical test shows that the results of both schemes are also very similar for non-stationary setups for the 2-d acoustic equations.

\begin{figure}
 \centering
 \includegraphics[width=0.98\textwidth]{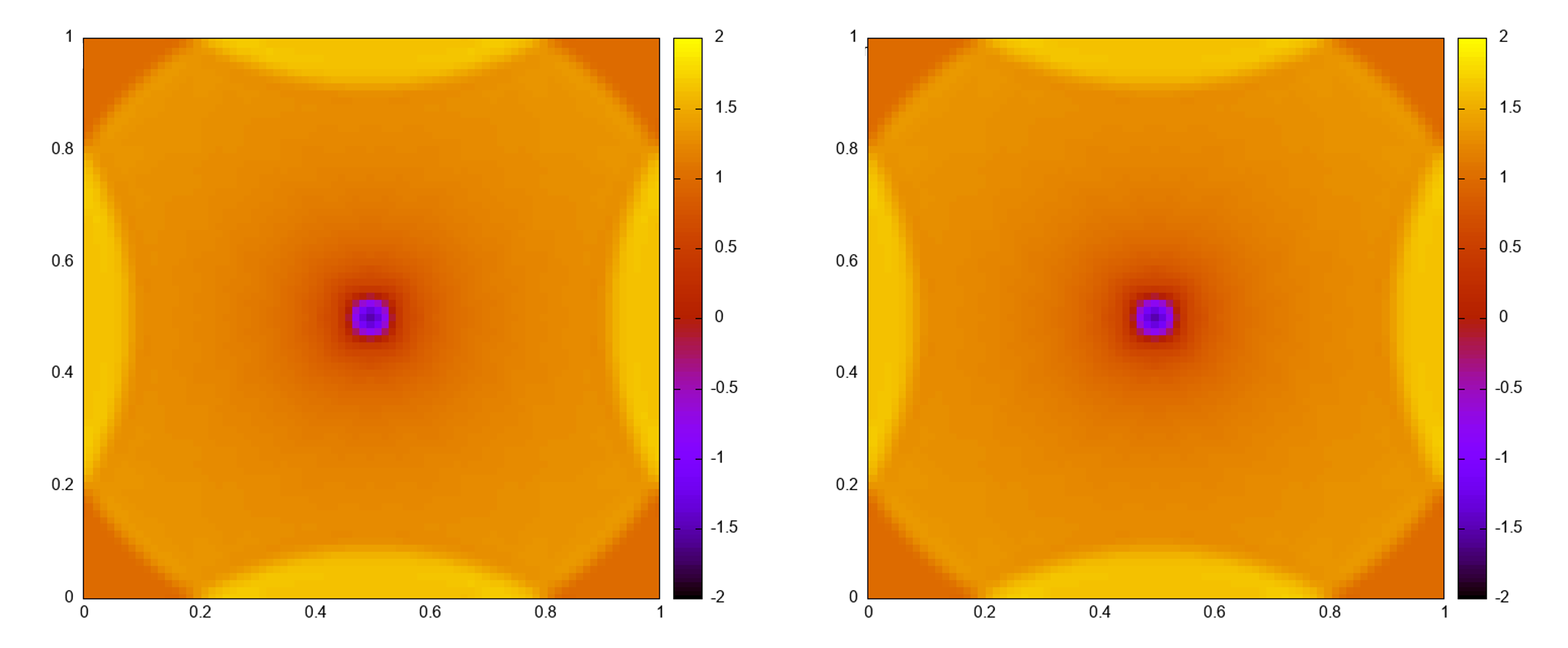} 
 \caption{Spherical Riemann problem. Numerical results for the pressure $p$ at time 0.3 on an 80$\times$80 grid and for CFL number 0.1. \emph{Left}: Classical active flux. \emph{Right}: Generalized active flux.}
 \label{fig:SphericalRiemannColour}
\end{figure}

\begin{figure}
 \centering
 \includegraphics[width=0.98\textwidth]{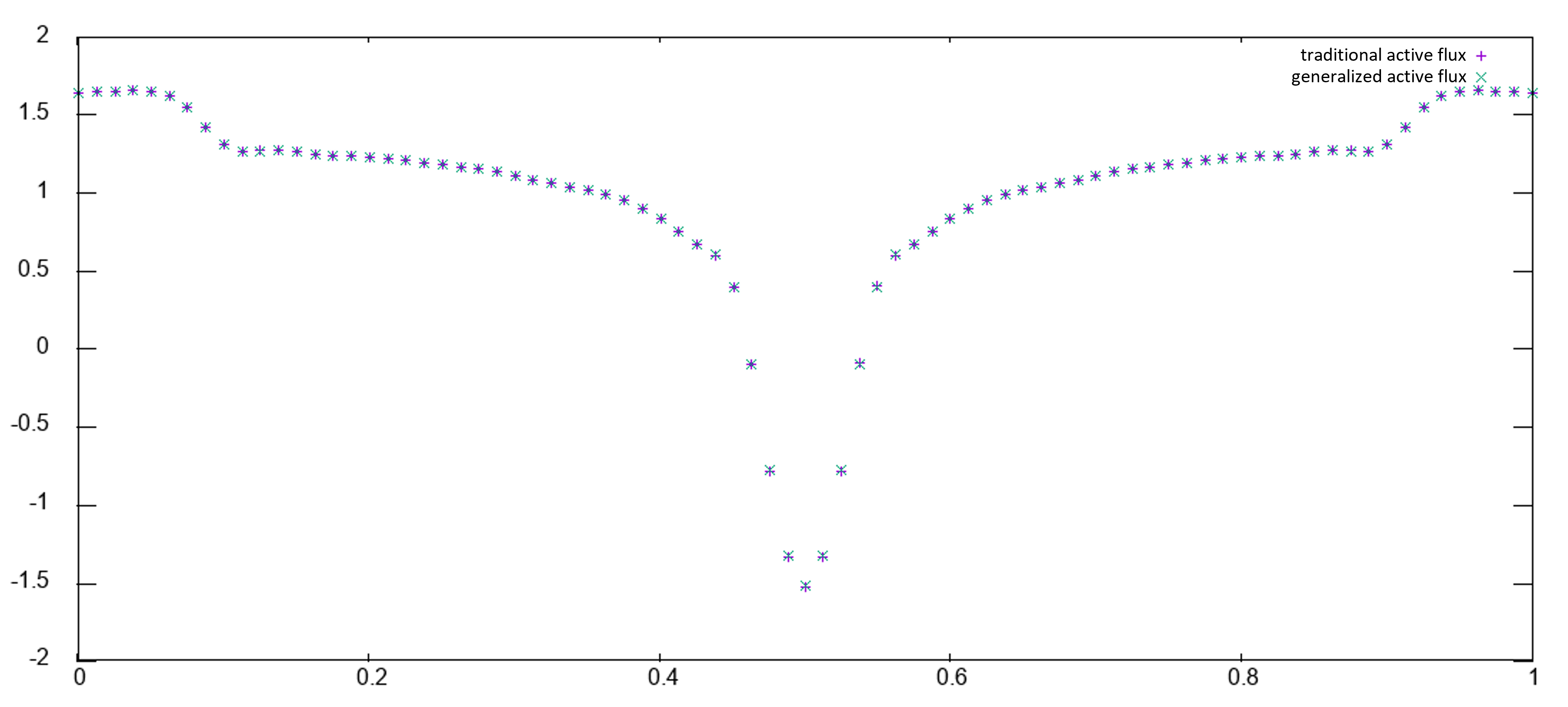} 
 \caption{Spherical Riemann problem. Numerical results for the pressure $p$ along the line $y=0.5$ at time 0.3 on an 80$\times$80 grid and for CFL number 0.1. }
 \label{fig:SphericalRiemannQuantitative}
\end{figure}

We consider the spherical Riemann problem
\begin{align}
\mathbf{v}_0(\mathbf{x})&=0\\
p_0(\mathbf{x})&=\begin{cases} 2 & \mathrm{if}\,\lvert\mathbf{x}\rvert^2<0.08\\ 1 & \mathrm{else}.\end{cases}
\end{align}
The simulation has been performed for different combinations of CFL numbers and grid sizes and with periodic boundary conditions. Figure \ref{fig:SphericalRiemannColour} shows a comparison of the numerical results of both methods for the pressure $p$ at time 0.3 on an 80$\times$80 grid with CFL number 0.1 as an example. No significant differences were found between the methods for any of these tested cases.
A more detailed comparison of the results of both methods along the line $y=0.5$ also shows no significant differences for $p$ at time 0.3. This is shown in Figure \ref{fig:SphericalRiemannQuantitative}, the grid size is again 80$\times$80 and the CFL number 0.1. 

\section{Conclusion}

We have shown that the two- and three-dimensional semi-discrete Active Flux method is stationarity preserving when applied to linear acoustics on Cartesian grids. This is the same conclusion as has been drawn for the classical Active Flux method before, the discrete stationary states are even exactly the same if upwind Jacobian splitting is used for the point value update. This is a consequence of the choice of the approximation/reconstruction space, common to both schemes. We also have not been able to find visible differences in the results of the two methods for non-stationary setups. However, the stability region of the semi-discrete method is smaller with a maximum CFL number only about half of that of the classical method (\cite{chudzik21}). The semi-discrete method can be applied more easily to nonlinear problems as it does not require an evolution operator, and it being structure preserving in very much the same way as the traditional method is encouraging. 

In the future we will study how this property can be extended to nonlinear conservation laws, including a theoretical analysis of low Mach number compliance for the Euler equations. Experimentally, it has been observed in \cite{barsukow24afeuler} that the method behaves well in this regime. Following \cite{morton01}, \cite{barsukow23nodal} we will investigate stationarity preservation of Active Flux on unstructured grids.

\section*{Acknowledgement}

WB, CK and LL acknowledge funding by the Deutsche Forschungsgemeinschaft (DFG, German Research Foundation) within \emph{SPP 2410 Hyperbolic Balance Laws in Fluid Mechanics: Complexity, Scales, Randomness (CoScaRa)}, project number 525941602.

\appendix

\section{Shape functions for Active Flux in two spatial dimensions} \label{app:basisfcts}

The basis function for the biparabolic reconstruction employed for the two-dimensional Active Flux are
\begin{align}
B_0(x, y)&=\frac94(-1+2\xi)(1+2\xi)(-1+2\eta)(1+2\eta).\label{basisfunktionenAF2Dneun}\\
B_{1}(x, y)&=\frac{1}{16}(-1+2\xi)(-1+2\eta)(12\xi\eta-2\eta-2\xi-1)\label{basisfunktionenAF2Deins}\\
B_{2}(x, y)&=\frac14(1-2\xi)(1+2\xi)(-1+2\eta)(1+6\eta)\label{basisfunktionenAF2Dzwei}\\
B_{3}(x, y)&=\frac{1}{16}(1+2\xi)(-1+2\eta)(12\xi\eta+2\eta-2\xi+1)\label{basisfunktionenAF2Ddrei}\\
B_{4}(x, y)&=\frac{1}{4}(1+2\xi)(-1+6\xi)(1-2\eta)(1+2\eta)\label{basisfunktionenAF2Dvier}\\
B_{5}(x, y)&=\frac{1}{16}(1+2\xi)(1+2\eta)(12\xi\eta+2\eta+2\xi-1)\label{basisfunktionenAF2Dfunf}\\
B_{6}(x, y)&=\frac14(1-2\xi)(1+2\xi)(1+2\eta)(-1+6\eta)\label{basisfunktionenAF2Dsechs}\\
B_{7}(x, y)&=\frac{1}{16}(-1+2\xi)(1+2\eta)(12\xi\eta-2\eta+2\xi+1)\label{basisfunktionenAF2Dsieben}\\
B_{8}(x, y)&=\frac{1}{4}(-1+2\xi)(1+6\xi)(1-2\eta)(1+2\eta)\label{basisfunktionenAF2Dacht}\\
\end{align}
with $\xi = x/\Delta x$, $\eta = y/\Delta y$. For the numbering, see Figure \ref{fig:distributionpointvalues}.

\section{Update equation for 2-d upon the Fourier transform}

\label{app:fourierequations}

The update of the nodal point value can be written as
\begin{align}
0
=\frac{\mathrm{d}}{\mathrm{d}t}\hat{q}\node&+\left[J_x^+\frac{1}{\Delta x}\left(3+\frac{1}{t_x}\right)+J_x^-\frac{1}{\Delta x}\left(-3-t_x\right) +J_y^+\frac{1}{\Delta y}\left(\frac{1}{t_y}+3\right) \right . \\ \nonumber & \left.+J_y^-\frac{1}{\Delta y}\left(-t_y-3\right)\right]\hat{q}\node+\left[\frac{-4}{\Delta x}(J_x^+-J_x^- t_x)\right]\hat{q}\eh \\ \nonumber &+\left[\frac{-4}{\Delta y}(J_y^+-J_y^- t_y)\right]\hat{q}\ev.
\end{align}

The update equation for the point value at a vertical edge is 
\begin{align}
0
=\frac{\mathrm{d}}{\mathrm{d}t}\hat{q}\ev&+\left[-\frac{9}{\Delta x}(J_x^+-J_x^-t_x)\right]\hat{q}\avg+\left[\frac{1}{\Delta x}\left(J_x^+\left(1+\frac{1}{t_y}\right)+J_x^-\left(-\frac{t_x}{t_y}-t_x\right)\right)\right]\hat{q}\eh\\\nonumber&+\left[\frac{2}{\Delta x}\left(J_x^+\left(2+\frac{1}{t_x}\right)+J_x^-\left(-2t_x-2\right)\right)\right]\hat{q}\ev\\\nonumber&+\bigg[J_x^+\frac{1}{4\Delta x}\left(\frac{1}{t_xt_y}+\frac{1}{t_y}+1+\frac{1}{t_x}\right)+J_x^-\frac{1}{4\Delta x}\left(-\frac{1}{t_y}-\frac{t_x}{t_y}-t_x-1\right)\\\nonumber&\qquad\qquad+J_y^+\frac{1}{\Delta y}\left(1-\frac{1}{t_y}\right)+J_y^-\frac{1}{\Delta y}\left(1-\frac{1}{t_y}\right)\bigg]\hat{q}\node
\end{align}
and for the update of the point value at a horizontal edge
\begin{align}
0
=\frac{\mathrm{d}}{\mathrm{d}t}\hat{q}\eh&+\left[-\frac{9}{\Delta y}(J_y^+-J_y^-t_y)\right]\hat{q}\avg+\left[\frac{2}{\Delta y}\left(J_y^+\left(\frac{1}{t_y}+2\right)+J_y^-\left(-2-t_y\right)\right)\right]\hat{q}\eh\\\nonumber&+\left[\frac{1}{\Delta y}\left(J_y^+\left(1+\frac{1}{t_x}\right)+J_y^-\left(-t_y-\frac{t_y}{t_x}\right)\right)\right]\hat{q}\ev\\\nonumber&+\bigg[J_x^+\frac{1}{\Delta x}\left(1-\frac{1}{t_x}\right)+J_x^-\frac{1}{\Delta x}\left(1-\frac{1}{t_x}\right)+J_y^+\frac{1}{4\Delta y}\left(\frac{1}{t_xt_y}+\frac{1}{t_y}+1+\frac{1}{t_x}\right)\\\nonumber&\qquad\qquad+J_y^-\frac{1}{4 \Delta y}\left(-\frac{1}{t_x}-1-t_y-\frac{t_y}{t_x}\right)\bigg]\hat{q}\node.
\end{align}
The update of the average is computed by following the instructions from Section \ref{sec:degreesOfFreedom}
\begin{align}
0
&=\frac{\mathrm{d}}{\mathrm{d}t}\hat{q}\avg+\left[\frac{2}{3\Delta y}J_y\left(1-\frac{1}{t_y}\right)\right]\hat{q}\eh+\left[\frac{2}{3\Delta x}J_x\left(1-\frac{1}{t_x}\right)\right]\hat{q}\ev\\\nonumber
&\qquad\quad\,\,\,+\left[\frac{1}{6\Delta x}J_x\left(1-\frac{1}{t_x}+\frac{1}{t_y}-\frac{1}{t_xt_y}\right)+\frac{1}{6\Delta y}J_y\left(1-\frac{1}{t_y}+\frac{1}{t_x}-\frac{1}{t_x t_y}\right)\right]\hat{q}\node.
\end{align}

\section{Evolution matrices of Active Flux in two spatial dimensions} \label{sec:2devomatricesacoustics}
       
The submatrices of the evolution matrix \eqref{eq:evolutionMatrixBlock} for the acoustic equations \eqref{acousticEquationsOhneMatrizen} in two spatial dimensions are
\begin{align}
\mathcal{E}_{\avgonly\avgonly}&=0_{3\times 3} & \mathcal{E}_{\nodeonly\avgonly}&=0_{3\times 3} 
\end{align}
\begin{equation}
\mathcal{E}_{\avgonly\ehonly}=\frac{2}{3\Delta y}\left(\begin{matrix} 0&0&0\\ 0&0&1-\frac{1}{t_y} \\ 0&1-\frac{1}{t_y} &0\end{matrix}\right)
\end{equation}
\begin{equation}
\mathcal{E}_{\avgonly\evonly}=\frac{2}{3\Delta x}\left(\begin{matrix} 0&0&1-\frac{1}{t_x}\\ 0&0&0 \\ 1-\frac{1}{t_x}&0 &0\end{matrix}\right)
\end{equation}
\begin{align}
\mathcal{E}_{\avgonly\nodeonly}&=\left(\begin{matrix} 0&0&\frac{1}{6 \Delta x}\frac{(-1+t_x)(1+t_y)}{t_xt_y}\\ 0&0& \frac{1}{6 \Delta y}\frac{(1+t_x)(-1+t_y)}{t_xt_y}\\\frac{1}{6 \Delta x}\frac{(-1+t_x)(1+t_y)}{t_xt_y}& \frac{1}{6 \Delta y}\frac{(1+t_x)(-1+t_y)}{t_xt_y}&0\end{matrix}\right)
\end{align}   
\begin{equation}
\mathcal{E}_{\ehonly\avgonly}=-\frac{9}{2\Delta y}\left(\begin{matrix} 0&0&0\\0&(1+t_y)&(1-t_y)\\0&(1-t_y)&(1+t_y)\end{matrix}\right)
\end{equation}
\begin{equation}
\mathcal{E}_{\ehonly\ehonly}=\frac{1}{\Delta y}\left(\begin{matrix} 0&0&0\\0&(\frac{1}{t_y}+4+t_y)&(\frac{1}{t_y}-t_y)\\0&(\frac{1}{t_y}-t_y)&(\frac{1}{t_y}+4+t_y)\end{matrix}\right)
\end{equation}
\begin{equation}
\mathcal{E}_{\ehonly\evonly}=\frac{1}{2\Delta y}\left(\begin{matrix} 0&0&0\\0&\frac{(1+t_x)(1+t_y)}{t_x}&-\frac{(1+t_x)(-1+t_y)}{t_x}\\0&-\frac{(1+t_x)(-1+t_y)}{t_x}&\frac{(1+t_x)(1+t_y)}{t_x})\end{matrix}\right)
\end{equation}
\begin{equation}
\mathcal{E}_{\ehonly\nodeonly}=\left(\begin{matrix} 0&0&\frac{1}{\Delta x}\left(1-\frac{1}{t_x}\right) \\0&\frac{1}{8\Delta y}\frac{(1+t_x)(1+t_y)^2}{t_xt_y}&-\frac{1}{8 \Delta y}\frac{(1+t_x)(-1+t_y^2)}{t_xt_y}\\\frac{1}{\Delta x}\left(1-\frac{1}{t_x}\right)&-\frac{1}{8\Delta y}\frac{(1+t_x)(1+t_y)^2}{t_xt_y}&\frac{1}{8 \Delta y}\frac{(1+t_x)(-1+t_y^2)}{t_xt_y}\end{matrix}\right)
\end{equation}
\begin{equation}
\mathcal{E}_{\evonly\avgonly}=-\frac{9}{2\Delta y}\left(\begin{matrix} (1+t_x)&0&(1-t_x)\\0&0&0\\(1-t_x)&0&(1+t_x)\end{matrix}\right)
\end{equation}
\begin{equation}
\mathcal{E}_{\evonly\ehonly}=\frac{1}{2\Delta x}\left(\begin{matrix} \frac{(1+t_x)(1+t_y)}{t_y}&0&-\frac{(-1+t_x)(1+t_y)}{t_y}\\0&0&0\\-\frac{(-1+t_x)(1+t_y)}{t_y}&0&\frac{(1+t_x)(1+t_y)}{t_y})\end{matrix}\right)
\end{equation}
\begin{equation}
\mathcal{E}_{\evonly\evonly}= \frac{1}{\Delta x}\left(\begin{matrix}(\frac{1}{t_x}+4+t_x)&0&(\frac{1}{t_x}-t_x)\\0&0&0\\(\frac{1}{t_x}-t_x)&0&(\frac{1}{t_x}+4+t_x)\end{matrix}\right)
\end{equation}
\begin{equation}
\mathcal{E}_{\evonly\nodeonly}=\left(\begin{matrix} \frac{1}{8\Delta x}\frac{(1+t_x)^2(1+t_y)}{t_xt_y}&0&-\frac{1}{8\Delta x}\frac{(-1+t_x^2)(1+ t_y)}{t_xt_y} \\0&0&\frac{1}{\Delta y}\left(1-\frac{1}{t_y}\right)\\ -\frac{1}{8\Delta x}\frac{(-1+t_x^2)(1+ t_y)}{t_xt_y}&\frac{1}{\Delta y}\left(1-\frac{1}{t_y}\right)&\frac{1}{8\Delta x}\frac{(1+t_x)^2(1+t_y)}{t_xt_y}\end{matrix}\right)
\end{equation}

\begin{equation}
\mathcal{E}_{\nodeonly\ehonly}= -\frac{2}{\Delta x}\left(\begin{matrix}(1+t_x)&0&(1-t_x)\\0&0&0\\(1-t_x)&0&(1+t_x)\end{matrix}\right),
\end{equation}
\begin{equation}
\mathcal{E}_{\nodeonly\evonly}=-\frac{2}{\Delta y}\left(\begin{matrix} 0&0&0\\0&(1+t_y)&(1-t_y)\\0&(1-t_y)&(1+t_y)\end{matrix}\right)
\end{equation}
\begin{equation}
\mathcal{E}_{\nodeonly\nodeonly}=\left(\begin{matrix} \frac{6+\frac{1}{t_x}+t_x}{2\Delta x} &0&\frac{\frac{1}{t_x}-t_x}{2\Delta x} \\ 0& \frac{1}{2\Delta y}\left(6+\frac{1}{t_y}+t_y\right) &\frac{1}{2\Delta y}\left(\frac{1}{t_y}-t_y\right) \\ \frac{1}{2\Delta x}\left(\frac{1}{t_x}-t_x\right)& \frac{1}{2\Delta y}\left(\frac{1}{t_y}-t_y\right)&\frac{1}{2}\left[\frac{6+\frac{1}{t_x}+t_x}{\Delta x}+\frac{6+\frac{1}{t_y}+t_y}{\Delta y}\right]\end{matrix}\right),
\end{equation}

\section{Evolution matrix and its kernel for Active Flux in three spatial dimensions} \label{sec:3dkernelacoustics}
The variables are ordered as $\hat{Q}=(\hat{q}\avg, \hat{q}\ex, \hat{q}\ey, \hat{q}\ez, \hat{q}\fx, \hat{q}\fy, \hat{q}\fz, \hat{q}\node)$, where $\hat{q}=(\hat{u}, \hat{v}, \hat{w}, \hat{p}).$



The evolution matrix reads:
\begin{align}
 \mathcal E &= \left( \begin{matrix} \mathcal E_{\avgonly\avgonly} & \mathcal E_{\avgonly\exonly} & \mathcal E_{\avgonly\eyonly} & \mathcal E_{\avgonly\ezonly} & \mathcal E_{\avgonly\fxonly} &\mathcal E_{\avgonly\fyonly} &\mathcal E_{\avgonly\fzonly} & \mathcal E_{\avgonly\nodeonly} \\
 \mathcal E_{\exonly\avgonly} & \mathcal E_{\exonly\exonly} & \mathcal E_{\exonly\eyonly} & \mathcal E_{\exonly\ezonly} & \mathcal E_{\exonly\fxonly} &\mathcal E_{\exonly\fyonly} &\mathcal E_{\exonly\fzonly} & \mathcal E_{\exonly\nodeonly} \\
 \mathcal E_{\eyonly\avgonly} & \mathcal E_{\eyonly\exonly} & \mathcal E_{\eyonly\eyonly} & \mathcal E_{\eyonly\ezonly} & \mathcal E_{\eyonly\fxonly} &\mathcal E_{\eyonly\fyonly} &\mathcal E_{\eyonly\fzonly} & \mathcal E_{\eyonly\nodeonly} \\
 \mathcal E_{\ezonly\avgonly} & \mathcal E_{\ezonly\exonly} & \mathcal E_{\ezonly\eyonly} & \mathcal E_{\ezonly\ezonly} & \mathcal E_{\ezonly\fxonly} &\mathcal E_{\ezonly\fyonly} &\mathcal E_{\ezonly\fzonly} & \mathcal E_{\ezonly\nodeonly} \\
 \mathcal E_{\fxonly\avgonly} & \mathcal E_{\fxonly\exonly} & \mathcal E_{\fxonly\eyonly} & \mathcal E_{\fxonly\ezonly} & \mathcal E_{\fxonly\fxonly} &\mathcal E_{\fxonly\fyonly} &\mathcal E_{\fxonly\fzonly} & \mathcal E_{\fxonly\nodeonly} \\
 \mathcal E_{\fyonly\avgonly} & \mathcal E_{\fyonly\exonly} & \mathcal E_{\fyonly\eyonly} & \mathcal E_{\fyonly\ezonly} & \mathcal E_{\fyonly\fxonly} &\mathcal E_{\fyonly\fyonly} &\mathcal E_{\fyonly\fzonly} & \mathcal E_{\fyonly\nodeonly} \\
 \mathcal E_{\fzonly\avgonly} & \mathcal E_{\fzonly\exonly} & \mathcal E_{\fzonly\eyonly} & \mathcal E_{\fzonly\ezonly} & \mathcal E_{\fzonly\fxonly} &\mathcal E_{\fzonly\fyonly} &\mathcal E_{\fzonly\fzonly} & \mathcal E_{\fzonly\nodeonly} \\
 \mathcal E_{\nodeonly\avgonly} & \mathcal E_{\nodeonly\exonly} & \mathcal E_{\nodeonly\eyonly} & \mathcal E_{\nodeonly\ezonly} & \mathcal E_{\nodeonly\fxonly} &\mathcal E_{\nodeonly\fyonly} &\mathcal E_{\nodeonly\fzonly} & \mathcal E_{\nodeonly\nodeonly} \\\end{matrix} \right )
\end{align}

with the blocks 
\begin{align}
 \mathcal E_{\avgonly\avgonly}   &= \left(\begin{matrix} 0 & 0 & 0 & 0\\ 0 & 0 & 0 & 0\\  0 & 0 & 0 & 0 \\ 0 & 0 & 0 & 0  \end{matrix} \right ) \\
 \mathcal E_{\avgonly\exonly} &= \left(\begin{matrix} 0 & 0 & 0 & 0\\ 0 & 0 & 0 & \frac{c (t_y-1) (1+t_z)}{9 \Delta y t_y t_z}\\
 0 & 0 & 0 & \frac{c (1+t_y) (t_z-1)}{9 \Delta z t_y t_z}\\ 0 & \frac{c (t_y-1) (1+t_z)}{9 \Delta y t_y t_z} & \frac{c (1+t_y) (t_z-1)}{9 \Delta z t_y t_z} & 0
 \end{matrix} \right )\\
 \mathcal E_{\avgonly\eyonly} &= \left(\begin{matrix} 0 & 0 & 0 & \frac{c (t_x-1) (1+t_z)}{9 \Delta x t_x t_z}\\0 & 0 & 0 & 0 \\
 0 & 0 & 0 & \frac{c (1+t_x) (t_z-1)}{9 \Delta z t_x t_z}\\ \frac{c (t_x-1) (1+t_z)}{9 \Delta x t_x t_z} & 0 & \frac{c (1+t_x) (t_z-1)}{9 \Delta z t_x t_z} & 0 \end{matrix} \right )\\
 \mathcal E_{\avgonly\ezonly} &= \left(\begin{matrix} 0 & 0 & 0 & \frac{c (t_x-1) (1+t_y)}{9 \Delta x t_x t_y} \\ 0 & 0 & 0 & \frac{c (1+t_x) (t_y-1)}{9 \Delta y t_x t_y}\\0 & 0 & 0 & 0\\\frac{c (t_x-1) (1+t_y)}{9 \Delta x t_x t_y} & \frac{c (1+t_x) (t_y-1)}{9 \Delta y t_x t_y} & 0 & 0 
 \end{matrix} \right )\\
 \mathcal E_{\avgonly\fxonly} &= \left(\begin{matrix} 0 & 0 & 0 & \frac{4 c (t_x-1)}{9 \Delta x t_x}\\ 0 & 0 & 0 & 0 \\0 & 0 & 0 & 0\\\frac{4 c (t_x-1)}{9 \Delta x t_x} & 0 & 0 & 0
 \end{matrix} \right )
 \end{align}

\begin{align}
 \mathcal E_{\avgonly\fyonly} &= \left(\begin{matrix} 0 & 0 & 0 & 0\\ 0 & 0 & 0 & \frac{4 c (t_y-1)}{9 \Delta y t_y}\\0 & 0 & 0 & 0\\0 & \frac{4 c (t_y-1)}{9 \Delta y t_y} & 0 & 0
 \end{matrix} \right )\\
 \mathcal E_{\avgonly\fzonly} &= \left(\begin{matrix} 0 & 0 & 0 & 0\\  0 & 0 & 0 & 0 \\0 & 0 & 0 & \frac{4 c (t_z-1)}{9 \Delta z t_z}\\0 & 0 & \frac{4 c (t_z-1)}{9 \Delta z t_z} & 0
 \end{matrix} \right )\\
 \mathcal E_{\avgonly\nodeonly}   &= \left(\begin{matrix} 0 & 0 & 0 & \frac{c (t_x-1) (1+t_y) (1+t_z)}{36 \Delta x t_x t_y t_z}\\  0 & 0 & 0 & \frac{c (1+t_x) (t_y-1) (1+t_z)}{36 \Delta y t_x t_y t_z}\\0 & 0 & 0 & \frac{c (1+t_x) (1+t_y) (t_z-1)}{36 \Delta z t_x t_y t_z}\\\frac{c (t_x-1) (1+t_y) (1+t_z)}{36 \Delta x t_x t_y t_z} & \frac{c (1+t_x) (t_y-1) (1+t_z)}{36 \Delta y t_x t_y t_z} & \frac{c (1+t_x) (1+t_y) (t_z-1)}{36 \Delta z t_x t_y t_z} & 0
 \end{matrix} \right )
\\
 \mathcal E_{\exonly\avgonly}   &= \left(\begin{matrix} 0 & 0 & 0 & 0\\ 0 & 0 & 0 & 0\\ 0 & 0 & 0 & 0 \\0 & 0 & 0 & 0 \end{matrix} \right )\\
 \mathcal E_{\exonly\exonly} &= \left(\begin{matrix} 0 & 0 & 0 & 0\\ 0 & \frac{c+c t_y (6+t_y)}{2 \Delta y t_y} & 0 & -\frac{c (t_y^2-1)}{2 \Delta y t_y}\\0 & 0 & \frac{c+c t_z (6+t_z)}{2 \Delta z t_z} & -\frac{c (t_z^2-1)}{2 \Delta z t_z}\\0 & -\frac{c (t_y^2-1)}{2 \Delta y t_y} & -\frac{c (t_z^2-1)}{2 \Delta z t_z} & \frac12 c \left(\frac{6+1/t_y+t_y}{\Delta y}+\frac{1+t_z (6+t_z)}{\Delta z t_z}\right)
 \end{matrix} \right )\\
 \mathcal E_{\exonly\eyonly} &= \left(\begin{matrix} 0 & 0 & 0 & 0 \\ 0 & 0 & 0 & 0\\0 & 0 & 0 & 0\\0 & 0 & 0 & 0
 \end{matrix} \right )\\
 \mathcal E_{\exonly\ezonly} &= \left(\begin{matrix} 0 & 0 & 0 & 0\\ 0 & 0 & 0 & 0 \\0 & 0 & 0 & 0 \\0 & 0 & 0 & 0
 \end{matrix} \right )\\
 \mathcal E_{\exonly\fxonly} &= \left(\begin{matrix} 0 & 0 & 0 & 0\\ 0 & 0 & 0 & 0 \\ 0 & 0 & 0 & 0\\ 0 & 0 & 0 & 0 \end{matrix} \right )
 \end{align}

\begin{align}
 \mathcal E_{\exonly\fyonly} &= \left(\begin{matrix} 0 & 0 & 0 & 0 \\ 0 & 0 & 0 & 0\\ 0 & 0 & -\frac{2 c (1+t_z)}{\Delta z} & \frac{2 c (t_z-1)}{\Delta z} \\ 0 & 0 & \frac{2 c (t_z-1)}{\Delta z} & -\frac{2 c (1+t_z)}{\Delta z}\end{matrix} \right )\\
 \mathcal E_{\exonly\fzonly} &= \left(\begin{matrix} 0 & 0 & 0 & 0 \\ 0 & -\frac{2 c (1+t_y)}{\Delta y} & 0 & \frac{2 c (t_y-1)}{\Delta y} \\ 0 & 0 & 0 & 0 \\ 0 & \frac{2 c (t_y-1)}{\Delta y} & 0 & -\frac{2 c (1+t_y)}{\Delta y}\end{matrix} \right )\\
 \mathcal E_{\exonly\nodeonly}   &= \left(\begin{matrix} 0 & 0 & 0 & \frac{c (t_x-1)}{\Delta x t_x} \\0 & 0 & 0 & 0 \\ 0 & 0 & 0 & 0 \\  \frac{c (t_x-1)}{\Delta x t_x} & 0 & 0 & 0  \end{matrix} \right )
\\
 \mathcal E_{\eyonly\avgonly}   &= \left(\begin{matrix} 0 & 0 & 0 & 0\\ 0 & 0 & 0 & 0\\ 0 & 0 & 0 & 0 \\ 0 & 0 & 0 & 0\end{matrix} \right )\\
 \mathcal E_{\eyonly\exonly} &= \left(\begin{matrix} 0 & 0 & 0 & 0\\ 0 & 0 & 0 & 0 \\ 0 & 0 & 0 & 0 \\0 & 0 & 0 & 0\end{matrix} \right )\\
 \mathcal E_{\eyonly\eyonly} &= \left(\begin{matrix} \frac{c+c t_x (6+t_x)}{2 \Delta x t_x} & 0 & 0 & -\frac{c (t_x^2-1)}{2 \Delta x t_x}\\ 0 & 0 & 0 & 0\\  0 & 0 & \frac{c+c t_z (6+t_z)}{2 \Delta z t_z} & -\frac{c (t_z^2-1)}{2 \Delta z t_z} \\ -\frac{c (t_x^2-1)}{2 \Delta x t_x} & 0 & -\frac{c (t_z^2-1)}{2 \Delta z t_z} & \frac12 c \left(\frac{6+1/t_x+t_x}{\Delta x}+\frac{1+t_z (6+t_z)}{\Delta z t_z}\right)\end{matrix} \right )\\
 \mathcal E_{\eyonly\ezonly} &= \left(\begin{matrix} 0 & 0 & 0 & 0\\ 0 & 0 & 0 & 0\\  0 & 0 & 0 & 0 \\ 0 & 0 & 0 & 0\end{matrix} \right )\\
 \mathcal E_{\eyonly\fxonly} &= \left(\begin{matrix} 0 & 0 & 0 & 0\\ 0 & 0 & 0 & 0\\  0 & 0 & -\frac{2 c (1+t_z)}{\Delta z} & \frac{2 c (t_z-1)}{\Delta z} \\ 0 & 0 & \frac{2 c (t_z-1)}{\Delta z} & -\frac{2 c (1+t_z)}{\Delta z}\end{matrix} \right )
 \end{align}

\begin{align}
 \mathcal E_{\eyonly\fyonly} &= \left(\begin{matrix} 0 & 0 & 0 & 0\\ 0 & 0 & 0 & 0 \\  0 & 0 & 0 & 0 \\ 0 & 0 & 0 & 0\end{matrix} \right )\\
 \mathcal E_{\eyonly\fzonly} &= \left(\begin{matrix} -\frac{2 c (1+t_x)}{\Delta x} & 0 & 0 & \frac{2 c (t_x-1)}{\Delta x}\\ 0 & 0 & 0 & 0\\  0 & 0 & 0 & 0 \\ \frac{2 c (t_x-1)}{\Delta x} & 0 & 0 & -\frac{2 c (1+t_x)}{\Delta x}\end{matrix} \right )\\
 \mathcal E_{\eyonly\nodeonly}   &= \left(\begin{matrix} 0 & 0 & 0 & 0\\ 0 & 0 & 0 & \frac{c (t_y-1)}{\Delta y t_y}\\ 0 & 0 & 0 & 0 \\ 0 & \frac{c (t_y-1)}{\Delta y t_y} & 0 & 0  \end{matrix} \right )
\\
 \mathcal E_{\ezonly\avgonly}   &= \left(\begin{matrix} 0 & 0 & 0 & 0 \\ 0 & 0 & 0 & 0 \\ 0 & 0 & 0 & 0 \\ 0 & 0 & 0 & 0 \end{matrix} \right )\\
 \mathcal E_{\ezonly\exonly} &= \left(\begin{matrix} 0 & 0 & 0 & 0 \\ 0 & 0 & 0 & 0 \\ 0 & 0 & 0 & 0 \\ 0 & 0 & 0 & 0  \end{matrix} \right )\\
 \mathcal E_{\ezonly\eyonly} &= \left(\begin{matrix} 0 & 0 & 0 & 0 \\ 0 & 0 & 0 & 0 \\ 0 & 0 & 0 & 0 \\ 0 & 0 & 0 & 0  \end{matrix} \right )\\
 \mathcal E_{\ezonly\ezonly} &= \left(\begin{matrix} \frac{c+c t_x (6+t_x)}{2 \Delta x t_x} & 0 & 0 & -\frac{c (t_x^2-1)}{2 \Delta x t_x} \\ 0 & \frac{c+c t_y (6+t_y)}{2 \Delta y t_y} & 0 & -\frac{c (t_y^2-1)}{2 \Delta y t_y} \\ 0 & 0 & 0 & 0 \\  -\frac{c (t_x^2-1)}{2 \Delta x t_x} & -\frac{c (t_y^2-1)}{2 \Delta y t_y} & 0 & \frac12 c \left(\frac{6+1/t_x+t_x}{\Delta x}+\frac{1+t_y (6+t_y)}{\Delta y t_y}\right) \end{matrix} \right )\\
 \mathcal E_{\ezonly\fxonly} &= \left(\begin{matrix}  0 & 0 & 0 & 0\\ 0 & -\frac{2 c (1+t_y)}{\Delta y} & 0 & \frac{2 c (t_y-1)}{\Delta y} \\ 0 & 0 & 0 & 0 \\  0 & \frac{2 c (t_y-1)}{\Delta y} & 0 & -\frac{2 c (1+t_y)}{\Delta y} \end{matrix} \right )
 \end{align}

\begin{align}
 \mathcal E_{\ezonly\fyonly} &= \left(\begin{matrix} -\frac{2 c (1+t_x)}{\Delta x} & 0 & 0 & \frac{2 c (t_x-1)}{\Delta x} \\ 0 & 0 & 0 & 0 \\ 0 & 0 & 0 & 0 \\  \frac{2 c (t_x-1)}{\Delta x} & 0 & 0 & -\frac{2 c (1+t_x)}{\Delta x} \end{matrix} \right )\\
 \mathcal E_{\ezonly\fzonly} &= \left(\begin{matrix} 0 & 0 & 0 & 0 \\ 0 & 0 & 0 & 0 \\ 0 & 0 & 0 & 0 \\  0 & 0 & 0 & 0 \end{matrix} \right )\\
 \mathcal E_{\ezonly\nodeonly}   &= \left(\begin{matrix} 0 & 0 & 0 & 0 \\ 0 & 0 & 0 & 0 \\ 0 & 0 & 0 & \frac{c (t_z-1)}{\Delta z t_z} \\ 0 & 0 & \frac{c (t_z-1)}{\Delta z t_z} & 0  \end{matrix} \right )
\\
 \mathcal E_{\fxonly\avgonly}   &= \left(\begin{matrix} -\frac{27 c (1+t_x)}{4 \Delta x} & 0 & 0 & \frac{27 c (t_x-1)}{4 \Delta x} \\0 & 0 & 0 & 0  \\ 0 & 0 & 0 & 0 \\ \frac{27 c (t_x-1)}{4 \Delta x} & 0 & 0 & -\frac{27 c (1+t_x)}{4 \Delta x}  \end{matrix} \right )\\
 \mathcal E_{\fxonly\exonly} &= \left(\begin{matrix} \frac{c (1+t_x) (1+t_y) (1+t_z)}{8 \Delta x t_y t_z} & 0 & 0 & -\frac{c (t_x-1) (1+t_y) (1+t_z)}{8 \Delta x t_y t_z} \\ 0 & 0 & 0 & 0 \\ 0 & 0 & 0 & 0 \\  -\frac{c (t_x-1) (1+t_y) (1+t_z)}{8 \Delta x t_y t_z} & 0 & 0 & \frac{c (1+t_x) (1+t_y) (1+t_z)}{8 \Delta x t_y t_z} \end{matrix} \right )\\
 \mathcal E_{\fxonly\eyonly} &= \left(\begin{matrix} \frac{c (1+t_x)^2 (1+t_z)}{8 \Delta x t_x t_z} & 0 & 0 & -\frac{c (t_x^2-1) (1+t_z)}{8 \Delta x t_x t_z} \\ 0 & 0 & 0 & 0 \\ 0 & 0 & 0 & \frac{c (t_z-1)}{\Delta z t_z} \\  -\frac{c (t_x^2-1) (1+t_z)}{8 \Delta x t_x t_z} & 0 & \frac{c (t_z-1)}{\Delta z t_z} & \frac{c (1+t_x)^2 (1+t_z)}{8 \Delta x t_x t_z}  \end{matrix} \right )\\
 \mathcal E_{\fxonly\ezonly} &= \left(\begin{matrix} \frac{c (1+t_x)^2 (1+t_y)}{8 \Delta x t_x t_y} & 0 & 0 & -\frac{c (t_x^2-1) (1+t_y)}{8 \Delta x t_x t_y} \\ 0 & 0 & 0 & \frac{c (t_y-1)}{\Delta y t_y} \\ 0 & 0 & 0 & 0 \\  -\frac{c (t_x^2-1) (1+t_y)}{8 \Delta x t_x t_y} & \frac{c (t_y-1)}{\Delta y t_y} & 0 & \frac{c (1+t_x)^2 (1+t_y)}{8 \Delta x t_x t_y} \end{matrix} \right )\\
 \mathcal E_{\fxonly\fxonly} &= \left(\begin{matrix} \frac{c+c t_x (4+t_x)}{\Delta x t_x} & 0 & 0 & \frac{c-c t_x^2}{\Delta x t_x}\\ 0 & 0 & 0 & 0 \\ 0 & 0 & 0 & 0 \\ \frac{c-c t_x^2}{\Delta x t_x} & 0 & 0 & \frac{c+c t_x (4+t_x)}{\Delta x t_x}  \end{matrix} \right )
 \end{align}

\begin{align}
 \mathcal E_{\fxonly\fyonly} &= \left(\begin{matrix} \frac{c (1+t_x) (1+t_y)}{2 \Delta x t_y} & 0 & 0 & -\frac{c (t_x-1) (1+t_y)}{2 \Delta x t_y} \\ 0 & 0 & 0 & 0 \\ 0 & 0 & 0 & 0 \\ -\frac{c (t_x-1) (1+t_y)}{2 \Delta x t_y} & 0 & 0 & \frac{c (1+t_x) (1+t_y)}{2 \Delta x t_y} \end{matrix} \right )\\
 \mathcal E_{\fxonly\fzonly} &= \left(\begin{matrix} \frac{c (1+t_x) (1+t_z)}{2 \Delta x t_z} & 0 & 0 & -\frac{c (t_x-1) (1+t_z)}{2 \Delta x t_z} \\ 0 & 0 & 0 & 0 \\ 0 & 0 & 0 & 0 \\ -\frac{c (t_x-1) (1+t_z)}{2 \Delta x t_z} & 0 & 0 & \frac{c (1+t_x) (1+t_z)}{2 \Delta x t_z} \end{matrix} \right )\\
 \mathcal E_{\fxonly\nodeonly}   &= \left(\begin{matrix} \frac{c (1+t_x)^2 (1+t_y) (1+t_z)}{32 \Delta x t_x t_y t_z} & 0 & 0 & -\frac{c (t_x^2-1) (1+t_y) (1+t_z)}{32 \Delta x t_x t_y t_z}  \\ 0 & 0 & 0 & 0 \\ 0 & 0 & 0 & 0 \\  -\frac{c (t_x^2-1) (1+t_y) (1+t_z)}{32 \Delta x t_x t_y t_z} & 0 & 0 & \frac{c (1+t_x)^2 (1+t_y) (1+t_z)}{32 \Delta x t_x t_y t_z} \end{matrix} \right )
\\
 \mathcal E_{\fyonly\avgonly}   &= \left(\begin{matrix} 0 & 0 & 0 & 0 \\ 0 & -\frac{27 c (1+t_y)}{4 \Delta y} & 0 & \frac{27 c (t_y-1)}{4 \Delta y} \\ 0 & 0 & 0 & 0 \\ 0 & \frac{27 c (t_y-1)}{4 \Delta y} & 0 & -\frac{27 c (1+t_y)}{4 \Delta y}  \end{matrix} \right )\\
 \mathcal E_{\fyonly\exonly} &= \left(\begin{matrix} 0 & 0 & 0 & 0 \\ 0 & \frac{c (1+t_y)^2 (1+t_z)}{8 \Delta y t_y t_z} & 0 & -\frac{c (t_y^2-1) (1+t_z)}{8 \Delta y t_y t_z} \\ 0 & 0 & 0 & \frac{c (t_z-1)}{\Delta z t_z} \\ 0 & -\frac{c (t_y^2-1) (1+t_z)}{8 \Delta y t_y t_z} & \frac{c (t_z-1)}{\Delta z t_z} & \frac{c (1+t_y)^2 (1+t_z)}{8 \Delta y t_y t_z}  \end{matrix} \right )\\
 \mathcal E_{\fyonly\eyonly} &= \left(\begin{matrix} 0 & 0 & 0 & 0 \\ 0 & \frac{c (1+t_x) (1+t_y) (1+t_z)}{8 \Delta y t_x t_z} & 0 & -\frac{c (1+t_x) (t_y-1) (1+t_z)}{8 \Delta y t_x t_z} \\ 0 & 0 & 0 & 0 \\  0 & -\frac{c (1+t_x) (t_y-1) (1+t_z)}{8 \Delta y t_x t_z} & 0 & \frac{c (1+t_x) (1+t_y) (1+t_z)}{8 \Delta y t_x t_z}  \end{matrix} \right )\\
 \mathcal E_{\fyonly\ezonly} &= \left(\begin{matrix} 0 & 0 & 0 & \frac{c (t_x-1)}{\Delta x t_x} \\ 0 & \frac{c (1+t_x) (1+t_y)^2}{8 \Delta y t_x t_y} & 0 & -\frac{c (1+t_x) (t_y^2-1)}{8 \Delta y t_x t_y} \\ 0 & 0 & 0 & 0 \\  \frac{c (t_x-1)}{\Delta x t_x} & -\frac{c (1+t_x) (t_y^2-1)}{8 \Delta y t_x t_y} & 0 & \frac{c (1+t_x) (1+t_y)^2}{8 \Delta y t_x t_y} \end{matrix} \right )\\
 \mathcal E_{\fyonly\fxonly} &= \left(\begin{matrix} 0 & 0 & 0 & 0\\ 0 & \frac{c (1+t_x) (1+t_y)}{2 \Delta y t_x} & 0 & -\frac{c (1+t_x) (t_y-1)}{2 \Delta y t_x} \\0 & 0 & 0 & 0  \\ 0 & -\frac{c (1+t_x) (t_y-1)}{2 \Delta y t_x} & 0 & \frac{c (1+t_x) (1+t_y)}{2 \Delta y t_x}  \end{matrix} \right )
 \end{align}

\begin{align}
 \mathcal E_{\fyonly\fyonly} &= \left(\begin{matrix} 0 & 0 & 0 & 0 \\ 0 & \frac{c+c t_y (4+t_y)}{\Delta y t_y} & 0 & \frac{c-c t_y^2}{\Delta y t_y} \\ 0 & 0 & 0 & 0 \\ 0 & \frac{c-c t_y^2}{\Delta y t_y} & 0 & \frac{c+c t_y (4+t_y)}{\Delta y t_y} \end{matrix} \right )\\
 \mathcal E_{\fyonly\fzonly} &= \left(\begin{matrix} 0 & 0 & 0 & 0 \\ 0 & \frac{c (1+t_y) (1+t_z)}{2 \Delta y t_z} & 0 & -\frac{c (t_y-1) (1+t_z)}{2 \Delta y t_z} \\ 0 & 0 & 0 & 0 \\ 0 & -\frac{c (t_y-1) (1+t_z)}{2 \Delta y t_z} & 0 & \frac{c (1+t_y) (1+t_z)}{2 \Delta y t_z} \end{matrix} \right )\\
 \mathcal E_{\fyonly\nodeonly}   &= \left(\begin{matrix} 0 & 0 & 0 & 0 \\ 0 & \frac{c (1+t_x) (1+t_y)^2 (1+t_z)}{32 \Delta y t_x t_y t_z} & 0 & -\frac{c (1+t_x) (t_y^2-1) (1+t_z)}{32 \Delta y t_x t_y t_z} \\ 0 & 0 & 0 & 0 \\ 0 & -\frac{c (1+t_x) (t_y^2-1) (1+t_z)}{32 \Delta y t_x t_y t_z} & 0 & \frac{c (1+t_x) (1+t_y)^2 (1+t_z)}{32 \Delta y t_x t_y t_z}  \end{matrix} \right )
\\
 \mathcal E_{\fzonly\avgonly}   &= \left(\begin{matrix} 0 & 0 & 0 & 0 \\ 0 & 0 & 0 & 0 \\ 0 & 0 & -\frac{27 c (1+t_z)}{4 \Delta z} & \frac{27 c (t_z-1)}{4 \Delta z} \\ 0 & 0 & \frac{27 c (t_z-1)}{4 \Delta z} & -\frac{27 c (1+t_z)}{4 \Delta z}  \end{matrix} \right )\\
 \mathcal E_{\fzonly\exonly} &= \left(\begin{matrix} 0 & 0 & 0 & 0 \\  0 & 0 & 0 & \frac{c (t_y-1)}{\Delta y t_y}\\ 0 & 0 & \frac{c (1+t_y) (1+t_z)^2}{8 \Delta z t_y t_z} & -\frac{c (1+t_y) (t_z^2-1)}{8 \Delta z t_y t_z} \\ 0 & \frac{c (t_y-1)}{\Delta y t_y} & -\frac{c (1+t_y) (t_z^2-1)}{8 \Delta z t_y t_z} & \frac{c (1+t_y) (1+t_z)^2}{8 \Delta z t_y t_z } \end{matrix} \right )\\
 \mathcal E_{\fzonly\eyonly} &= \left(\begin{matrix} 0 & 0 & 0 & \frac{c (t_x-1)}{\Delta x t_x} \\ 0 & 0 & 0 & 0 \\ 0 & 0 & \frac{c (1+t_x) (1+t_z)^2}{8 \Delta z t_x t_z} & -\frac{c (1+t_x) (t_z^2-1)}{8 \Delta z t_x t_z} \\  \frac{c (t_x-1)}{\Delta x t_x} & 0 & -\frac{c (1+t_x) (t_z^2-1)}{8 \Delta z t_x t_z} & \frac{c (1+t_x) (1+t_z)^2}{8 \Delta z t_x t_z}  \end{matrix} \right )\\
 \mathcal E_{\fzonly\ezonly} &= \left(\begin{matrix} 0 & 0 & 0 & 0 \\ 0 & 0 & 0 & 0 \\ 0 & 0 & \frac{c (1+t_x) (1+t_y) (1+t_z)}{8 \Delta z t_x t_y} & -\frac{c (1+t_x) (1+t_y) (t_z-1)}{8 \Delta z t_x t_y} \\  0 & 0 & -\frac{c (1+t_x) (1+t_y) (t_z-1)}{8 \Delta z t_x t_y} & \frac{c (1+t_x) (1+t_y) (1+t_z)}{8 \Delta z t_x t_y} \end{matrix} \right )\\
 \mathcal E_{\fzonly\fxonly} &= \left(\begin{matrix} 0 & 0 & 0 & 0\\ 0 & 0 & 0 & 0 \\ 0 & 0 & \frac{c (1+t_x) (1+t_z)}{2 \Delta z t_x} & -\frac{c (1+t_x) (t_z-1)}{2 \Delta z t_x} \\ 0 & 0 & -\frac{c (1+t_x) (t_z-1)}{2 \Delta z t_x} & \frac{c (1+t_x) (1+t_z)}{2 \Delta z t_x}  \end{matrix} \right )
 \end{align}

\begin{align}
 \mathcal E_{\fzonly\fyonly} &= \left(\begin{matrix} 0 & 0 & 0 & 0 \\ 0 & 0 & 0 & 0 \\ 0 & 0 & \frac{c (1+t_y) (1+t_z)}{2 \Delta z t_y} & -\frac{c (1+t_y) (t_z-1)}{2 \Delta z t_y} \\ 0 & 0 & -\frac{c (1+t_y) (t_z-1)}{2 \Delta z t_y} & \frac{c (1+t_y) (1+t_z)}{2 \Delta z t_y} \end{matrix} \right )\\
 \mathcal E_{\fzonly\fzonly} &= \left(\begin{matrix} 0 & 0 & 0 & 0 \\ 0 & 0 & 0 & 0 \\ 0 & 0 & \frac{c+c t_z (4+t_z)}{\Delta z t_z} & \frac{c-c t_z^2}{\Delta z t_z} \\ 0 & 0 & \frac{c-c t_z^2}{\Delta z t_z} & \frac{c+c t_z (4+t_z)}{\Delta z t_z} \end{matrix} \right )\\
 \mathcal E_{\fzonly\nodeonly}   &= \left(\begin{matrix} 0 & 0 & 0 & 0 \\ 0 & 0 & 0 & 0 \\ 0 & 0 & \frac{c (1+t_x) (1+t_y) (1+t_z)^2}{32 \Delta z t_x t_y t_z} & -\frac{c (1+t_x) (1+t_y) (t_z^2-1)}{32 \Delta z t_x t_y t_z} \\ 0 & 0 & -\frac{c (1+t_x) (1+t_y) (t_z^2-1)}{32 \Delta z t_x t_y t_z} & \frac{c (1+t_x) (1+t_y) (1+t_z)^2}{32 \Delta z t_x t_y t_z}  \end{matrix} \right )
\\
 \mathcal E_{\nodeonly\avgonly}   &= \left(\begin{matrix} 0 & 0 & 0 & 0 \\ 0 & 0 & 0 & 0 \\ 0 & 0 & 0 & 0 \\ 0 & 0 & 0 & 0  \end{matrix} \right )\\
 \mathcal E_{\nodeonly\exonly} &= \left(\begin{matrix} -\frac{2 c (1+t_x)}{\Delta x} & 0 & 0 & \frac{2 c (t_x-1)}{\Delta x} \\ 0 & 0 & 0 & 0 \\ 0 & 0 & 0 & 0 \\ \frac{2 c (t_x-1)}{\Delta x} & 0 & 0 & -\frac{2 c (1+t_x)}{\Delta x}  \end{matrix} \right )\\
 \mathcal E_{\nodeonly\eyonly} &= \left(\begin{matrix} 0 & 0 & 0 & 0 \\0 & -\frac{2 c (1+t_y)}{\Delta y} & 0 & \frac{2 c (t_y-1)}{\Delta y}  \\ 0 & 0 & 0 & 0 \\ 0 & \frac{2 c (t_y-1)}{\Delta y} & 0 & -\frac{2 c (1+t_y)}{\Delta y}   \end{matrix} \right )\\
 \mathcal E_{\nodeonly\ezonly} &= \left(\begin{matrix} 0 & 0 & 0 & 0 \\ 0 & 0 & 0 & 0 \\ 0 & 0 & -\frac{2 c (1+t_z)}{\Delta z} & \frac{2 c (t_z-1)}{\Delta z} \\ 0 & 0 & \frac{2 c (t_z-1)}{\Delta z} & -\frac{2 c (1+t_z)}{\Delta z}  \end{matrix} \right )\\
 \mathcal E_{\nodeonly\fxonly} &= \left(\begin{matrix} 0 & 0 & 0 & 0\\  0 & 0 & 0 & 0\\ 0 & 0 & 0 & 0 \\ 0 & 0 & 0 & 0  \end{matrix} \right )\end{align}

\begin{align}
 \mathcal E_{\nodeonly\fyonly} &= \left(\begin{matrix} 0 & 0 & 0 & 0 \\ 0 & 0 & 0 & 0 \\ 0 & 0 & 0 & 0 \\ 0 & 0 & 0 & 0 \end{matrix} \right )\\
 \mathcal E_{\nodeonly\fzonly} &= \left(\begin{matrix} 0 & 0 & 0 & 0 \\ 0 & 0 & 0 & 0 \\ 0 & 0 & 0 & 0 \\ 0 & 0 & 0 & 0 \end{matrix} \right )\\
 \mathcal E_{\nodeonly\nodeonly}   &= \left(\begin{matrix} \frac{c+c t_x (6+t_x)}{2 \Delta x t_x} & 0 & 0 & -\frac{c (t_x^2-1)}{2 \Delta x t_x} \\ 0 & \frac{c+c t_y (6+t_y)}{2 \Delta y t_y} & 0 & -\frac{c (t_y^2-1)}{2 \Delta y t_y} \\ 0 & 0 & \frac{c+c t_z (6+t_z)}{2 \Delta z t_z} & -\frac{c (t_z^2-1)}{2 \Delta z t_z} \\ -\frac{c (t_x^2-1)}{2 \Delta x t_x} & -\frac{c (t_y^2-1)}{2 \Delta y t_y} & -\frac{c (t_z^2-1)}{2 \Delta z t_z} & \frac12 c \left( \frac{6+1/t_x+t_x}{\Delta x}+\frac{6+1/t_y+t_y}{\Delta y}+\frac{1+t_z (6+t_z)}{\Delta z t_z}\right)  \end{matrix} \right )
\end{align}

The kernel of $\mathcal E$ are spanned by $(\hat Q_1, \hat Q_2, \hat Q_3, \hat Q_4, \hat Q_5)$ given by

\begin{align}
\hat Q_1 = \Big ( &-\frac{\Delta x (1+t_x (4+t_x)) (t_y-1) (t_z^2+4 t_x t_z^2+t_x^2 (2+t_z (8+3 t_z)))}{54 \Delta y t_x^2 (t_x^2-1) t_y t_z^2}, \nonumber\\\nonumber & \frac{(1+t_y (4+t_y)) (t_z^2+4 t_x t_z^2+t_x^2 (2+t_z (8+3 t_z)))}{54 t_x^2 t_y (1+t_y) t_z^2},0,0,
\\\nonumber  &-\frac{\Delta x (1+t_x (6+t_x)) (t_y-1)}{4 \Delta y (t_x-1) t_x (1+t_y)},0, 
 \frac{\Delta z (1+t_x) (t_y-1) (1+t_z)}{2 \Delta y t_x (1+t_y) (t_z-1)},0,
\\\nonumber  & 0, 
\frac14+\frac{1}{4 t_y}+\frac{1}{1+t_y}, -\frac{\Delta z (t_y-1) (1+t_z)}{2 \Delta y t_y (t_z-1)},0, 
\\\nonumber  &-\frac{\Delta x (t_y-1) (4 t_x^2+t_x (3+19 t_x) t_z+(1+t_x) (2+9 t_x) t_z^2)}{12 \Delta y (t_x-1) t_x (1+t_y) t_z^2}, 
\\\nonumber &\frac{4 t_x^2+t_x (3+19 t_x) t_z+(1+t_x) (2+9 t_x) t_z^2}{12 t_x (1+t_x) t_z^2},0,0, 
\\\nonumber  &-\frac{\Delta x (t_y-1) (2 t_z^2+t_x^2 (4+t_z) (1+3 t_z)+t_x t_z (-3+5 t_z))}{24 \Delta y (t_x-1) t_x t_y t_z^2}, \\\nonumber  & \frac{(1+t_y (6+t_y)) (4 t_x^2+t_x (3+19 t_x) t_z+(1+t_x) (2+9 t_x) t_z^2)}{48 t_x (1+t_x) t_y (1+t_y) t_z^2}, \\\nonumber &-(\frac{\Delta z (t_y-1) (1+t_z (6+t_z))}{8 \Delta y t_y (t_z-1) t_z}),0, 
\\\nonumber  & -\frac{\Delta x (1+t_x (6+t_x)) (t_y-1) (4 t_x^2+t_x (3+19 t_x) t_z+(1+t_x) (2+9 t_x) t_z^2)}{48 \Delta y t_x^2 (t_x^2-1) (1+t_y) t_z^2},
\\\nonumber & \frac{2 t_z^2+t_x^2 (4+t_z) (1+3 t_z)+t_x t_z (-3+5 t_z)}{24 t_x^2 t_z^2},\frac{\Delta z (1+t_x) (t_y-1) (1+t_z (6+t_z))}{8 \Delta y t_x (1+t_y) (t_z-1) t_z},0, \\\nonumber 
& 0,0,0,0, \\  & -\frac{\Delta x (1+t_x) (t_y-1)}{\Delta y (t_x-1) (1+t_y)},1,0,0 \Big)
\end{align}
\begin{align}
\hat Q_2 = \Big ( &-\frac{\Delta x (1+t_x (4+t_x)) (1+t_y) (t_z-1) (2 t_z^2+t_x^2 (t_z-2) (1+3 t_z)+t_x t_z (3+11 t_z))}{108 \Delta z t_x^2 (t_x^2-1) t_y t_z^2 (1+t_z)}, \nonumber\\\nonumber  &
-\frac{\Delta y (1+t_y (4+t_y)) (t_x-t_z) (t_z-1) (t_x+t_z+4 t_x t_z)}{54 \Delta z t_x^2 (t_y-1) t_y t_z^2 (1+t_z)}, \\\nonumber & \frac{(1+t_x) (1+t_y) (1+t_z (4+t_z))}{36 t_x t_y t_z (1+t_z)},0,
\\\nonumber  &-\frac{(\Delta x (1+t_x (6+t_x)) (t_z-1)}{4 \Delta z (t_x-1) t_x (1+t_z)},0,\frac{1+t_x}{2 t_x},0, \\\nonumber 
& 0,0,0,0, \\\nonumber 
& -\frac{\Delta x (t_z-1) (t_z^2+t_x t_z (3+7 t_z)+t_x^2 (-1+t_z (-1+3 t_z)))}{6 \Delta z (t_x-1) t_x t_z^2 (1+t_z)}, \\\nonumber & -\frac{\Delta y (1+t_y) (t_x-t_z) (t_z-1) (t_x+t_z+4 t_x t_z)}{6 \Delta z t_x (1+t_x) (t_y-1) t_z^2 (1+t_z)},\frac14+\frac{1}{4 t_z}+\frac{1}{1+t_z},0,
\\\nonumber  & \frac{\Delta x (1+t_y) (t_x-t_z) (t_z-1) (t_x+t_z+4 t_x t_z)}{12 \Delta z (t_x-1) t_x t_y t_z^2 (1+t_z)},\\\nonumber  & -\frac{\Delta y (1+t_y (6+t_y)) (t_x-t_z) (t_z-1) (t_x+t_z+4 t_x t_z)}{24 \Delta z t_x (1+t_x) (t_y-1) t_y t_z^2 (1+t_z)},0,0, \\\nonumber 
& -\frac{\Delta x (1+t_x (6+t_x)) (t_z-1) (t_z^2+t_x t_z (3+7 t_z)+t_x^2 (-1+t_z (-1+3 t_z)))}{24 \Delta z t_x^2 (t_x^2-1) t_z^2 (1+t_z)}, \\\nonumber  & -\frac{\Delta y (1+t_y) (t_x-t_z) (t_z-1) (t_x+t_z+4 t_x t_z)}{12 \Delta z t_x^2 (t_y-1) t_z^2 (1+t_z)},\frac{(1+t_x) (1+t_z (6+t_z))}{8 t_x t_z (1+t_z)},0, \\\nonumber  & 0,0,0,0, \\ 
&-\frac{\Delta x (1+t_x) (t_z-1)}{\Delta z (t_x-1) (1+t_z)},0,1,0 \Big )
\end{align}
\begin{align}
\hat Q_3 = \Big ( & -\frac{2 \Delta x (1+t_x (4+t_x)) (t_z-1) }{9 \Delta z (t_x^2-1) t_z},0,\frac{2 (1+t_z (4+t_z)) }{9 t_z (1+t_z)},0,\nonumber \\\nonumber  &
0,-\frac{2 \Delta y t_y (t_z-1) }{\Delta z (t_y-1) (1+t_z))},\frac{2 t_y }{1+t_y},0,\\\nonumber  
& 0,0,0,0, \\\nonumber  & -\frac{2 \Delta x t_x t_y (t_z-1) }{\Delta z (t_x-1) (1+t_y) t_z},\frac{2 \Delta y t_x t_y (t_z-1) }{\Delta z (1+t_x) (t_y-1) t_z},0,0, \\\nonumber 
& -\frac{\Delta x t_x (t_z-1) }{\Delta z (t_x-1) t_z},\frac{\Delta y t_x (1+t_y (6+t_y)) (t_z-1) }{2 \Delta z (1+t_x) (t_y^2-1) t_z},0,0,\\\nonumber 
&-\frac{\Delta x (1+t_x (6+t_x)) t_y (t_z-1) }{2 \Delta z (t_x^2-1) (1+t_y) t_z},0,\frac{t_y (1+t_z (6+t_z)) }{2 (1+t_y) t_z (1+t_z)},0, \\ 
&0,-\frac{\Delta y (1+t_y (6+t_y)) (t_z-1) }{2 \Delta z (t_y^2-1) (1+t_z)},1,0, 0,0,0,0 \Big ) 
\end{align}
\begin{align}
\hat Q_4 = \Big ( & \frac{8 \Delta x (1+t_x (4+t_x)) (t_y^2-1) (t_x-t_z) (t_x+t_z+4 t_x t_z) }{27 \Delta y (t_x-1) t_x (1+t_x)^2 (1+t_y (6+t_y)) t_z^2},\nonumber \\\nonumber  &  -\frac{8 (1+t_y (4+t_y)) (t_x-t_z) (t_x+t_z+4 t_x t_z) }{27 t_x (1+t_x) (1+t_y (6+t_y)) t_z^2},0,0,\\\nonumber 
&0,\frac{4 t_y (1+t_y) }{1+t_y (6+t_y)},-\frac{4 \Delta z (t_y-1) t_y (1+t_z)}{\Delta y (1+t_y (6+t_y)) (t_z-1)},0, \\\nonumber 
&-\frac{4 \Delta x t_x (t_y^2-1) }{\Delta y (t_x-1) (1+t_y (6+t_y))},0,\frac{4 \Delta z t_x (t_y^2-1) (1+t_z) }{\Delta y (1+t_x) (1+t_y (6+t_y)) (t_z-1)},0, \\\nonumber 
&\frac{4 \Delta x (t_y-1) t_y (2 t_x^2+t_x (3+11 t_x) t_z+(-2+t_x) (1+3 t_x) t_z^2) }{3 \Delta y (t_x^2-1) (1+t_y (6+t_y)) t_z^2}, \\\nonumber  & -\frac{4 t_y (1+t_y) (2 t_x^2+t_x (3+11 t_x) t_z+(-2+t_x) (1+3 t_x) t_z^2) }{3 (1+t_x)^2 (1+t_y (6+t_y)) t_z^2},0,0, \\\nonumber 
&\frac{4 \Delta x (t_y^2-1) (t_x-t_z) (t_x+t_z+4 t_x t_z) }{3 \Delta y (t_x^2-1) (1+t_y (6+t_y)) t_z^2}, \\\nonumber &-\frac{(2 t_x^2+t_x (3+11 t_x) t_z+(-2+t_x) (1+3 t_x) t_z^2) }{3 (1+t_x)^2 t_z^2}, \\\nonumber  & \frac{\Delta z t_x (t_y^2-1) (1+t_z (6+t_z)) }{\Delta y (1+t_x) (1+t_y (6+t_y)) (t_z-1) t_z},0, \\\nonumber  
&\frac{\Delta x (1+t_x (6+t_x)) (t_y-1) t_y (2 t_x^2+t_x (3+11 t_x) t_z+(-2+t_x) (1+3 t_x) t_z^2) }{3 \Delta y (t_x-1) t_x (1+t_x)^2 (1+t_y (6+t_y)) t_z^2}, \\\nonumber  &   -\frac{4 t_y (1+t_y) (t_x-t_z) (t_x+t_z+4 t_x t_z) }{3 t_x (1+t_x) (1+t_y (6+t_y)) t_z^2},-\frac{\Delta z (t_y-1) t_y (1+t_z (6+t_z)) }{\Delta y (1+t_y (6+t_y)) (t_z-1) t_z)},0, \\  &
-\frac{\Delta x (1+t_x (6+t_x)) (t_y^2-1) }{\Delta y (t_x^2-1) (1+t_y (6+t_y))},1,0,0,0,0,0,0 \Big ) 
\end{align}
\begin{align}
\hat Q_5 = \Big ( & 0,0,0,0,\quad 0,0,0,0, \nonumber \\\nonumber  
&\frac{4 t_x (1+t_x)}{1+t_x (6+t_x)},-\frac{4 \Delta y t_x (t_x^2-1) (1+t_y (6+t_y))}{\Delta x (1+t_x) (1+t_x (6+t_x)) (t_y^2-1)},\\\nonumber  &\frac{4 \Delta z t_x (t_x^2-1) (1+t_z)}{\Delta x (1+t_x) (1+t_x (6+t_x)) (t_z-1)},0, \\\nonumber 
&-\frac{4 t_x (1+t_x) (t_y-1) t_y (1+t_z)}{(1+t_x (6+t_x)) (t_y^2-1) t_z},-\frac{4 \Delta y t_x (t_x^2-1) t_y (1+t_y) (1+t_z)}{\Delta x (1+t_x) (1+t_x (6+t_x)) (t_y^2-1) t_z}, \\\nonumber  &
\frac{4 \Delta z t_x (t_x^2-1) (t_y-1) t_y (1+t_z (6+t_z))}{\Delta x (1+t_x) (1+t_x (6+t_x)) (t_y^2-1) (t_z-1) t_z},0, \\\nonumber 
& 0,-\frac{\Delta y t_x (t_x^2-1) (1+t_y (6+t_y)) (1+t_z)}{\Delta x (1+t_x) (1+t_x (6+t_x)) (t_y^2-1) t_z},\frac{\Delta z t_x (t_x^2-1) (1+t_z (6+t_z))}{\Delta x (1+t_x) (1+t_x (6+t_x)) (t_z-1) t_z},0,\\\nonumber 
&-\frac{(t_y-1) t_y (1+t_z)}{(t_y^2-1) t_z},\frac{2 \Delta y (t_x-1) (1+t_x) t_y (1+t_y) (1+t_z)}{\Delta x (1+t_x (6+t_x)) (t_y^2-1) t_z},0,0,\\\nonumber 
&1,0,-\frac{2 \Delta z (t_x-1) (1+t_x) (1+t_z)}{\Delta x (1+t_x (6+t_x)) (t_z-1)},0,\\\nonumber
&0,-\frac{16 \Delta y t_x (t_x^2-1) t_y (1+t_y)}{\Delta x (1+t_x) (1+t_x (6+t_x)) (t_y^2-1)}, \\ &\frac{16 \Delta z t_x (t_x^2-1) (t_y-1) t_y (1+t_z)}{\Delta x (1+t_x) (1+t_x (6+t_x)) (t_y^2-1) (t_z-1)},0 \Big )
\end{align}

To implement one discrete Fourier mode in Section \ref{sec:stationarymode3d} we use 
\begin{align}
 \sum_{r = 1}^5 a_r \hat Q_r
\end{align}
with {\footnotesize
\begin{align}
 a_1 &= \Big ( (t_x-1) (2 \Delta y (1+t_y+8 t_x t_y) (t_x-t_z) (t_z-1) (t_x+t_z+4 t_x t_z) \\\nonumber 
  &+t_x (t_y^2-1) (-32 t_x^2 t_y-3 (1+t_x+32 t_x^2 t_y) t_z+8 (4 t_y+t_x (1+t_x) (3 \Delta z+16 t_y)) t_z^2 \\\nonumber 
  & +(3 (1+t_x)-32 (1+4 t_x) t_y) t_z^3)) \Big ) / \Big (2 \Delta x \Delta z (t_y-1) (t_z^2+4 t_x t_z^2+t_x^2 (2+t_z (8+3 t_z))) \Big )\\
 a_2 &= -\Big((t_x-1) (1+t_z) (-2 (1+t_y) (-1+2 \Delta y+t_y) (t_z-1) t_z^2\\\nonumber 
 &-t_x (t_z-1) t_z (3-3 t_y^2+16 \Delta y (1+2 t_y) t_z+(t_y-1) (1+t_y) (5+32 t_y) t_z)\\\nonumber 
 &+8 t_x^3 (-3 \Delta z t_z^2+3 \Delta z t_y^2 t_z^2+2 (-2+\Delta y) t_y (t_z-1) (1+4 t_z) +4 t_y^3 (t_z-1) (1+4 t_z)) \\\nonumber 
 &+t_x^2 (-2 \Delta y (t_z-1) (1+t_y+4 (1+t_y) t_z+(3+35 t_y) t_z^2)\\\nonumber 
 &-(t_y^2-1) (-4+t_z (-9+t_z (10-24 \Delta z+128 t_y (t_z-1)+3 t_z))))) \Big)\\\nonumber & /\Big (2 \Delta x \Delta y (1+t_y) (t_z-1) (t_z^2+4 t_x t_z^2 +t_x^2 (2+t_z (8+3 t_z)))\Big )\\
 a_3 &= \Big ((t_x^2-1) (1+t_z) (-2 (1+t_y) (-1+\Delta y+t_y) (t_z-1) t_z^2 \\\nonumber 
 &-t_x (1+t_y) (t_z-1) t_z (3+(-5+8 \Delta y) t_z+t_y (-3+(-27+32 t_y) t_z))\\\nonumber 
 &+8 t_x^3 (-3 \Delta z t_z^2+3 \Delta z t_y^2 t_z^2+4 t_y^3 (t_z-1) (1+4 t_z)+2 t_y (t_z-1) (-2-8 t_z+3 \Delta y (1+t_z (4+t_z))))\\\nonumber &+t_x^2 (1+t_y) (2 \Delta y (t_z-1) (1+4 t_z)-(t_y-1) (-4+t_z (-9+t_z (10-24 \Delta z+128 t_y (t_z-1)+3 t_z)))))\Big )\\\nonumber  &
 /\Big (16 \Delta x \Delta y t_x t_y (t_z-1) (t_z^2+4 t_x t_z^2+t_x^2 (2+t_z (8+3 t_z)))\Big )\\
 a_4 &= \Big ((t_x-1) (1+t_x) (1+t_y (6+t_y)) (-2 (1+t_y) (-1+\Delta y+t_y) (t_z-1) t_z^2\\\nonumber 
 &-t_x (1+t_y) (t_z-1) t_z (3-5 t_z+8 \Delta y t_z+t_y (-3+5 t_z)) \\\nonumber 
 &+24 t_x^3 (-\Delta z t_z^2+\Delta z t_y^2 t_z^2+2 (-2+\Delta y) t_y (t_z-1) (1+t_z (4+t_z))\\\nonumber 
 &+4 t_y^3 (t_z-1) (1+t_z (4+t_z)))+t_x^2 (1+t_y) (2 \Delta y (t_z-1) (1+4 t_z)\\\nonumber 
 & -(t_y-1) (-4+t_z (-9+t_z (10-24 \Delta z+3 t_z)))))\Big)\\\nonumber 
 &/\Big(32 \Delta x \Delta z t_x (t_y-1) t_y (1+t_y) (t_z^2+4 t_x t_z^2+t_x^2 (2+t_z (8+3 t_z)))\Big )\\
 a_5 &= \Big ((1+t_x (6+t_x)) (-2 (1+t_y) (-1+\Delta y+t_y) (t_z-1) t_z^2 \\\nonumber 
 &-t_x (t_z-1) t_z (3-3 t_y^2+8 \Delta y (1+3 t_y) t_z+(t_y-1) (1+t_y) (5+32 t_y) t_z)\\\nonumber 
 &+8 t_x^3 (-3 \Delta z t_z^2+3 \Delta z t_y^2 t_z^2+2 (-2+\Delta y) t_y (t_z-1) (1+4 t_z)\\\nonumber 
 &+4 t_y^3 (t_z-1) (1+4 t_z))+t_x^2 (-2 \Delta y (t_z-1) (-1-4 t_z+t_y (-1+4 t_z) (1+8 t_z))\\\nonumber 
 &-(t_y^2-1) (-4+t_z (-9+t_z (10-24 \Delta z+128 t_y (t_z-1)+3 t_z)))))\Big )\\\nonumber &/\Big (32 \Delta y \Delta z t_x t_y (t_z^2+4 t_x t_z^2+t_x^2 (2+t_z (8+3 t_z)))\Big )
\end{align}}

\begin{landscape}
This gives the following mode: {\tiny
\begin{align}
 \sum_{r = 1}^5 a_r \hat Q_r = \Big ( & -\frac{\Delta x (1+4 t_x+t_x^2) ( \Delta y/\Delta x (1+t_y) (t_z^2-1)+1/\Delta x (t_y^2-1) (t_z^2-1)+8 t_x ( \Delta z/\Delta x (t_y^2-1) t_z+ \Delta y/\Delta x t_y (t_z^2-1)))}{36 \Delta y \Delta z t_x t_y t_z}, \\ \nonumber
 & \frac{(t_x-1) (1+t_x) (1+4 t_y+t_y^2) (8 /\Delta x t_x t_z+1/\Delta x/\Delta z (t_z^2-1))}{36 t_x t_y t_z},\frac{(t_x-1) (1+t_x) (1/\Delta x+1/\Delta x t_y+8 /\Delta x t_x t_y) (1+4 t_z+t_z^2)}{36 t_x t_y t_z},0,\\\nonumber 
 &-\frac{\Delta x (1+6 t_x+t_x^2) (\Delta y/\Delta x +1/\Delta x (t_y-1)) (t_z-1)}{4 \Delta y \Delta z t_x},4 1/\Delta x/\Delta z (t_x-1) (1+t_x) t_y (1+t_y) (t_z-1), \\\nonumber 
 &\frac{(t_x-1) (1+t_x) (1/\Delta x \Delta y+\Delta z (t_y-1) (1/\Delta x/\Delta z-8 /\Delta x/\Delta z t_x t_y)) (1+t_z)}{2 \Delta y t_x},0,\\\nonumber 
 &-\frac{2 \Delta x t_x (1+t_x) (\Delta y/\Delta x +2 /\Delta x (t_y^2-1)) (t_z-1)}{\Delta y \Delta z},\frac{1/\Delta x/\Delta z (t_x-1) (1+6 t_y+t_y^2) (t_z-1)}{4 t_y},\\\nonumber  &\frac{(t_x-1) (1/\Delta x (1- t_y^2)+4 t_x t_y (\Delta y/\Delta x +2 /\Delta x (t_y^2-1))) (1+t_z)}{2 \Delta y t_y},0,\\\nonumber 
 &\frac{\Delta x (1+t_x) (1/\Delta x (\Delta y-\Delta y t_z^2)+\Delta z (t_y-1) (1/\Delta x/\Delta z-4 /\Delta x t_x t_z-1/\Delta x/\Delta z t_z^2+8 1/\Delta x/\Delta z t_x t_y (t_z^2-1)))}{2 \Delta y \Delta z t_z},\\\nonumber 
 &-\frac{(t_x-1) (1+t_y) (1/\Delta x/\Delta z-4 /\Delta x t_x t_z-1/\Delta x/\Delta z t_z^2+8 /\Delta x/\Delta z t_x t_y (t_z^2-1))}{2 t_z},\frac{1/\Delta x (t_x-1) (1+6 t_z+t_z^2)}{4 t_z},0,\\\nonumber 
 &-\frac{\Delta x t_x (1+t_x) (\Delta z/\Delta x  (t_y^2-1) t_z+\Delta y/\Delta x  t_y (t_z^2-1))}{\Delta y \Delta z t_y t_z},\\\nonumber 
 &-\frac{(t_x-1) (1+6 t_y+t_y^2) (1/\Delta x/\Delta z-4 1/\Delta x t_x t_z-1/\Delta x/\Delta z t_z^2+8 /\Delta x/\Delta z t_x t_y (t_z^2-1))}{8 t_y t_z},\\\nonumber 
 &\frac{(t_x-1) (1/\Delta x (1- t_y^2)+4 t_x t_y (\Delta y/\Delta x +2 /\Delta x (t_y^2-1))) (1+6 t_z+t_z^2)}{8 \Delta y t_y t_z},0,\\\nonumber 
 &\frac{\Delta x (1+6 t_x+t_x^2) (1/\Delta x (\Delta y-\Delta y t_z^2)+\Delta z (t_y-1) (1/\Delta x/\Delta z-4 /\Delta x t_x t_z-1/\Delta x/\Delta z t_z^2+8 1/\Delta x/\Delta z t_x t_y (t_z^2-1)))}{8 \Delta y \Delta z t_x t_z},\frac{ (t_x-1) (1+t_x) (1+t_y)}{\Delta x},\\\nonumber 
 &\frac{(t_x-1) (1+t_x) (\Delta y/\Delta x +\Delta z (t_y-1) (1/\Delta x/\Delta z-8 /\Delta x/\Delta z t_x t_y)) (1+6 t_z+t_z^2)}{8 \Delta y t_x t_z},0,\\\nonumber 
 &-\frac{\Delta x (1+6 t_x+t_x^2) (\Delta y/\Delta x +2 /\Delta x (t_y^2-1)) (t_z-1)}{2 \Delta y \Delta z},\frac{ (t_x-1) (1+t_x) (1+6 t_y+t_y^2) (t_z-1)}{\Delta x \Delta z},\frac{ (t_x-1) (1+t_x) (1+t_z)}{\Delta x},0, \\\nonumber 
 &-\frac{\Delta x (1+t_x) (\Delta y/\Delta x +1/\Delta x (t_y-1)) (t_z-1)}{\Delta y \Delta z}, \frac{(t_x-1) (1+t_y) (t_z-1)}{\Delta x \Delta z},\frac{ (t_x-1) (1+t_z)}{\Delta x},0 \Big )
\end{align}
}
\end{landscape}

\end{document}